\documentclass[preprint,12pt]{elsarticle}

\usepackage[utf8]{inputenc} 
\usepackage[export]{adjustbox}
\usepackage[hang]{footmisc} 
\usepackage[english]{babel}
\usepackage[section]{algorithm} 
\usepackage{algpseudocode} 
\usepackage{lineno}   
\usepackage{amsmath}
\usepackage{amssymb}
\usepackage{amsfonts}
\usepackage{amsthm}
\usepackage{subfig}
\usepackage{tabularx}  
\usepackage{makecell}  
\usepackage{booktabs}  
\usepackage{multirow}  
\usepackage{adjustbox} 
\usepackage{color}
\usepackage{xcolor}
\usepackage{gensymb}
\usepackage{bm}
\usepackage{hyperref}
\hypersetup{
  colorlinks=false,
  frenchlinks=false,
  pdfborder={0 0 0},
  naturalnames=false,
  hypertexnames=false,
  breaklinks
}
\graphicspath{{Figures/}} 

\makeatletter
\let\OldStatex\Statex
\renewcommand{\Statex}[1][3]{%
  \setlength\@tempdima{\algorithmicindent}%
  \OldStatex\hskip\dimexpr#1\@tempdima\relax}
\makeatother

\newcommand*\patchAmsMathEnvironmentForLineno[1]{%
  \expandafter\let\csname old#1\expandafter\endcsname\csname #1\endcsname
  \expandafter\let\csname oldend#1\expandafter\endcsname\csname end#1\endcsname
  \renewenvironment{#1}%
  {\linenomath\csname old#1\endcsname}%
  {\csname oldend#1\endcsname\endlinenomath}}%
  \newcommand*\patchBothAmsMathEnvironmentsForLineno[1]{%
  \patchAmsMathEnvironmentForLineno{#1}%
  \patchAmsMathEnvironmentForLineno{#1*}}%
  \AtBeginDocument{%
  \patchBothAmsMathEnvironmentsForLineno{equation}%
  \patchBothAmsMathEnvironmentsForLineno{align}%
  \patchBothAmsMathEnvironmentsForLineno{flalign}%
  \patchBothAmsMathEnvironmentsForLineno{alignat}%
  \patchBothAmsMathEnvironmentsForLineno{gather}%
  \patchBothAmsMathEnvironmentsForLineno{multline}%
}

\newcommand{\ra}[1]{\renewcommand{\arraystretch}{#1}}
\newcommand{\J}{\mathcal{J}}
\newcommand{\br}{\mathbf{r}}
\newcommand{\bv}{\mathbf{v}}
\newcommand{\bw}{\mathbf{w}}
\newcommand{\bbv}{\overline{\bv}}

\newcommand{\gamg}{GAMG}
\newcommand{\amgname}{aSP-AMG}
\newcommand{\aspdef}{aSP-AMG (Def)}
\newcommand{\aspopt}{aSP-AMG (Opt)}
\renewcommand{\vec}[1]{\mathbf{#1}}
\providecommand{\vecT}[1]{\mathbf{#1}^{T}}
\providecommand{\mat}[1]{#1}
\providecommand{\matT}[1]{{#1}^{T}}
\providecommand{\matI}[1]{{#1}^{-1}}
\providecommand{\matIT}[1]{{#1}^{-T}}

\providecommand{\twonorm}[1]{\left|\left| #1 \right|\right|_{2}}

\journal{Comp. Methods in Appl. Mech. and Eng.}
\begin{document}
\begin{frontmatter}

\title{A robust adaptive algebraic multigrid linear solver for structural mechanics}
\author[label1]{Andrea Franceschini}
\author[label2]{Victor A. Paludetto Magri}
\author[label2]{Gianluca Mazzucco}
\author[label3]{Nicol\`{o} Spiezia}
\author[label3]{Carlo Janna\corref{cor1}}
\address[label1]{Department of Energy Resources, Stanford University, Stanford, CA, United States of America}
\address[label2]{Department ICEA, University of Padova, Via Marzolo, 9 - 35131 Padova, Italy}
\address[label3]{M$^3$E s.r.l., Via Giambellino, 7 - 35129 Padova, Italy}
\cortext[cor1]{Corresponding author: c.janna@m3eweb.it}

\begin{abstract}
The numerical simulation of structural mechanics applications via finite elements usually
requires the solution of large-size and ill-conditioned linear systems, especially when
accurate results are sought for derived variables interpolated with lower order functions,
like stress or deformation fields. Such task represents the most time-consuming kernel in
commercial simulators; thus, it is of significant interest the development of robust and
efficient linear solvers for such applications. In this context, direct solvers, which are
based on LU factorization techniques, are often used due to their robustness and easy
setup; however, they can reach only superlinear complexity, in the best case, thus, have
limited applicability depending on the problem size. On the other hand, iterative solvers
based on algebraic multigrid (AMG) preconditioners can reach up to linear complexity for
sufficiently regular problems but do not always converge and require more knowledge from
the user for an efficient setup. In this work, we present an adaptive AMG method
specifically designed to improve its usability and efficiency in the solution of
structural problems. We show numerical results for several practical applications with
millions of unknowns and compare our method with two state-of-the-art linear solvers
proving its efficiency and robustness.
\end{abstract}

\begin{keyword}
AMG \sep linear solver \sep preconditioning \sep ill-conditioning \sep structural problems
\sep parallel computing
\end{keyword}

\end{frontmatter}

\section{Introduction}
\label{sec:Intro}

Across a broad range of structural mechanics applications, the demand for more accurate,
complex and reliable numerical simulations is increasing exponentially. In this context,
the Finite Element Method (FEM) remains the most widely used approach and the number of
new extremely challenging applications is countless, e.g., modeling of fractal formation
in macroscopic elasto-plasticity in three-dimensional bodies \cite{LiSahKorOst12},
simulation of cardiac mechanics \cite{PavScaZam15}, submarine landslides
\cite{ZhaOnaAnd19}, stir welding processes \cite{DiaChiCer17}, concrete gravity dam
\cite{WanWanLu17}, just to name some recent applications. In all the problems mentioned
above, the underlying partial differential equations (PDEs) are discretized to approximate
the continuous solution in an algebraic system of equations of the form:
\begin{equation}
  \mat{A} \vec{x} = \vec{b},
\label{eq:system}
\end{equation}
with $\vec{b}$ and $\vec{x}$, the right hand side and solution vectors, respectively,
belonging to $\mathbb{R}^{n}$ and $\mat{A} \in \mathbb{R}^{n \times n}$ the matrix
deriving from the Finite Element discretization. For small strain mechanical problems, the
linear system in \eqref{eq:system} assumes the notation $\mat{K} \vec{u} = \vec{f}$, where
the vectors $\vec{u}$ and $\vec{f}$ represent the unknown nodal displacements and the
applied nodal forces, respectively, and the matrix $\mat{K}$ arises from the numerical
integration of $\int_V \matT{B} \mat{D} \mat{B} dV $, where $\mat{D}$ stands for the
constitutive matrix and $\mat{B} = \partial{\mat{N}}/ \partial{\mat{x}}$ contains the
spatial derivatives of the element shape functions $\mat{N}$. If the material is linear
elastic, as assumed in this work, the matrix $\mat{D}$ is Symmetric Positive Definite
(SPD) and it is fully defined once the Young modulus $E$ and the Poisson ratio $\nu$ are
assigned to each element composing the grid.

In numerical simulations, the solution of the linear system in \eqref{eq:system} is the
most time-consuming part of the entire process, taking up $70\%-80\%$ of the total
computational time \cite{KorGup16}. With the constant demand for larger simulations,
involving up to billions of unknowns, the linear equation solver is the most significant
bottleneck, and an inefficient algorithm may dramatically slow down the simulations
process, forcing either to simplify the problem or to wait for a very long time for the
results.

Linear equation solvers may be grouped into two broad categories: direct methods and
iterative methods. The main advantages of direct methods are their generality and
robustness. However, the matrix factors created to solve the system are often
significantly denser than the original matrix. This issue leads to memory shortage for
both forming and storing the factors, in particular for large-scale system arising from
the discretization of three-dimensional problems. On the other hand, iterative methods are
much more suitable for the solution of large systems of equations and hence are gaining
more and more attention as far as the dimension of the computational grid
increases. However, iterative methods do not guarantee the convergence to the scheme,
unless a proper preconditioner is adopted. In fact, preconditioning the system is
essential not only to accelerate the convergence rate, but also to avoid divergence in
ill-conditioned cases.

In linear elastic problems, as those investigated in this work, ill-conditioned systems
may arise due to several conditions, such as the presence of multiple materials or highly
distorted elements, two of the most common source of ill-conditioness in structural
problems.  The first one is when the model is characterized by materials with
significantly different constitutive characteristics and therefore the entries of the
linear system matrix present large jumps.  This is the case for example when a soft and
hard tissues are connected \cite{CotPeaYou16,HasAdaTak16,TruMas17} or when a solid matrix
presents soft inclusions \cite{ZhaLiWan15,ZhoSonLu17,MazAntPel18,PalKup18}. The second
case appears when the model is characterized by a highly distorted mesh, with poorly
proportioned elements. This is the case when complex geometries characterize the model,
which cannot be discretized with a structured mesh. This undesired situation is becoming
more frequent since it is increasingly common to create meshes directly from complex CAD
solid geometries \cite{BazBanBra15,SchKluBar16} or from an X-Ray tomography
\cite{YouMaiGae05,HuaYanRen15,SenYanWan16,WuXiaWen16}. Moreover, newly developed
technologies, e.g., additive manufacturing, allow to create elaborated solids, with
consequent demanding numerical simulations \cite{BikStaChr16,CheGuiGan17,GouMicDen18}.
Besides, also the presence of materials close to the incompressibility limit or a loosely
constrained body can lead to an ill-conditioned system. In the last decades, several types
of preconditioners have been developed, such as incomplete factorizations
\cite{Saa94,LinMor99,Ben02}, sparse approximate inverses
\cite{BenMeyTum96,Tan99,Huc03,JanFerSarGam15}, domain decomposition
\cite{DolJolNat15,Zam16,BadMarPri16,LiSaa17} and Algebraic Multigrid (AMG) methods. This
paper focuses on the last category.

AMG methods are built on a hierarchy of levels associated with linear problems of
decreasing size. Such methods are defined by the choice of interpolation operators, which
transfers information between different levels; the coarsening strategy, which guides the
definition of new levels; the smoothing technique, which solves high frequencies
components of the error on the given level and, lastly, the application strategy, which
defines how the multigrid cycle is applied.

Several families of AMG methods can be found in the literature. The first multigrid
strategy, namely the classical AMG \cite{Stu83,Bra86}, was proposed in the early 1980s for
the efficient solution of M-matrices and Poisson models. These methods are built on the
knowledge that the near-kernel of the operator is well approximated by the constant vector,
which is a limiting hypothesis in linear elasticity problems. In the latter, a larger
near-kernel, usually well represented by the rigid body modes (RBMs), is needed to obtain
good results. The first attempt to overcome such limitation came in the early 1990s with
the Smoothed Aggregation AMG (SA-AMG) method, where coarsening is done via aggregation
of nodes and interpolation is built column-wise starting from a tentative operator
spanning the RBMs in the aggregates \cite{Van92,Van96}.
The element-based AMG family is composed by the energy-minimization AMGe \cite{BreCleFal01},
element-free AMGe \cite{HenVas01} and spectral AMGe \cite{ChaFalHenJonManMccRugVas03}.
Here, the coarse spaces are constructed via an energy minimization process with the aim of
improving robustness by alleviating the heuristics based on M-matrices properties
implemented in classical AMG. More recently, the adaptive and Bootstrap AMG ($\alpha$AMG
and BAMG, respectively) were designed for the solution of more difficult problems where
the classical and smoothed AMG may fail or show poor convergence
\cite{BreFalMacManMccRug05,BraBraKahLiv11,BraBraKahLiv15,DamFilVas18}. In such methods, no
preliminary assumption is made about the near-null space of $\mat{A}$, but it is
approximated adaptively during the AMG hierarchy construction. For a comprehensive review
of algebraic multigrid variants, we refer the reader to \citet{XuZik17}.

In the context of elasticity problems, \citet{Yoo03} presented a W-cycle method for
solving linear elasticity problems in the nearly incompressible limit and \citet{GriOelSch03}
presented a generalization of the classical AMG where a block interpolation method is proposed
and showed to reproduce the RBMs components in a multilevel strategy.
\citet{BakKolYan10} investigated several approaches for improving
convergence of classical AMG when solving linear elasticity problems by incorporating the
RBMs in the range of interpolation.

This work presents an extension of the adaptive Smoothing and Prolongation based Algebraic
Multigrid method (\amgname{}), proposed by \cite{MagFraJan19}, with the aim of
specifically improving its performance for the solution of linear elasticity
problems. Such method follows the path of bootstrap and adaptive AMG which is to build an
approximation of the near-null space components of the problem at hand
automatically. Also, this method proposes new strategies for computing the interpolation
operators in a least-squares sense and introduces the Adaptive Factorized Sparse
Approximate Inverse (aFSAI) as a flexible smoother, with its accuracy automatically tuned
during set-up. This is an essential advantage with respect to the commonly used weighted
Jacobi and Gauss-Seidel methods.  Moreover, an approach to improve the test space
computation by considering a priori knowledge of the RBMs in case of elasticity is
proposed and, finally, the smoother calculation is enhanced by introducing a heuristic
towards an optimal balance between quality and setup cost. To demonstrate the potentiality
of the method and its broad applicability, the developed linear solver has been applied in
the solution of linear systems arising from the discretization of a comprehensive range of
real-world structural problems, spanning from biomechanics, reservoir geomechanics,
material mechanics, solid mechanics and so on. These problems have been chosen not only
for their large size (millions of unknowns) but also because each of them presents one or
more source of ill-conditioness. The results prove the efficiency and robustness of the
method and show how in several cases the proposed algorithm outperforms state-of-the-art
AMG linear solvers such as BoomerAMG in the HYPRE package \cite{FalYan02} and GAMG in the
PETSc library \cite{PETSc}. Even more important, the results show how the proposed linear
solver gives good results even assuming a default set of parameters, making it fully
adoptable by a user not keen in AMG parameters fine tuning.

The paper is organized as follows. Section \ref{sec:AMG} presents the fundamental aspects
of the AMG algorithm. Section \ref{sec:Improving} describes the numerical methods
developed to allow AMG to tackle structural problems. Section \ref{sec:numexp} presents
the computational performance of the enhanced \amgname{}, tested on several real-world
challenging examples. Finally, Section \ref{sec:conclusions} draws some conclusions.

\section{Algebraic Multigrid}
\label{sec:AMG}

Algebraic multigrid refers to an efficient class of algorithms for solving linear systems
of equations. These methods, which can be used as solvers or preconditioners for
Krylov-subspace iterative schemes, are characterized by tackling, in a multilevel fashion,
different components of the error $\vec{e}^{k} = \matI{A} \vec{b} - \vec{x}^{k}$
corresponding to the $k$-th solution guess $\vec{x}^{k}$. Particularly, correction vectors
are computed via the approximate solution of linear systems of reduced size derived
algebraically from $\mat{A}$ and combined to form a global correction yielding an updated
solution $\vec{x}^{k+1}$. The main feature of AMG is that it can efficiently reduce
error components on both the low and high part of the eigenspectrum of $\mat{A}$.  If
designed correctly, AMG leads to a convergence rate which is independent of the mesh size
giving optimal execution times and enabling the solution of linear problems in the order
of billions of unknowns.

For the sake of simplicity, let us assume a two-level method. AMG starts by applying a
fixed-point iteration scheme, such as (damped) Jacobi or (block) Gauss-Seidel, to the
linear system \eqref{eq:system}. This is written mathematically by the following equation:
\begin{equation}
\vec{x}^{k + 1/2} = \left( \mat{I} - \matI{M} \mat{A} \right) \vec{x}^{k} + \matI{M}
\vec{b},
\label{eq:preSmootherApp}
\end{equation}
where $\matI{M}$ denotes the preconditioning operator that depends on the choice of the
scheme, e.g., the diagonal of $\mat{A}$ in case of Jacobi. Such process, called relaxation
or smoothing, stops after a few iterations because convergence degrades very quickly. This
is a well-known issue of fixed-point iteration schemes which are very effective on the
high-frequency error components but are unable to handle low frequencies. To reduce the
latter, AMG proceeds with the idea that they can be viewed as high-frequencies in a lower
resolution representation of $\mat{A}$, i.e., in a coarser projection $\mat{A}_c =
\matT{P} \mat{A} \mat{P}$ where $\mat{A}_c \in \mathbb{R}^{n_c \times n_c}$ is the
coarse-grid matrix and $\mat{P} \in \mathbb{R}^{n \times n_c}$ is a prolongation operator
that transfers information from a low-resolution (coarse) to a high-resolution (fine)
level. Following this idea, the smoothed residual $\vec{r}^{k + 1/2} = \vec{b} - \mat{A}
\vec{x}^{k + 1/2}$ is projected onto the coarse space where a correction vector is
computed using $\mat{A}_c$:
\begin{equation}
\vec{h}_c = \matI{A}_c \matT{P} \vec{r}^{k + 1/2},
\end{equation}
and interpolated back to the fine space through $\mat{P}$ to obtain the correction
$\vec{h} = \mat{P} \vec{h}_c $. Usually, the final approximation is obtained with another
relaxation step on the corrected solution:
\begin{equation}
\vec{x}^{k + 1} = \left( \mat{I} - \matIT{M} \mat{A} \right) \left( \vec{x}^{k + 1/2} +
\vec{h} \right) + \matI{M} \vec{b}.
\label{eq:postSmootherApp}
\end{equation}
This last step has the additional advantage of maintaining symmetry even when $\mat{M}$ is
not a symmetric operator, as required when using AMG for preconditioning in the conjugate
gradient method.

In the two-grid scheme, the combined action of relaxation and coarse-grid correction is
defined by the error propagation operator:
\begin{equation}
E_{TG} = \left( \mat{I} - \matIT{M} \mat{A} \right) \left( \mat{I} - \mat{P} \matI{A}_c
\mat{R} \mat{A} \right) \left( \mat{I} - \matI{M} \mat{A} \right).
\label{eq:twoGridErrorOp}
\end{equation}
An efficient AMG solver regarding convergence has the norm of $E_{TG}$ close to
zero. Such condition is met when smoothing and coarse-grid correction are complementary to
each other, i.e., error components not reduced by relaxation, also called algebraically
smooth, are properly reduced by coarse-grid correction. On the other hand, this depends on
the ability of the prolongation operator to span those components problematic for
relaxation \cite{XuZik17}.

Typically, the coarsening factor, i.e., the relation between $n_c$ and $n$, is in the order
of $10^{-1}$. Thus the application of $\matI{A}_c$ via a direct method remains costly. To
make the AMG solver efficient in terms of computational time, the action of $\matI{A}_c$
is approximated by a new two-grid method. Applying this idea recursively until a maximum
number of levels or a minimum coarse-grid size is reached, we obtain the
V$(\nu_1,\nu_2)$-cycle AMG in its multilevel format, where $\nu_1$ and $\nu_2$ are the
number of pre- and post-smoothing steps adopted.

Algorithms \ref{algo:cptAMG} and \ref{algo:applyAMGv} schematically provide the main steps
needed to compute and apply in a V$(\nu_1,\nu_2)$-cycle, respectively, a classical AMG
scheme. The ideas used in our approach to constructing the smoother $\matI{M}$, the
prolongation operator $\mat{P}$ and coarsening will be detailed in the next section.

\begin{algorithm}[t!]
\caption{Recursive AMG Setup}
\begin{algorithmic}[1]
\Procedure{AMG\_SetUp}{$\mat{A}_k$}
\State Define $\Omega_k$ as the set of the $n_k$ vertices of the adjacency graph of
  $\mat{A}_k$;
\If {$n_k$ is small enough to allow a direct factorization}
  \State Compute $\mat{A}_k = \mat{L}_k \mat{L}_k^T$;
\Else
  \State Compute $\mat{M}_k$ such that $\mat{M}_k^{-1} \simeq \matI{A}_k$;
  \State Define the smoother as $\mat{S}_k = \left(\mat{I}_k - \mat{M}_k^{-1} \mat{A}_k
    \right)$;
  \State Partition $\Omega_k$ into the two disjoint sets $\mathcal{C}_k$ and
    $\mathcal{F}_k$ via coarsening;
  \State Compute the prolongation matrix $\mat{P}_k$ from $\mathcal{C}_k$ to
    $\Omega_k$;
  \State Compute the new coarse level matrix $\mat{A}_{k + 1} = \matT{P}_k \mat{A}_k
    \mat{P}_k$;
  \State \text{AMG\_SetUp}$\left(\mat{A}_{k+1}\right)$;
\EndIf
\EndProcedure
\end{algorithmic}
\label{algo:cptAMG}
\end{algorithm}
\begin{algorithm}[t!]
\caption{AMG application in a V$(\nu_1,\nu_2)$-cycle}
\begin{algorithmic}[1]
\Procedure{AMG\_Apply}{$\mat{A}_k$, $\vec{y}_k$, $\vec{z}_k$}
\If{$k$ is the last level}
  \State Forward and Backward solve $\mat{L}_k \matT{L}_k \vec{z}_k = \vec{y}_k$;
\Else
  \State Smooth $\nu_1$ times $\mat{A}_k \vec{s}_k = \vec{y}_k$ starting from $\vec{0}$;
  \State Compute the residual $\vec{r}_k = \vec{y}_k - \mat{A}_k \vec{s}_k$;
  \State Restrict the residual to the coarse grid $\vec{r}_{k+1} = \matT{P}_k \vec{r}_k$;
  \State \text{AMG\_Apply}$\left(\mat{A}_{k+1},\vec{r}_{k + 1},\vec{h}_{k+1} \right)$;
  \State Prolongate the correction to the fine grid $\vec{h}_{k} = \mat{P}_k
    \vec{h}_{k+1}$;
  \State Update $\vec{s}_{k} \leftarrow \vec{s}_{k} + \vec{h}_{k}$;
  \State Smooth $\nu_2$ times $\mat{A}_k \vec{z}_k = \vec{y}_k$ starting from $\vec{s}_k$;
\EndIf
\EndProcedure
\end{algorithmic}
\label{algo:applyAMGv}
\end{algorithm}

\section{Numerical methods to improve AMG in structural problems}
\label{sec:Improving}

This section aims to illustrate the main components of \amgname{}, their behavior with
respect to the user-defined parameters and the connection with structural
problems. Moreover, we investigate the relationship among the smoother quality and the
overall AMG cycle effectiveness.

\subsection{Adaptive smoother computation}
\label{subsec:smoother}

One of the critical ingredients for a successful AMG method is the availability of an
effective and scalable smoother. This is particularly true for structural problems where
damping highest frequencies often requires the use of weights or Chebyshev polynomials
\cite{AdaBreHuTum03,BakFalKolYan11}. In \cite{MagFraJan19}, the aFSAI
\cite{JanFerSarGam15} is proposed as smoother and its effectiveness is assessed on an
extensive set of numerical experiments. aFSAI is designed for SPD matrices and, as
smoother, takes the following form:
\begin{equation}
\mat{S} = \mat{I} - \omega \matT{G} \mat{G} \mat{A},
\end{equation}
with $\matT{G} \mat{G}$ an approximation of $\matI{A}$ in factored form. The damping
parameter $\omega$ is used to ensure $\rho(\mat{S}) < 1$, being $\rho(\mat{S})$ the
spectral radius of $\mat{S}$. From a practical point of view, $\omega$ is chosen in such a
way that $\omega < \mbox{min}(1,{2}/{\lambda_1(\mat{G} \mat{A} \matT{G})})$, where the
maximum eigenvalue of $\mat{G} \mat{A} \matT{G}$, $\lambda_1$, is estimated with a few
Lanczos passes. The aFSAI setup exhibits a very high degree of parallelism with each
$\mat{G}$ row computed independently from the others. The number of entries of each
$\mat{G}$ row is adaptively increased in an iterative process until either a desired
accuracy or a maximal density is reached. One of the major downsides of this technique,
however, is that its setup may be expensive and its cost strongly depends on a few
user-specified parameters \cite{JanFerSarGam15}. In fact, the user controls the smoother
accuracy by specifying the number of steps of the iterative process, $k$, the number of
entries added per row at each step, $\rho$, and an exit tolerance, $\epsilon$, which is
based on an estimate of the Kaporin condition number of $\mat{G} \mat{A} \matT{G}$, e.g.,
see \cite{JanFerSarGam15}. Unfortunately, the aFSAI setup cost grows up very quickly with
its density and the user must pay careful attention in tuning $k$, $\rho$, and $\epsilon$
to avoid wasting resources unnecessarily. Due to the high number of parameters already
involved in the AMG setup, adding this extra effort in tuning the smoother may discourage
potential users. For this reason, we present a practical strategy to compute aFSAI able to
automatically tune the smoother's quality on the problem at hand, leaving all the setup
parameters to a default value. This strategy is mainly based on the observation that it is
possible to improve a previously computed $\mat{G}$ without losing the effort already
spent on it.

Let us call $\mat{G}_0$ the factor computed from scratch with parameters $k_0$, $\rho_0$
and $\epsilon_0$. We can improve its effectiveness by performing other $k_i$ adaptive
steps with size $\rho_i$ and exit tolerance $\epsilon_i$. If the smoother's quality is not
high enough, we can refine it by performing other steps of the adaptive procedure. The
main problem is how to estimate the smoother's quality inexpensively. Our experimentation
over a wide range of structural problems has shown that, when the $\omega$ value is close
to one, then there is no need to improve the smoother while, when $\omega$ is low,
improving the smoother may significantly decrease the overall AMG iterations. A complete
discussion on the relation among $\omega$ and the AMG effectiveness is presented in
Section \ref{sec:omega_AMG}. We can estimate $\omega$ with few Lanczos steps and use this
estimated value to stop the iteration whenever it falls below a prescribed
$\overline{\omega}$ value. Unfortunately, the aFSAI setup cost overgrows with its density,
so we also specify the parameter $\overline{\rho}$ representing the maximum allowed number
of average $G$ entries per row. The above procedure is summarized in Algorithm
\ref{algo:aFSAI}, where $\mbox{aFSAI}(A,G_s,k,\rho,\epsilon)$ is the function returning
the aFSAI factorization of $A$ obtained from the starting factor $G_s$ using the
parameters set $(k,\rho,\epsilon)$ and $\mbox{nnz}_r(\cdot)$ is the function returning the
average number of entries per row of the argument matrix. Note that, even though Algorithm
\ref{algo:aFSAI} formally needs a large number of input parameters, we show in the
numerical experiments that all of them can be safely set to a default value.

\begin{algorithm}[t!]
\caption{\bf aFSAI Smoother Set-up}
\begin{algorithmic}[1]
\Procedure{Smoother\_SetUp}{$A$,$k_0$,$\rho_0$,$\epsilon_0$,$k_i$,$\rho_i$,$\epsilon_i$,$\overline{\omega}$,$\overline{\rho}$}
\State Compute $\mat{G} \leftarrow \mbox{aFSAI}(\mat{A},\mat{I},k_0,\rho_0,\epsilon_0)$;
\State Compute $\omega \leftarrow \dfrac{2}{\lambda_1(\mat{G}\mat{A}\matT{G})}$;
\While {$\omega < \overline{\omega}$ {\bf and} $\mbox{nnz}_r(\mat{G}) < \overline{\rho}$}
  \State Compute $\mat{G} \leftarrow \mbox{aFSAI}(\mat{A},\mat{G},k_i,\rho_i,\epsilon_i)$;
  \State Compute $\omega \leftarrow \dfrac{2}{\lambda_1(GAG^T)}$;
\EndWhile
\EndProcedure
\end{algorithmic}
\label{algo:aFSAI}
\end{algorithm}

\subsection{Influence of $\omega$ on the AMG cycle}
\label{sec:omega_AMG}

Having a convergent smoother, that is $\mat{S}$ such that $\rho(\mat{S}) \leq 1$, is an
essential condition for the overall AMG-cycle effectiveness. In fact, it can be shown
that, in situations where $\rho(\mat{S}) > 1$, the preconditioned matrix may have negative
eigenvalues thus preventing PCG to converge. For the above reason, it is common to use a
damping factor, $\omega$, on those smoothers that do not guarantee $\rho(\mat{S}) \leq
1$. However, it can be observed from numerical experiments that, when the smoother needs a
too small relaxation factor, the overall AMG cycle quality decreases significantly. This
behavior can be explained as follow: the smoother, by its nature, is constructed to damp
the highest frequencies of $\mat{A}$ while the coarse grid correction selectively acts on
those eigenvectors associated with low frequencies, by eliminating the error components
belonging to the subspace of $\mathbb{R}^n$ spanned by $\mat{P}$. By damping $\mat{S}$
with a small $\omega$ value we do not alter the range of $\mat{P}$, but we reduce the
action of $\mat{S}$ on those eigenvectors that are not in $\mbox{range}(\mat{P})$. In
other words, when $\omega$ is too small a large part of the spectrum is pushed in a region
where neither the smoother nor the coarse grid correction are able to act effectively.

This behavior will be illustrated through a numerical example of a small matrix arising
from the solution of the equilibrium equation over a cube discretized with P1 linear
tetrahedra. This matrix has $1,773$ rows and $63,927$ non-zeroes and the corresponding
smoother is aFSAI computed with parameters $k=2$, $\rho=1$ and $\epsilon=10^{-4}$. With
this setup, the smoother is not convergent, i.e., $\rho(\mat{S}) > 1$, thus the use of
$\omega \leq 0.922$ is needed to ensure the overall convergence. With these settings, the
pure AMG solver without CG acceleration method needs $129$ iterations to solve a unitary
right-hand-side reducing the residual of $8$ orders of magnitude starting from the null
vector as the initial guess.

\begin{figure}[!htbp]
  \centering 
  \includegraphics[width=0.45\linewidth]{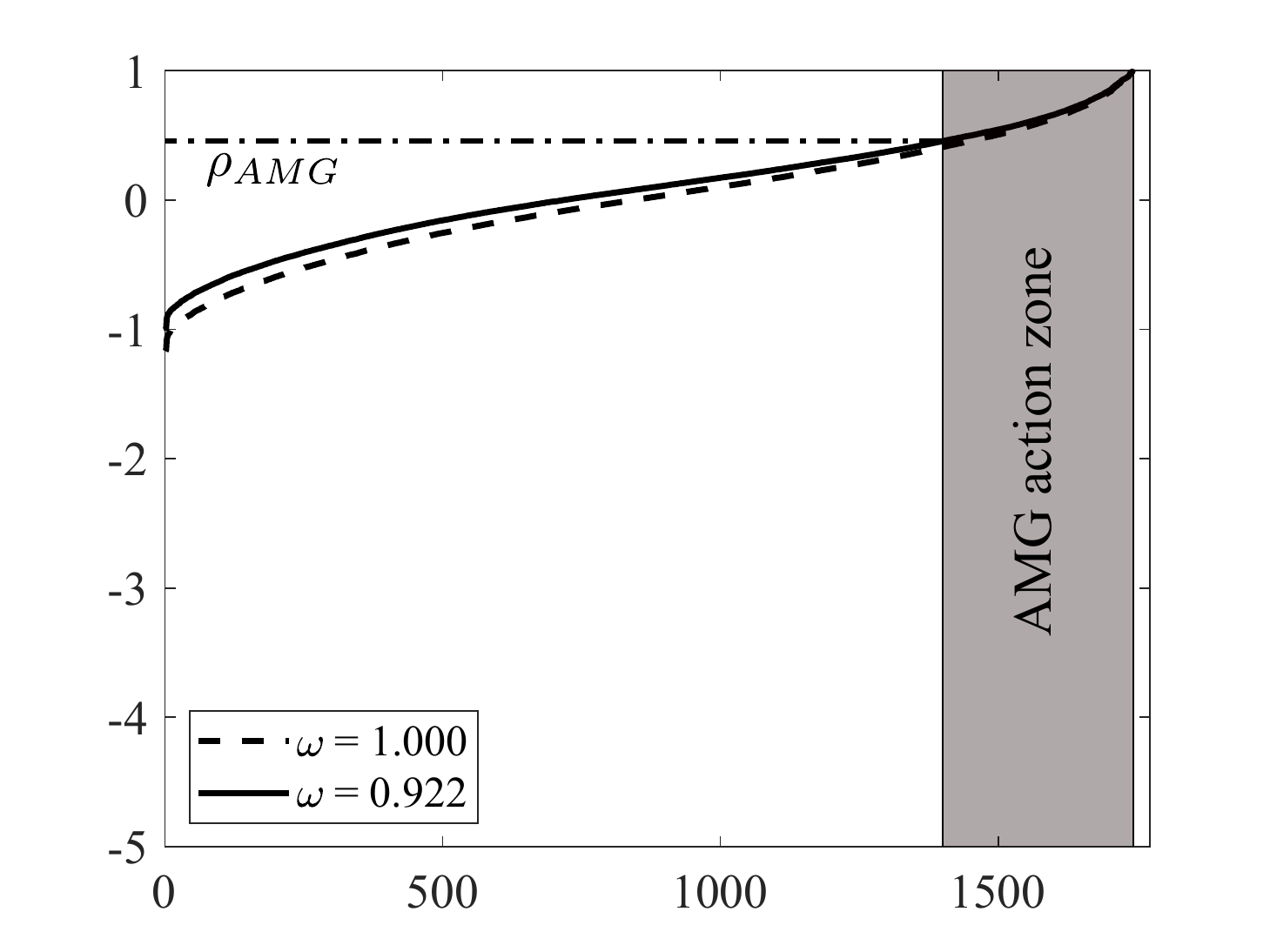}
  \includegraphics[width=0.45\linewidth]{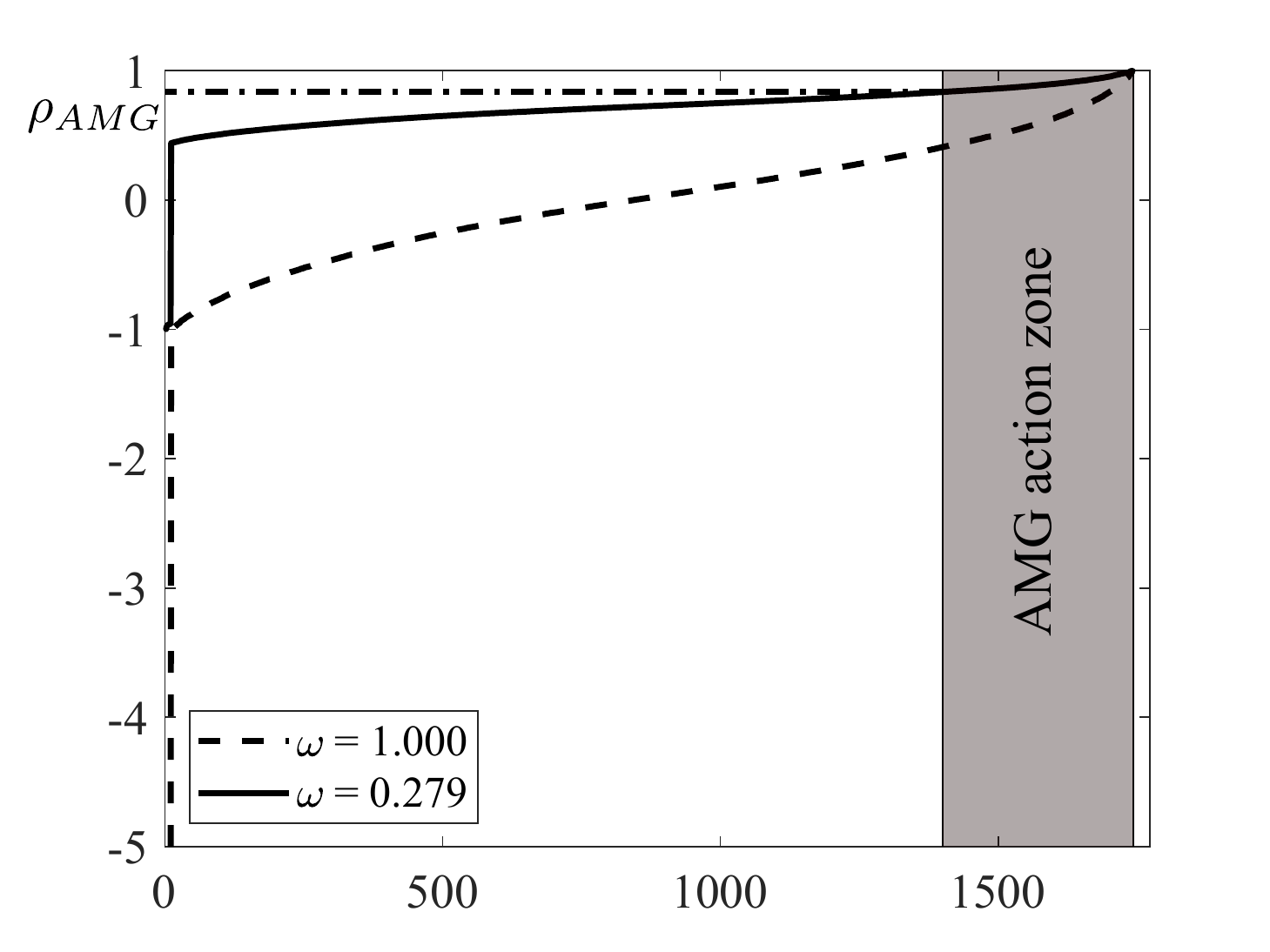}
  \caption{Eigenspecta for the two smoothers, $\mat{S} = \mat{I} - \omega \matT{G} \mat{G}
    \mat{A}$ and $\mat{S}_\alpha$, before (on the left) and after (on the right) the
    low-rank update. The darkened background highlights the part of the spectrum which is
    handled by the AMG hierarchy.}
  \label{fig:eig_SMO_upd}
\end{figure}

To show the impact of a small damping parameter we modify the above smoother $\mat{S} =
\mat{I} - \matT{G} \mat{G} \mat{A}$, leaving unaltered all the other AMG components, i.e.,
test space, prolongation and coarse grid correction. To this aim, we form the aFSAI
preconditioned matrix $\mat{G} \mat{A} \matT{G}$ and compute its $10$ largest eigenpairs
collecting the corresponding eigenvectors in a skinny matrix $\mat{U}$. This allow us to
create a low-rank update matrix $\mat{C} = \alpha \mat{U} \matT{U}$, with $\alpha$ a
multiplicative factor, such that the new preconditioned matrix has the same spectrum as
$\mat{G} \mat{A} \matT{G} + \mat{C}$ but same eigenvectors as $\mat{G} \mat{A}
\matT{G}$. This ensures that all the AMG hierarchy is the same as before. After some
algebra, the new smoother takes the following expression:
\begin{equation}
  \mat{S}_{\alpha} = \mat{I} - \left(\matT{G} \mat{G} + 
                   \alpha \matT{G} \mat{U} \matT{U} \matIT{G} \matI{A} \right) 
                   \mat{A},
\end{equation}
and the maximum damping factor allowed for $\mat{S}_{\alpha}$ is given by:
\begin{equation}
  \omega_\alpha = \frac{2}{\lambda_{1}(\mat{I} - \mat{S}_{\alpha})} 
               = \frac{2 \omega}{2 + \alpha \omega},
\end{equation}
where $\omega$ is the maximum damping factor allowed for the original smoother. In this
case, since $\omega = 0.922$, setting $\alpha = 5$ gives $\omega_\alpha = 0.279$ and the
the updated smoother with the same AMG hierarchy requires $276$ iterations instead of
$129$ to converge. Figure \ref{fig:eig_SMO_upd} shows the eigenspectra of the smoothers
before and after the low-rank update. It can be easily noted the effect of a lower
$\omega$ value causing a dramatic increase in the spectral radius of the AMG V-cycle.

\begin{figure}[!htbp]
  \centering 
  \includegraphics[width=0.45\linewidth]{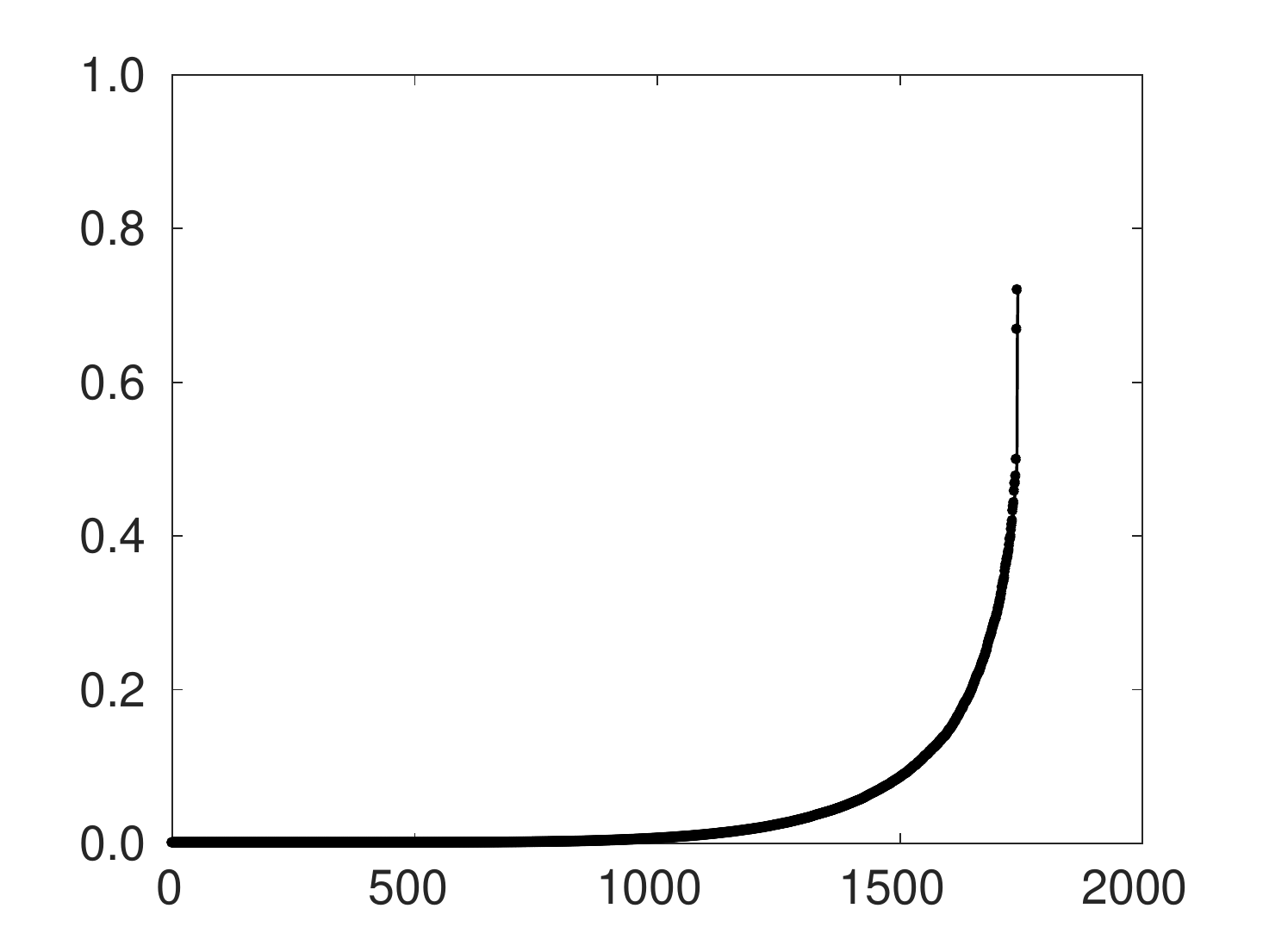}
  \includegraphics[width=0.45\linewidth]{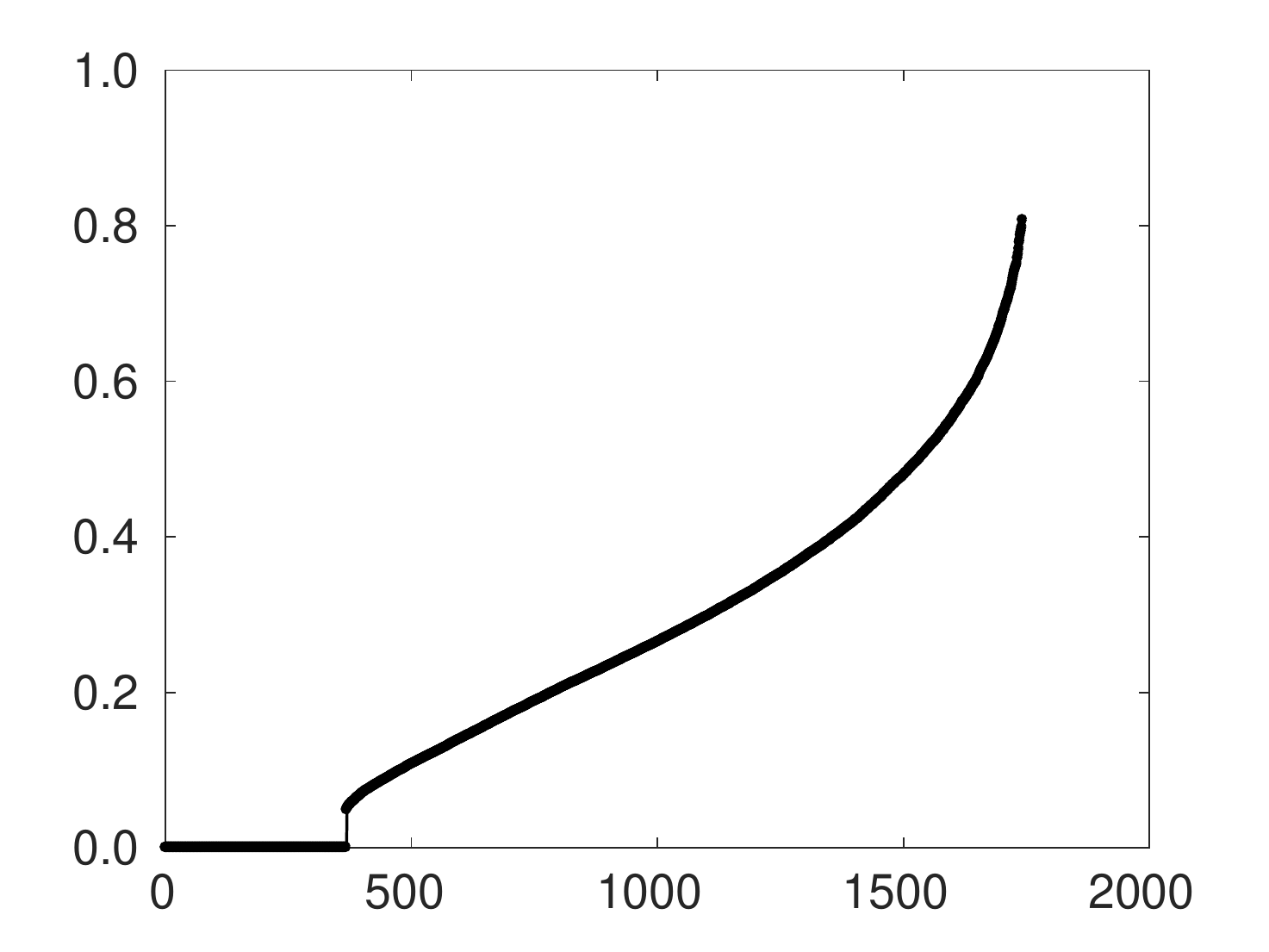}
  \caption{Eigenspecta for the two AMG iteration matrix, before (on the left) and after
    (on the right) the low-rank update.}
  \label{fig:eig_AMG_upd}
\end{figure}

The eigenvalue distribution of the two error propagation matrices for AMG relative to
$\mat{S}$ and $\mat{S}_\alpha$, as defined by Eq. \eqref{eq:twoGridErrorOp}, are plotted
in Figure \ref{fig:eig_AMG_upd}. The AMG iteration matrix, measuring the quality of the
V-cycle, has more clustered eigenvalues around zero with the larger $\omega$. Moreover,
using the original aFSAI smoother $\mat{G}$ the spectral radius of the whole AMG cycle is
about $0.50$, except for a couple of outliers, while for $\mat{S}_\alpha$ the distribution
is more uniform and the spectral radius is about $0.81$.

\subsection{Generation of the test space}
\label{subsec:tspace}

The central idea of adaptive AMG is that of uncovering algebraically smooth modes of the
system matrix by applying the smoother on a set of random vectors or by using some
eigensolution method. Smooth modes are by definition a set of vectors that are not damped
by relaxation and can be approximated by solving, usually with low accuracy, the
generalized eigenvalue problem:
\begin{equation}
\mat{A} \vec{v} = \lambda \mat{M} \vec{v},
\end{equation}
where $\mat{M}$ is the preconditioner, that in our case satisfies is in factored form
$\matI{M} = \matT{G} \mat{G}$. A basis for the test space, representing an approximation
of the near-kernel of $\mat{A}$, is collected in the $n \times n_{tv}$ matrix $\mat{V}$
and then used to generate the strength of connections (SoC) and to set up the prolongation
operator. Even if it is not required a high accuracy for $\mat{V}$, its computation may be
quite expensive and represents most of an adaptive AMG setup cost. For this reason, if an
approximation $\mat{V}_0$ of (or a part of) the near-kernel of $\mat{A}$ is available
beforehand, it is natural to use this information as an initial guess for building a
possibly larger test space $\mat{V}$. The knowledge of the physical problem originating
$\mat{A}$ usually gives a clue on the initial guess $\mat{V}_0$. When solving elasticity
problems, it is known from theory that rigid body modes (RBMs) constitute the kernel of
the differential operator before the application of boundary conditions; thus, they can be
used to define $\mat{V}_0$. As a matter of fact, such information is essential to the
effectiveness of several methods like aggregation-based AMG \cite{VanBreMan01,Ada02},
domain decomposition techniques \cite{FarRou91,XuZou98,FarPieLed00,DolJolNat15} and
deflation-based preconditioning \cite{FraVui01,JonGijVuiSca13,BagFraSpiJan17}, just to
cite a few. In 2D and 3D elasticity problems, there are 3 and 6 RBMs,
respectively. However, if a larger test space is needed for $\mat{A}$, the size of
$\mat{V}_0$ can be properly increased by padding it with random vectors.

Previous works on adaptive AMG suggest the use of power iteration or Lanczos algorithm for
the test space computation. However, the former generally provides slow convergence, while
the latter only accepts one vector as an initial guess. Given these limitations of the
power iteration or plain Lanczos, we propose a method based on the Rayleigh quotient
minimization. This method is called Simultaneous Rayleigh Quotient Minimization by
Conjugate Gradients (SRQCG) and was introduced by \cite{LonMcc80}. As the name suggests,
it simultaneously minimizes Rayleigh quotients over several linearly independent
vectors. A conjugate gradient technique is used to form each new vector through the linear
combination of current iterates and search directions. Moreover, as in other
subspace-based methods, SRQCG uses Ritz projections to accelerate convergence. Another
feature of SRQCG is that its initial convergence is generally faster than that of other
eigensolvers, e.g., Lanczos or Jacobi-Davidson, making it particularly attractive in
problems where accuracy is not a central issue \cite{BerMarPin06}.

The SRQCG method is described by Algorithm \ref{alg:SRQCG}. For a practical investigation
of this method through numerical experiments based on real-world applications, we refer
the reader to \cite{BerMarPin06,BerMarPin12,FerJanPin12b}.

\begin{algorithm}[!htbp]
\caption{Test Space generation by SRQCG.}
\label{alg:SRQCG}
\begin{algorithmic}[1]
\Procedure{SRQCG}{$k_{max}$, $k_{ritz}$, $A$, $G$, $V_0$}
\State Define $S$ as $G A G^{T}$;
\State $Z_0 \leftarrow \left(I - S\right) \mat{V}_0$;
\State Form $Z$ orthonormalizing the column vectors of $Z_0$;
\For{$i=1,n_t$}
  \State Compute Rayleigh quotient $q \leftarrow Z_i^T S Z_i$;
  \State Compute the residual $R_i \leftarrow S Z_i - q Z_i$;
  \State Set $P_i \leftarrow 2 R_i$ and $f_i \leftarrow 1$;
\EndFor
\For{$k=1,\dots,\mbox{convergence}$}
   \If {$\mbox{mod}(k,k_{ritz}) = 0$}
      \State Compute the Ritz projection solving the generalized
      \Statex[4] eigenproblem $\left(Z^T S Z\right) U = \left(Z^T Z\right) U \Lambda$;
      \State $Z \leftarrow Z U$;
   \EndIf
   \For{$i=1,n_t$}
     \State Compute scalar quantities $a \leftarrow P_i^T S Z_i$, $b \leftarrow P_i^T S
       P_i$,
       \Statex[4] $c \leftarrow P_i^T Z_i$, $d \leftarrow P_i^T P_i$ and $e \leftarrow
       Z_i^T S Z_i$;
     \State Compute $\Delta \leftarrow (de-bf_i)^2-4(bc-ad)(af_i-ce)$;
     \State Compute $\alpha \leftarrow (de-bf_i + \sqrt{\Delta})/(2(bc-ad))$;
     \State Update $Z_i \leftarrow Z_i + \alpha P_i$;
     \State Compute Rayleigh quotient $q \leftarrow Z_i^T S Z_i$ and scalar $f_i
       \leftarrow Z_i^T Z_i$;
     \State Compute the residual $R_i \leftarrow S Z_i - q Z_i$ and check convergence;
     \State Update $G_i \leftarrow 2 R_i / f_i$;
     \State Compute $\beta \leftarrow -G_i^T S P_i / b$;
     \State Update $P_i \leftarrow G_i + \beta P_i$;
   \EndFor
   \State {\bf if} $k = k_{max}$ {\bf exit}
\EndFor
\State $\mat{V} \leftarrow \matT{G} \mat{Z}$
\State Orthogonalize the column vectors of $\mat{V}$
\State \Return $\mat{V}$
\EndProcedure
\end{algorithmic}
\end{algorithm}

\subsection{Coarsening}
\label{subsec:coarsening}

Numerous coarsening algorithms have been developed over the years such as classical
Ruge-St\"{u}eben (RS), Cleary-Luby-Jones-Plasman (CLJP), Parallel Maximal Independent Set
(PMIS) and Hybrid MIS (HMIS) \cite{Yan06}. All of them rely on the concept of strength of
connection (SoC), which measures the influence exerted between two neighboring nodes. The
commonly used definition of strength of connection, however, is based on the assumption
that $\mat{A}$ is an M-matrix or it is applied to the M-matrix relative of $\mat{A}$ which
jeopardizes its applicability to more general discretizations. In \amgname{}, we employ
another definition of SoC based on the concept of {\em affinity}, recently introduced by
\cite{LivBra12}, that we believe more flexible and with a wider range of applicability.

Affinity-based SoC requires the availability of a suitable test space that we represent as
a matrix $\mat{V}$ collecting smooth modes on its columns. Let us denote as $\vec{v}_i^T$
the $i$th row vector of $\mat{V}$. Then, the connection strength between two adjacent
degrees of freedom $i$ and $j$ is given by:
\begin{equation}
\mat{SoC}[i,j] = \dfrac{ \left( \vecT{v}_i \vec{v}_j \right)^2} 
                       { \left( \vecT{v}_i \vec{v}_i \right) \left(\vecT{v}_j \vec{v}_j\right)}.
\label{:eq:SoCdef}
\end{equation}
With this definition, the SoC matrix is formed initially on the same pattern of $\mat{A}$
and then filtered by dropping weak connections. Coarse nodes are finally chosen by finding
a maximum independent set (MIS) of nodes on the filtered adjacency graph. Unfortunately,
affinity-based SoC gives usually rise to connections whose numerical values tend to
accumulate in a narrow interval close to one, so it is tricky to define an absolute
threshold for dropping. For this reason, we prefer to control the sparsity of the SoC
matrix, and thus the rapidity of coarsening, by specifying an integer parameter $\theta$,
representing the average number of connections per node retained after filtering.

\subsection{Adaptive prolongation}
\label{subsec:prolongation}

The last key component of our AMG method is the construction of suitable prolongation and
restriction operators. Following the idea proposed in \cite{BraBraKahLiv11} and
successively refined in \cite{MagFraJan19}, we choose to build an interpolation operator
fitting as close as possible the set of test vectors computed in the early setup stage.
More precisely, the prolongation weights $\beta_j$ are computed in order to minimize the
interpolation residual:
\begin{equation}
\| \overline{\bv}_i - \sum_{j \in \J_i} \beta_j \overline{\bv}_j \|_2 \rightarrow
\mbox{min},
\label{minDPLS}
\end{equation}
where $i$ is the node to be interpolated, $\J_i$ is the set of coarse nodes needed to
interpolate the fine node $i$ and $\overline{\bv}^T_k$ represents the $k$th row of the
$\mat{V}$ matrix collecting test vectors. Once the set $\J_i$ is determined, the weights
$\bw_i$ are computed by solving the following linear system of equations:
\begin{equation}
\mat{V}_{c,i} \bw_i = \overline{\bv}_i,
\end{equation}
where the $n_{tv} \times n_i$ matrix $V_{c,i}$ collects the $n_i$ vectors
$\overline{\bv}_k$ for $k \in \J_i$. The main feature of the DPLS algorithm described in
\cite{MagFraJan19} is the dynamic construction during setup of the prolongation pattern,
i.e., the sets $\J_i$ for $i \in \mathcal{C}$. This procedure is designed to reduce at
most \eqref{minDPLS} while preserving sparsity. DPLS iteratively adds entries to $\J_i$
and is controlled by two user-defined parameters: $d_p$, the maximum path length along
strong connections from a node to its interpolating coarse nodes, and $\epsilon_p$, a
relative tolerance controlling the stopping criterion on the interpolation residual. In
other words, the iterations stop whenever all the coarse nodes at a distance less or equal
to $d_p$ have been introduced in $\J_i$ or:
\begin{equation}
\nonumber \| \overline{\bv}_i - \sum_{j \in \J_i} \beta_j \overline{\bv}_j \|_2 \leq
\epsilon_p \| \overline{\bv}_i \|_2.
\end{equation}

From the user standpoint, choosing $d_p$ is relatively easy because, as experimentally
shown, it can take only small values, generally 1 or 2. In fact, the number of neighboring
coarse nodes proliferates with the allowed distance and condition \eqref{minDPLS} is
rapidly satisfied for most of the fine nodes, while choosing higher $d_p$ values only
burdens the setup with no significant benefits regarding convergence speed. On the other
hand, the prolongation quality is susceptible to $\epsilon_p$, making the proper choice of
this parameter difficult and extremely problem-dependent, especially in structural
mechanics. By contrast to intuition, it is not sufficient to prescribe a very small value
for $\epsilon_p$ to have a good prolongation operator, because, even if this choice
practically guarantees that $\mbox{range}(\mat{P})$ perfectly represents $\mat{V}$ when at
least $n_{tv}$ coarse neighbors have been selected for each fine node, the prolongation
operator becomes severely ill-conditioned deteorating the AMG effectiveness on the next
levels. Only relatively few coarse neighbors carry meaningful information to interpolate
$\overline{\bv}_i$ while the others are, in some sense, redundant. This is intrinsically
connected to the nature of rigid body modes that, in structural mechanics, give a good
approximation of the near-kernel space of $\mat{A}$. For several fine nodes, it may happen
that, due to an unfavorable geometrical placement of the mesh nodes, the matrix $V_{c,i}$
even if full-rank is severely ill-conditioned thus determining abnormally large weights in
$\mat{P}$. To control this unfortunate occurrence, we introduce a different stopping
criterion based on the user-defined tolerance $\kappa_p$, representing the maximum
allowable conditioning for $\mat{V}_{c,i}$. While including more coarse nodes in the
interpolation, we inexpensively monitor the condition number of $V_{c,i}$ during its
Householder (HH) QR factorization, e.g., see \cite{GolLoa13}, and, whenever
$\mbox{cond}(V_{c,i}) > \kappa_p$ we stop the iteration avoiding too large values of the
weights. The procedure to set up the prolongation is briefly summarized in Algorithm
\ref{alg:DPLS}.

\begin{algorithm}[t!]
\caption{\bf Dynamic Pattern Least Squares Prolongation}
\label{alg:DPLS}
\begin{algorithmic}[1]
\Procedure{DPLS\_SetUp}{$\kappa_p$,$S_c$,$V_f$,$V_c$,$W$}
\ForAll {nodes $i \in \mathcal{F}$}
\State $k \leftarrow 0$; $\mathcal{C}_i \leftarrow \varnothing$; $\br \leftarrow \bv_i$;
  $\bbv_i \leftarrow \bv_i$; $R \leftarrow \varnothing$;
\State Form $\mathcal{J}_i \leftarrow \{j \in \mathcal{C} \: \: | \:$ there is a
path from $i$ to $j$ in the $S_c$ graph$\}$;
\While {$\mbox{cond}(R) \leq \kappa_p$}
\State $k \leftarrow k + 1$;
\State Select $\bar{j} \in \mathcal{J}_i \setminus \mathcal{C}_i$ for which $\bbv_j$
has maximal affinity with $\br$;
\State Update $\mathcal{C}_i \leftarrow \mathcal{C}_i \cup \{ \bar{j} \}$;
\State Add $\bbv_{\bar{j}}$ as the last column of $V_{c,i}$;
\State Compute the HH reflection $\mat{Q}$ nullifying last $n_{tv}-k$ rows of $R$;
\State Compute $\br \leftarrow \mat{Q} \, \br$;
\ForAll {$j \in \mathcal{J}_i \setminus \mathcal{C}_i$}
\State Compute $\bbv_j \leftarrow \mat{Q} \, \bbv_j$;
\EndFor
\EndWhile
\State Compute $\bw_i \leftarrow R^{-1}\br$;
\EndFor
\EndProcedure
\end{algorithmic}
\end{algorithm}

\section{Numerical results on real-world structural problems}
\label{sec:numexp}

In this section, we evaluate the performance of the conjugate gradient preconditioned with
the algebraic multigrid method proposed in this work (\amgname{}/PCG) in the solution of
challenging real-world structural problems. In particular, this section initially presents
the CPU times required by the different solution phases of our method with the aim of
understanding the impact of the preconditioner input parameters on performance and
evaluating the influence of each parameter. These results allow identifying an optimal set
of parameters for the algorithm. Then, we compare our linear solver with two
state-of-the-art algebraic multigrid preconditioners used along with PCG: the first,
BoomerAMG, a classical AMG implemented in the HYPRE package \cite{FalYan02}, and the
second one, GAMG, an aggregation-based AMG implemented in the PETSc library \cite{PETSc}.
Lastly, we provide strong scalability results of our implementation for a shared memory
architecture.

\begin{figure}[!htbp]
  \centering 
  \includegraphics[trim={0 0 2cm 0},clip,width=0.32\linewidth]{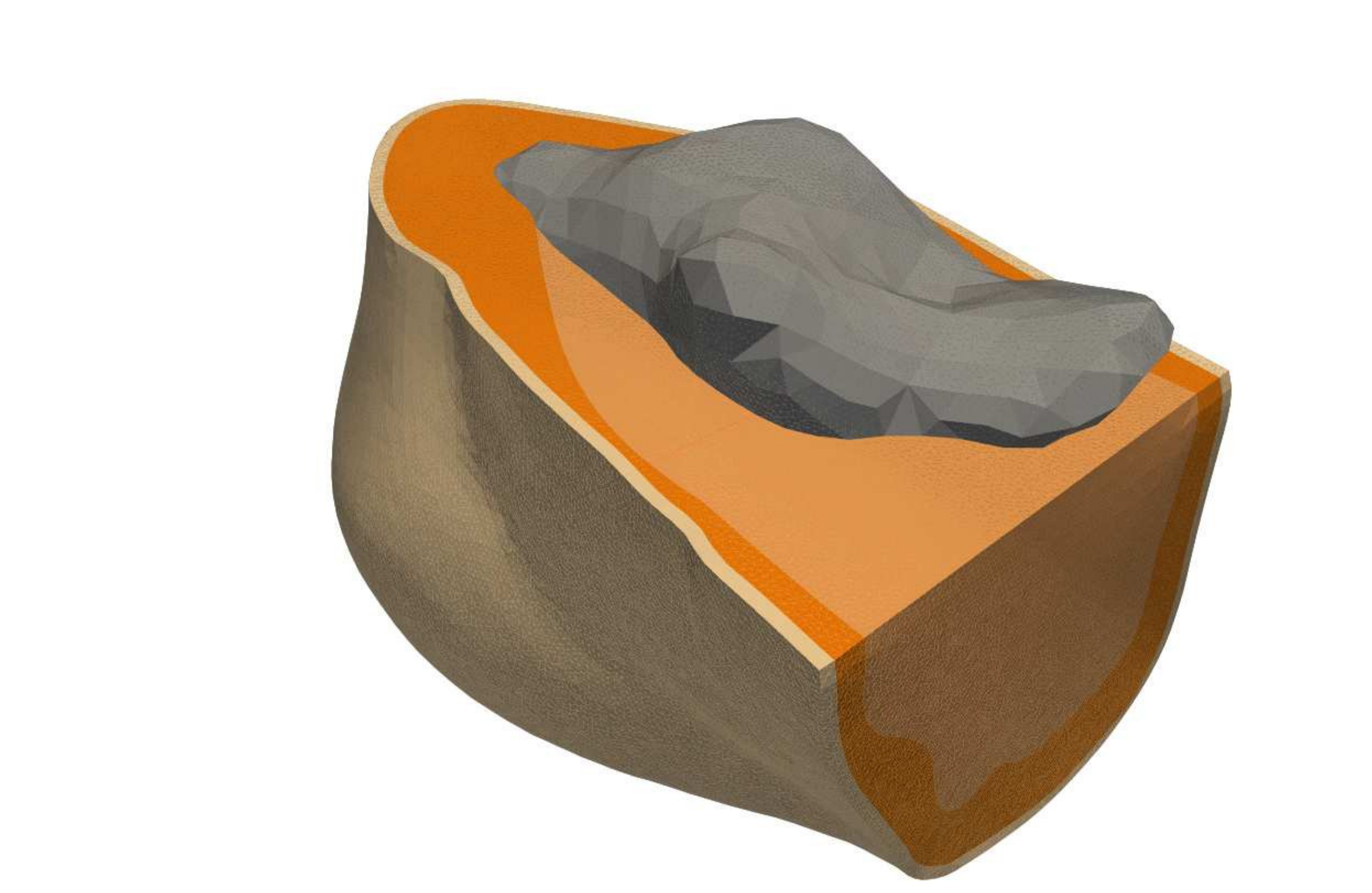}
  \includegraphics[trim={0 0 2cm 0},clip,width=0.32\linewidth]{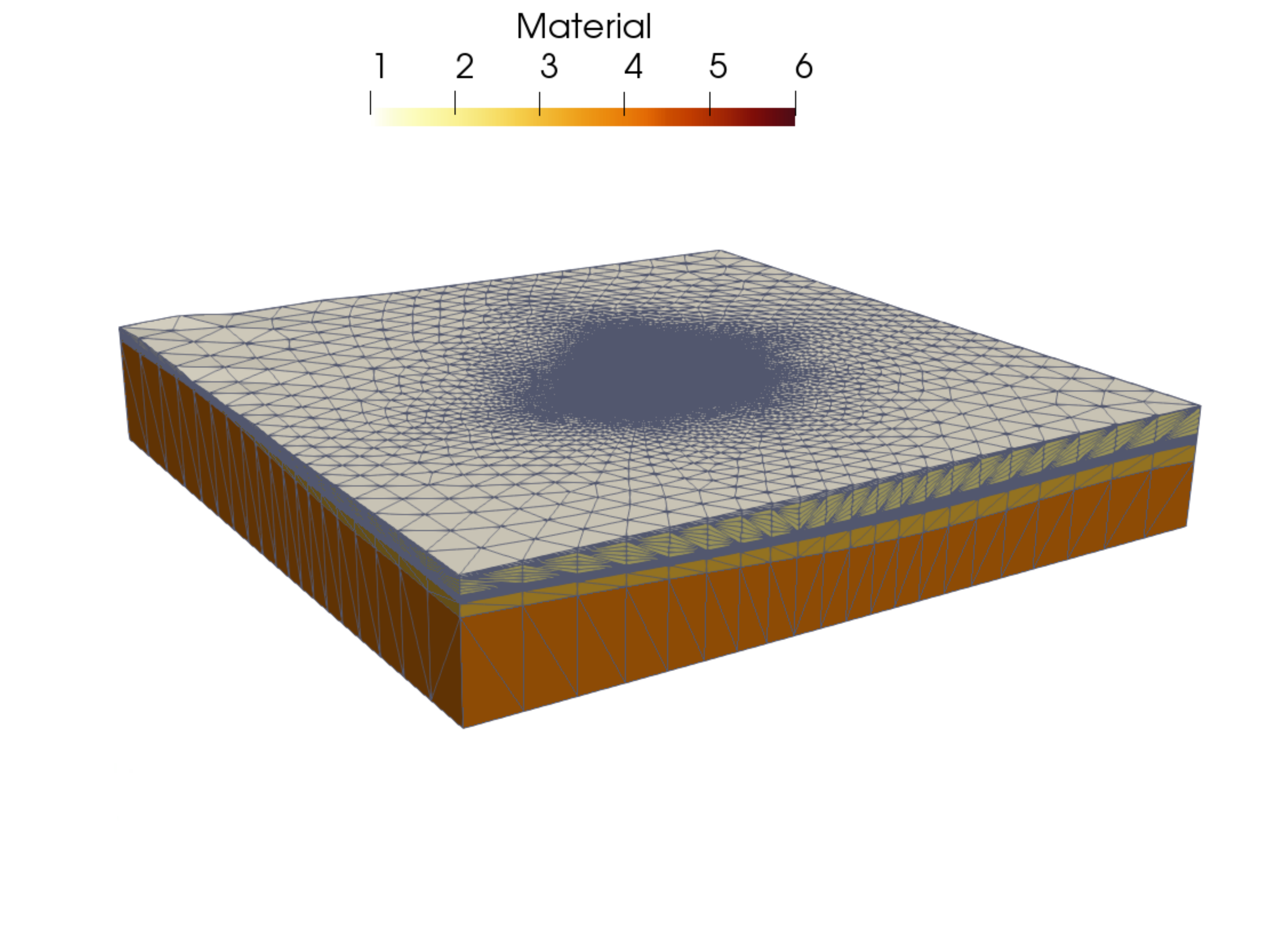}
  \includegraphics[trim={0 0 2cm 0},clip,width=0.32\linewidth]{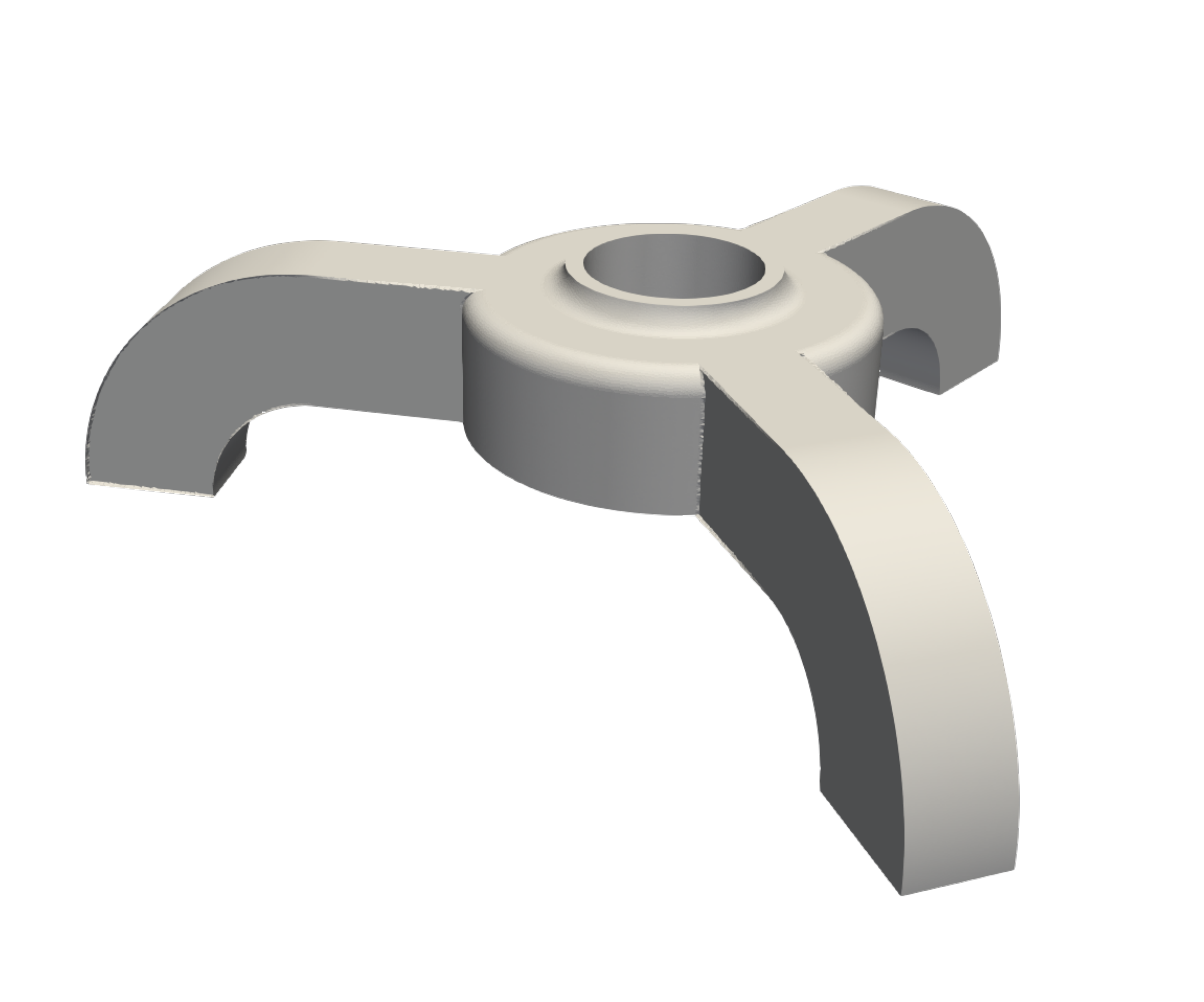} \\
  \includegraphics[trim={0 0 2cm 0},clip,width=0.32\linewidth]{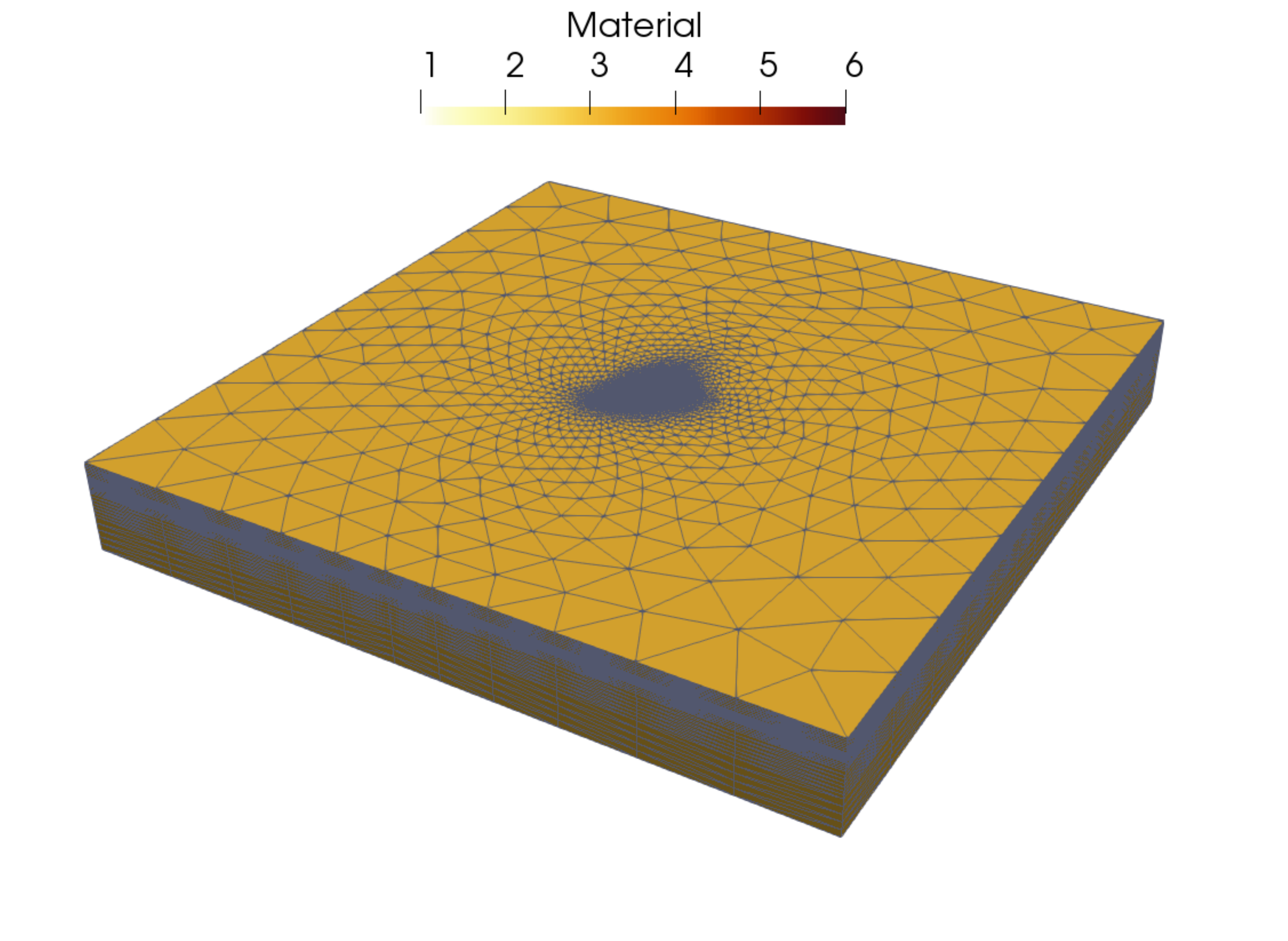}
  \includegraphics[trim={0 0 2cm 0},clip,width=0.32\linewidth]{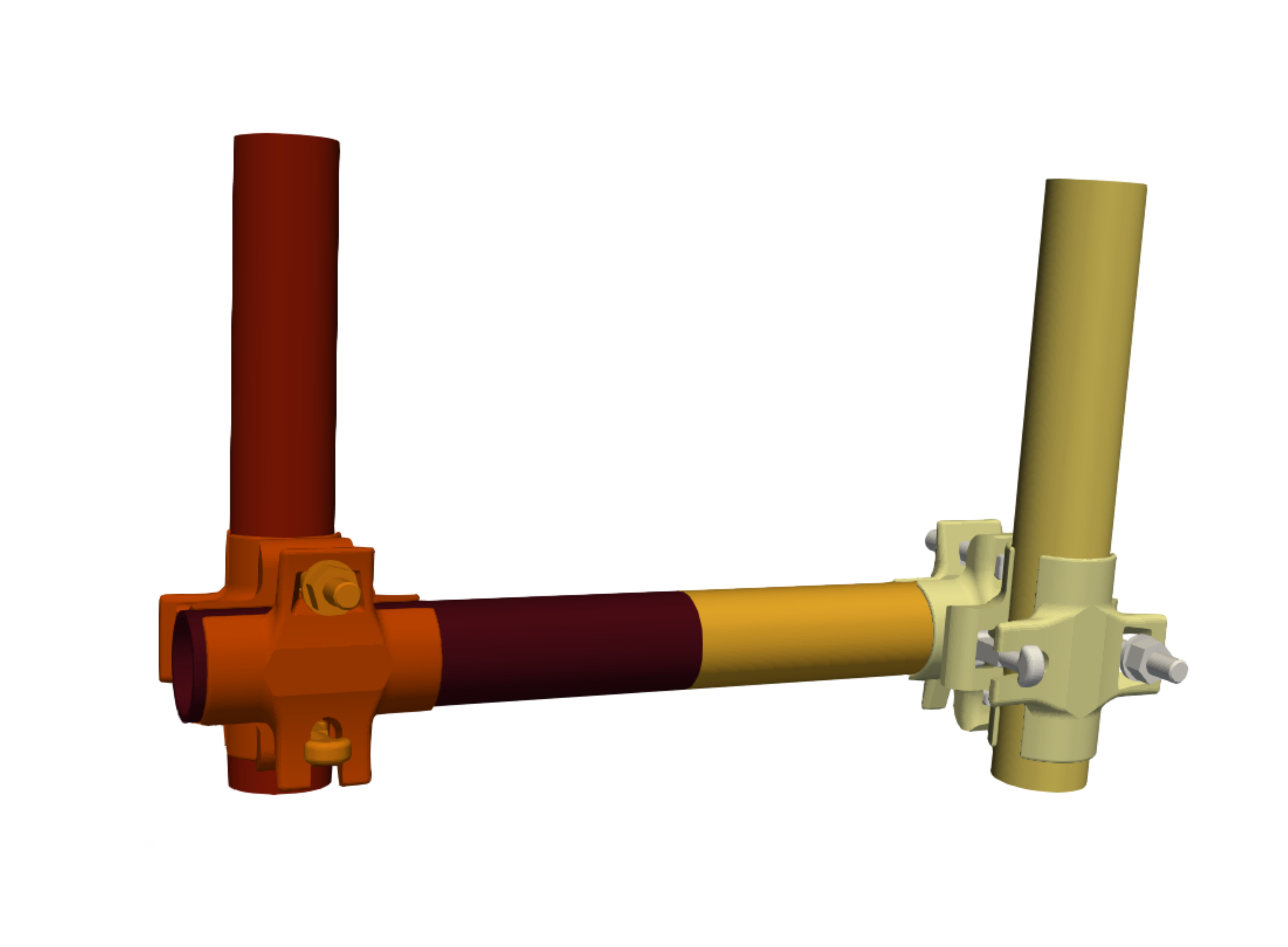}
  \includegraphics[trim={0 0 2cm 0},clip,width=0.32\linewidth]{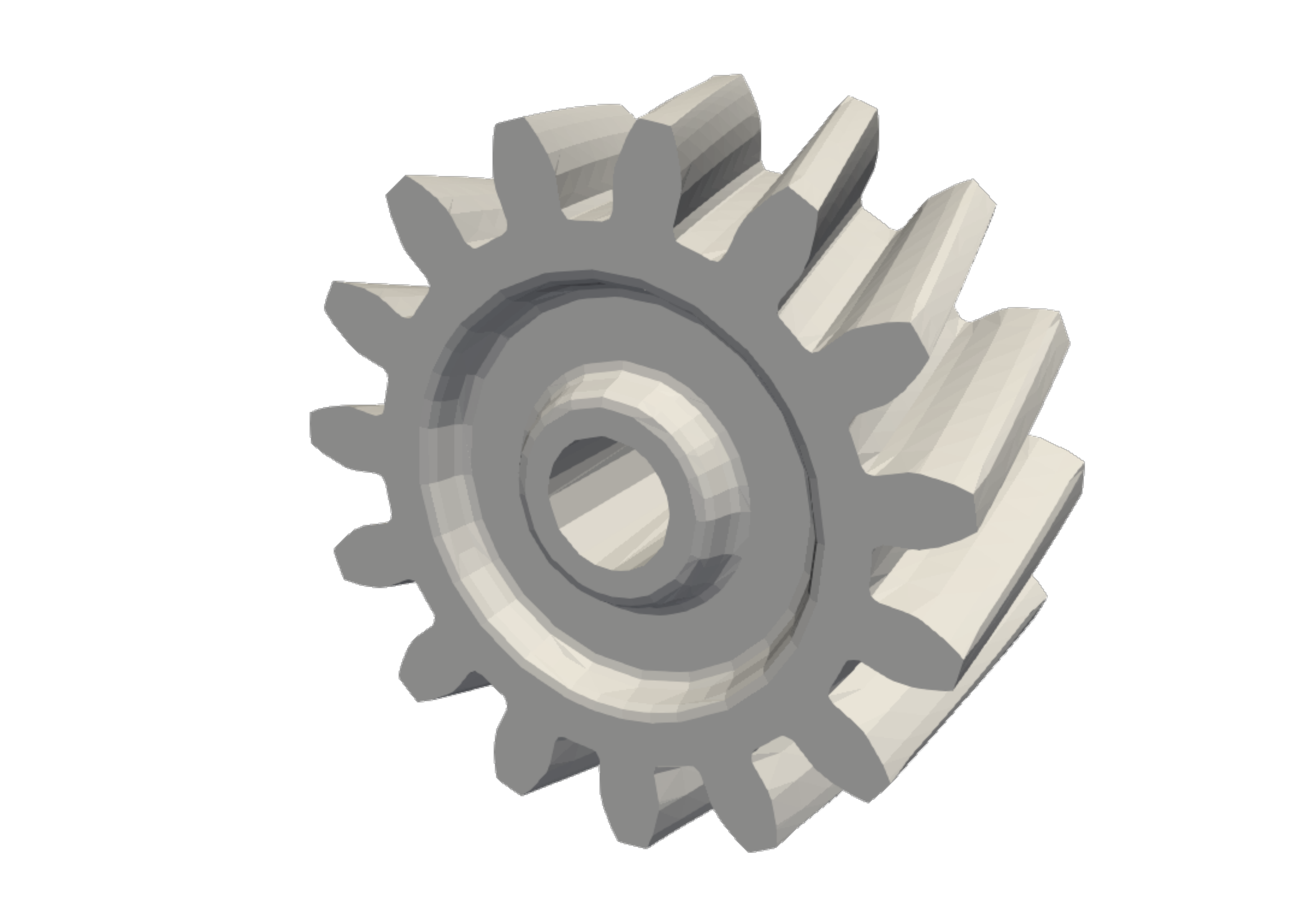} \\
  \includegraphics[trim={0 0 2cm 0},clip,width=0.32\linewidth]{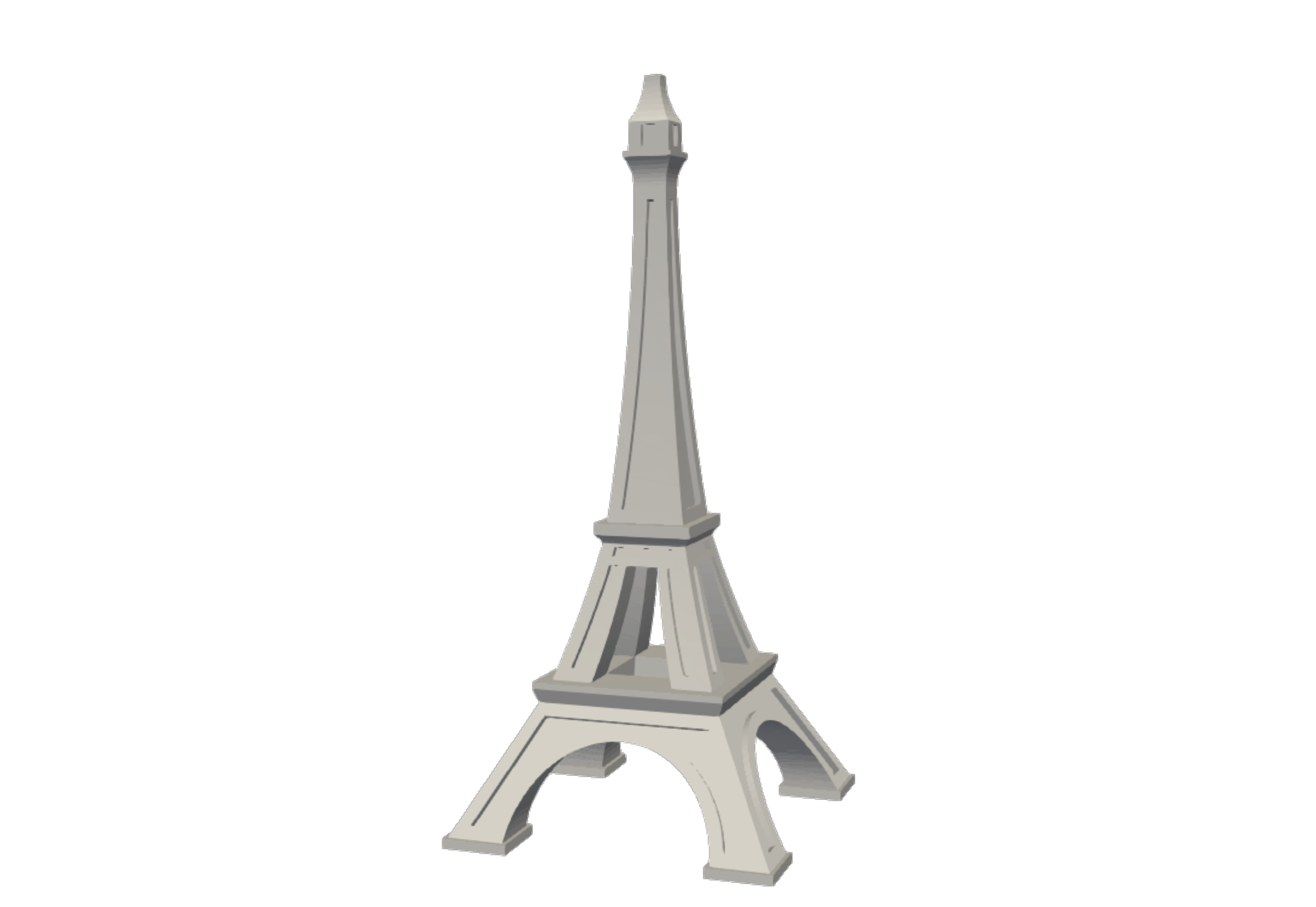}
  \includegraphics[trim={0 0 2cm 0},clip,width=0.32\linewidth]{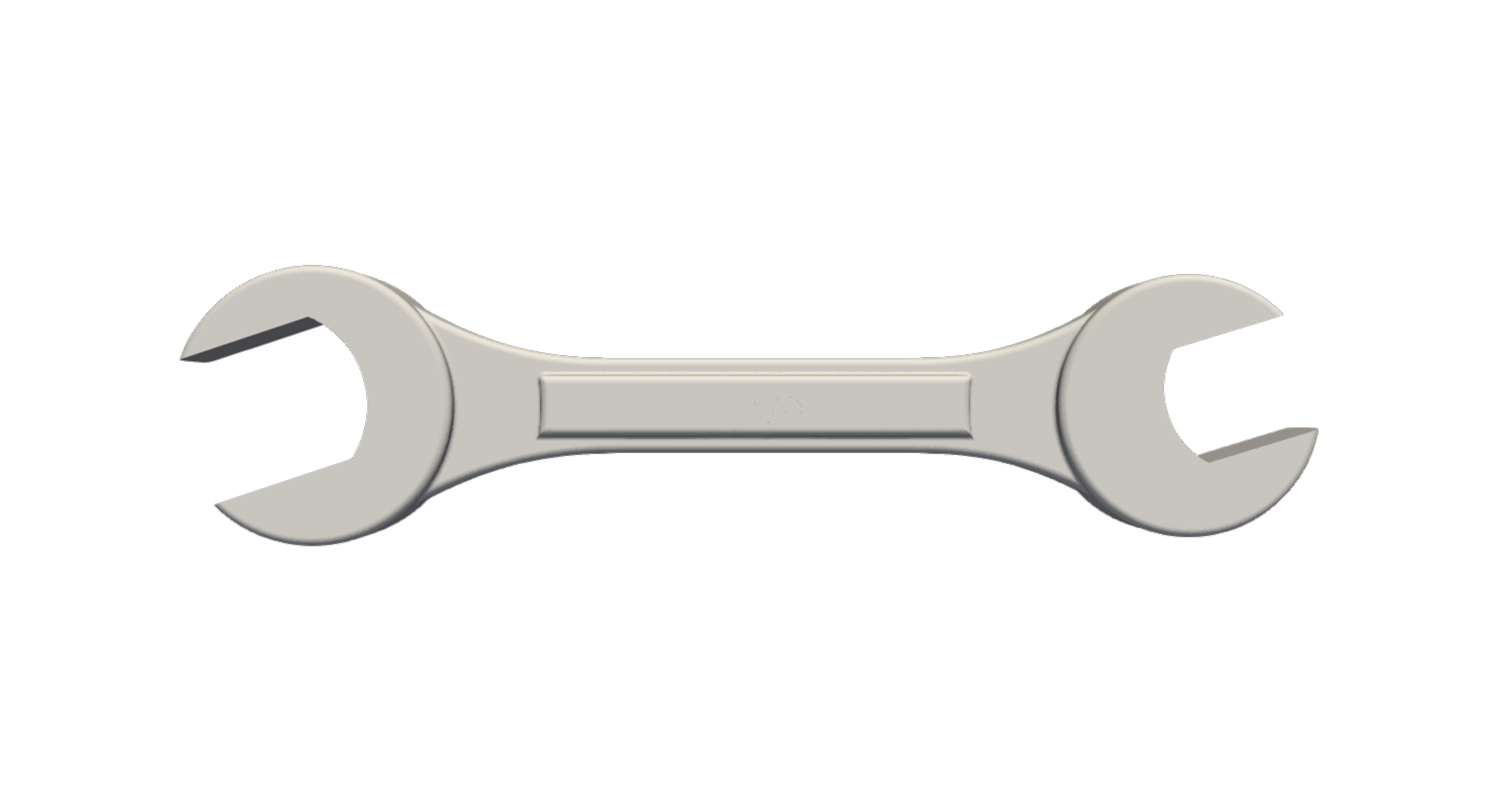}
  \includegraphics[trim={0 0 2cm 0},clip,width=0.32\linewidth]{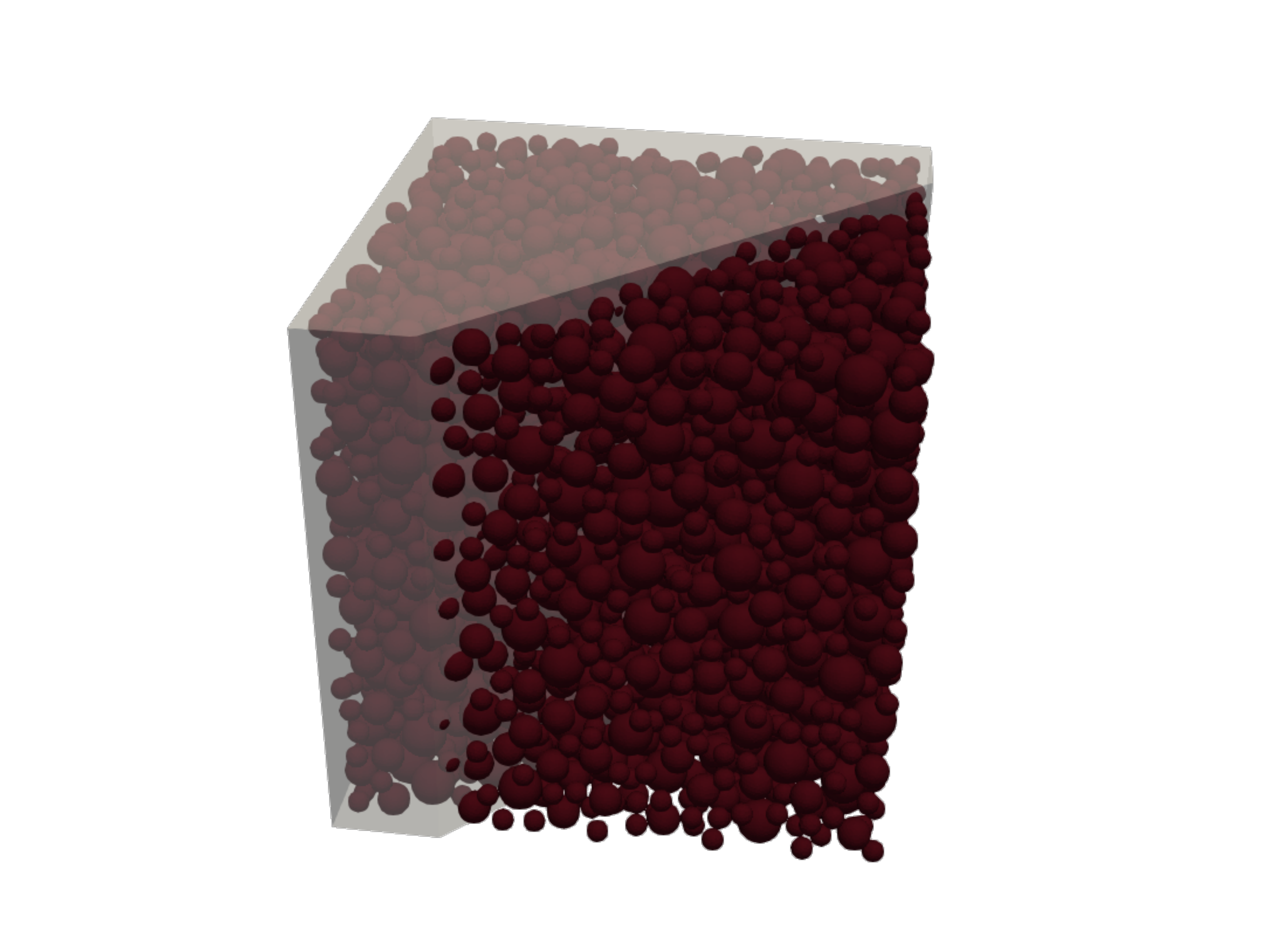}
  \caption{Global view of the meshes referent to the real-world test cases. From left to
    right: (top row) {\tt heel1138k}, {\tt guenda1446k}, {\tt tripod3239k}; (middle row)
    {\tt thdr3559k}, {\tt jpipe6557k}, {\tt gear8302k}; (bottom row) {\tt eiffel9213k},
    {\tt wrench13m}, {\tt agg14m}.}
\label{fig:meshall}
\end{figure}

The sparse matrices are derived from the finite element discretization of real-world
applications (some of them sketched in Figure \ref{fig:meshall}), modeled by the linear
elasticity PDE.  These are challenging test cases, not only for the high number of degrees
of freedom (DOFs), but also because they are characterized by several sources of
ill-conditioness, such as the presence of multiple materials ({\sf MM}), of highly
distorted elements ({\sf DE}), of material close to the incompressible limit, i.e., $\nu
\rightarrow 0.5$, ({\sf IM}) and finally due to loosely constrained body ({\sf LC}).
Table \ref{tab:realWorldProblems} shows a summary with the fundamental data of each test
case, while \ref{sec:meshGrimpse} describes accurately each test case and
\ref{sec:meshQuality} gives details about the mesh quality. The right-hand side vector
used for all test cases is the unitary vector.
\begin{table}[htbp]
\ra{1.4} \centering
\caption{Test matrices derived from the finite element discretization of real-world
  structural models. For each test problem, the table provides the size, $n(\mat{A})$, the
  number of non-zeroes, $nnz(\mat{A})$, and the kind of discretization, i.e., linear
  tetrahedra, P1, or hexahedra, Q1, along with the source/sources of ill-conditioning,
  ({\sf MM, DE, IM, LC}).}
\label{tab:realWorldProblems}
\resizebox{\textwidth}{!}{%
\begin{tabular}{@{}*{1}{l}|*{2}{r}|*{1}{r}|*{1}{r}|*{1}{r}|*{1}{r}@{}}
\toprule
\bf{Matrix name}  & $n(\mat{A})$ & $nnz(\mat{A})$ & \bf{FEM} & \bf{Application} & \bf{Description}      & \bf{Ill-conditioning} \\
\midrule
\tt{heel1138k}    &  1,138,443   &     51,677,937 &       P1 &  3D Biomedicine  & \ref{sec:heel1138k}   & {\sf MM} \\
\tt{guenda1446k}  &  1,446,624   &     64,374,678 &       P1 &  3D Geomechanics & \ref{sec:guenda1446k} & {\sf MM,DE} \\
\tt{hook1498k}    &  1,498,023   &     59,374,451 &       P1 &  3D Mechanical   & \ref{sec:hook1498k}   & {\sf LC} \\
\tt{tripod3239k}  &  3,239,649   &    142,351,857 &       P1 &  3D Mechanical   & \ref{sec:tripod3239k} & {\sf IM} \\
\tt{thdr3559k}    &  3,559,398   &     81,240,330 &       P1 &  3D Geomechanics & \ref{sec:thdr3559k}   & {\sf MM,DE} \\
\tt{beam6502k}    &  6,502,275   &    515,265,317 &       Q1 &  3D Mechanical   & \ref{sec:beam6502k}   & {\sf MM,LC} \\
\tt{jpipe6557k}   &  6,557,808   &    335,451,702 &       P1 &  3D Mechanical   & \ref{sec:jpipe6557k}  & {\sf MM,LC} \\
\tt{gear8302k}    &  8,302,032   &    672,580,800 &       P2 &  3D Mechanical   & \ref{sec:gear8302k}   & {\sf IM,DE} \\
\tt{eiffel9213k}  &  9,213,342   &    390,108,294 &       P1 &  3D Mechanical   & \ref{sec:eiffel9213k} & {\sf LC} \\
\tt{guenda11m}    & 11,452,398   &    512,484,300 &       P1 &  3D Geomechanics & \ref{sec:guenda1446k} & {\sf MM,DE} \\
\tt{wrench13m}    & 12,995,364   &    589,759,812 &       P1 &  3D Mechanical   & \ref{sec:wrench13m}   & {\sf IM,LC} \\
\tt{agg14m}       & 14,106,408   &    633,142,730 &       P1 &  3D Mesoscale    & \ref{sec:agg14m}      & {\sf MM} \\
\tt{M20}          & 20,056,050   &  1,634,926,088 &       P2 &  3D Mechanical   & \ref{sec:M20}         & {\sf DE,LC} \\
\bottomrule
\end{tabular}}
\end{table}
The linear systems are solved by PCG with initial solution equal to the null vector, and
convergence is achieved when the $l2$-norm of the iterative residual vector becomes
smaller than $10^{-8} \cdot \twonorm{\vec{b}}$.
Numerical experiments are run on the MARCONI-A2 cluster from the Italian Supercomputing
Center (CINECA), with a single node composed by one Intel Xeon Phi 7250 CPU (Knights
Landing) at 1.40 GHz with 68-cores, and 86GB of DDR4 RAM plus 16GB MCDRAM used in
cache/quadrant mode. Except for the scalability test, we always used all of the 68
computing cores via MPI, in case of BoomerAMG and GAMG, or via OpenMP, in case of
\amgname{}.

\subsection{Default parameters study}
\label{subsec:sensitivity}

\begin{figure}[!htbp]
  \centering 
  \includegraphics[trim={0 0 2cm 0},clip,width=0.32\linewidth]{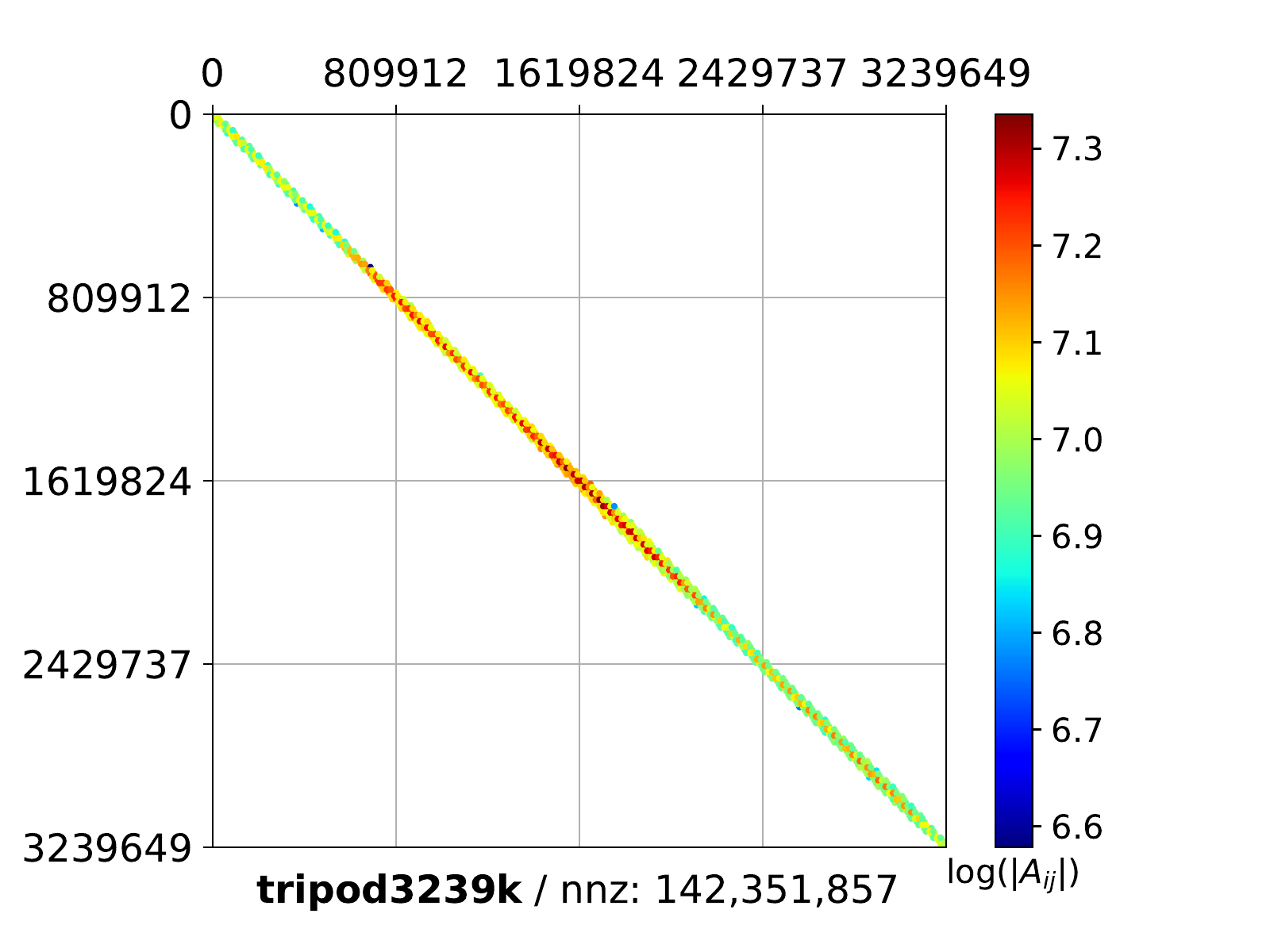}
  \includegraphics[trim={0 0 2cm 0},clip,width=0.32\linewidth]{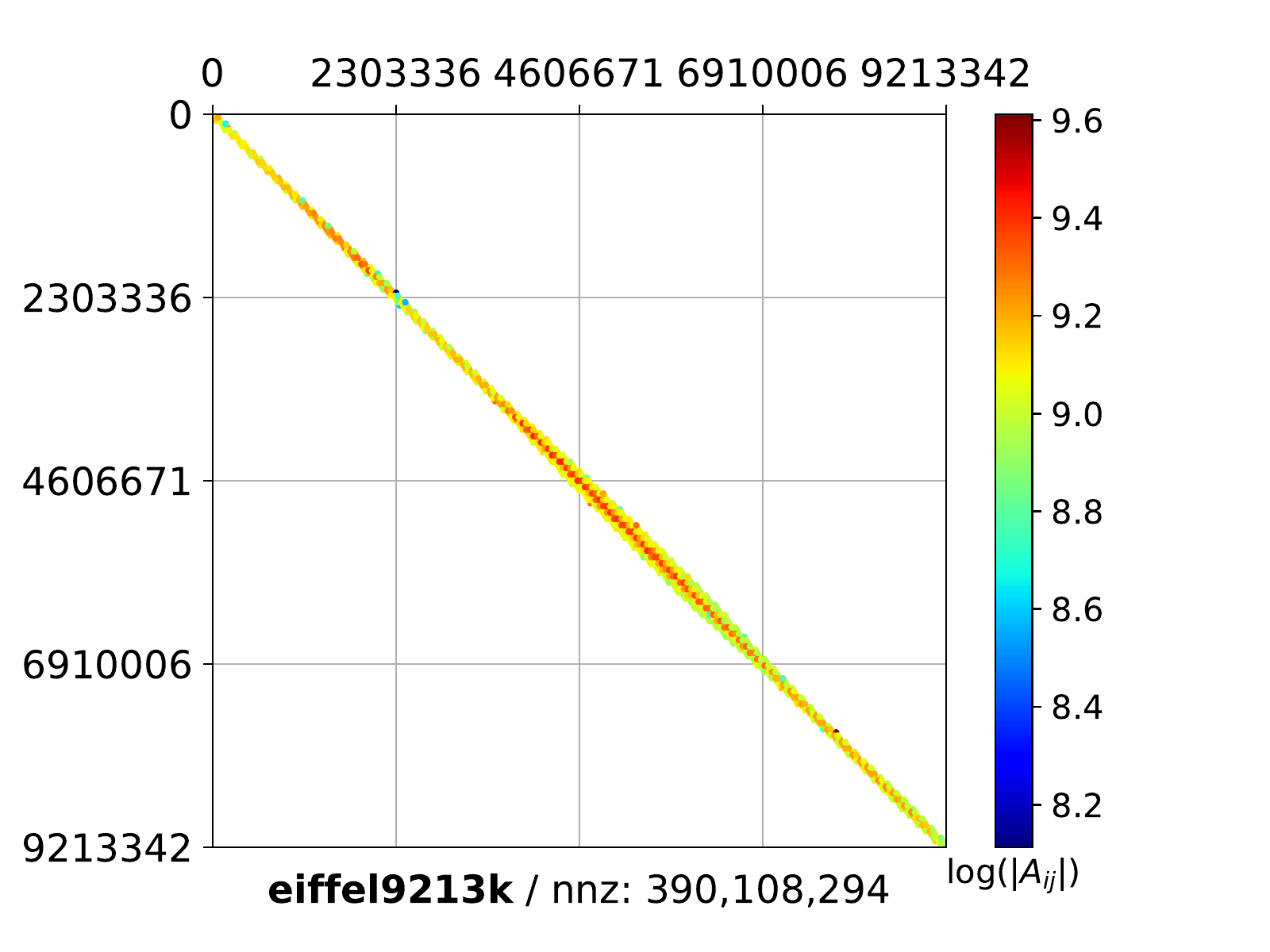}
  \includegraphics[trim={0 0 2cm 0},clip,width=0.32\linewidth]{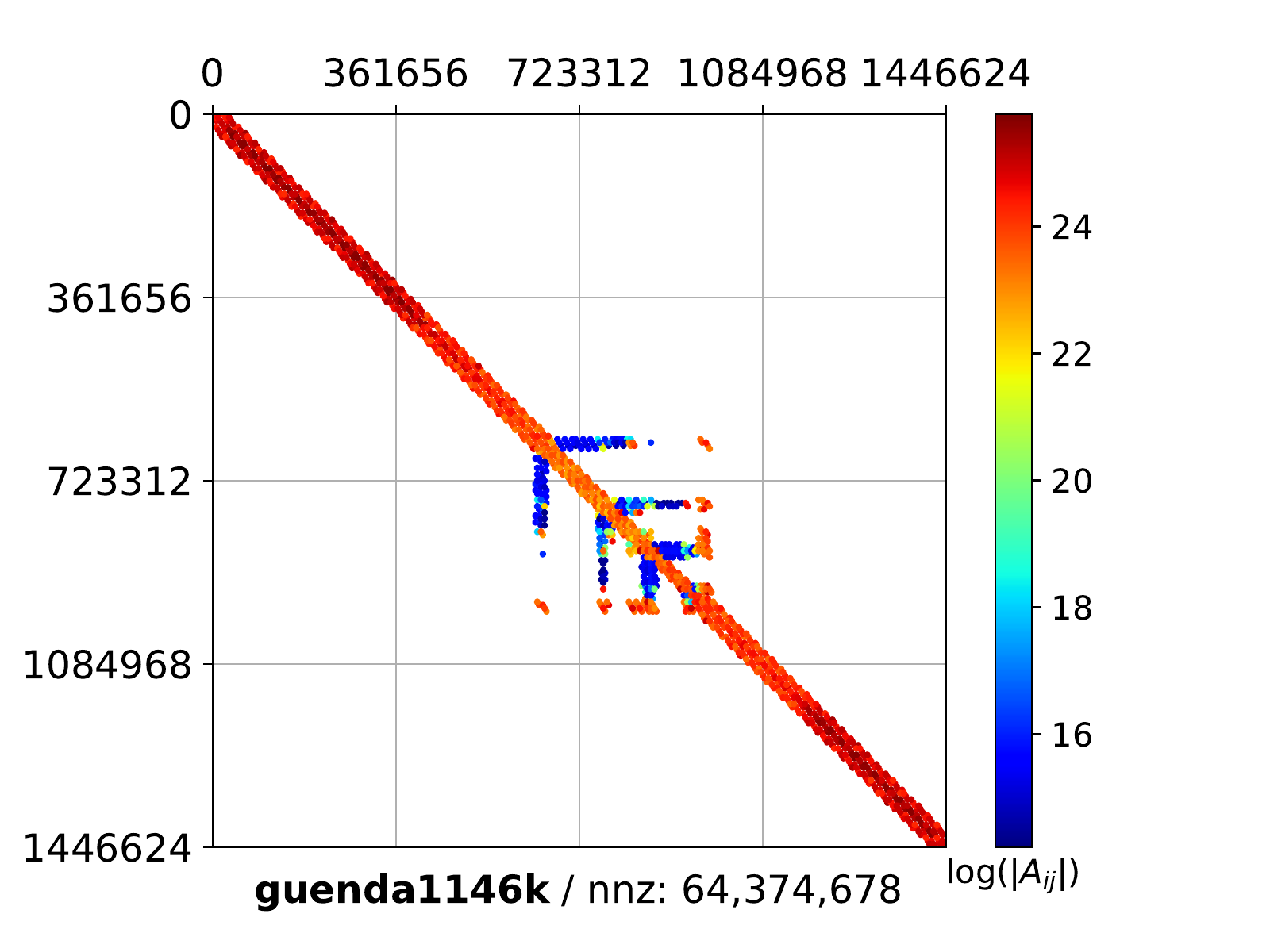}
  \caption{Sparsity pattern plots of the test matrices {\tt tripod3239k}, {\tt
      eiffel9213k}, and {\tt guenda1446k}, from left to right.}
  \label{fig:SpPatt}
\end{figure}

The \amgname{} solver depends upon eight user-defined parameters. Thus, it is desirable to
understand which of them plays a significant role in the construction of an efficient
preconditioner, as well as their range of suitable values. For this reason, we select a
smaller portion of test cases from Table \ref{tab:realWorldProblems} that are
representative of different sources of ill-conditioning for structural problems: {\tt
  tripod3239k} ({\bf T1}) for materials close to the incompressible limit ({\sf IM}); {\tt
  eiffel9213k} ({\bf T2}) for loosely constrained bodies ({\sf LC}), and {\tt guenda1446k}
({\bf T3}) for multiple materials ({\sf MM}) and meshes with highly distorted elements
({\sf DE}). The sparsity pattern of the above matrices is shown in Figure
\ref{fig:SpPatt}.  Then, we run a sensitivity analysis on \amgname{} by testing different
setup configurations starting from a default configuration. To better understand the
relative impact of each parameter, we vary only one or few of them at a time, leaving the
others fixed to the ``default'' values as listed in Table \ref{tab:TestParam}.

\begin{table}[!htbp]
\ra{1.2}
\centering
\caption{Default input parameters for \amgname{}.}
\label{tab:TestParam}
\begin{tabular}{@{}*{1}{l}|*{2}{l}@{}}
\toprule
\bf{Phase}                    &  \bf{Param.} &    \bf{Value} \\
\midrule
\multirow{3}{*}{Smoother}     &        $k_g$ &           $4$ \\
                              &     $\rho_g$ &           $4$ \\
                              & $\epsilon_g$ &      $10^{-3}$ \\ \hline
\multirow{2}{*}{Test Space}   &      $n_{tv}$ &          $10$ \\
                              &      $n_{rq}$ &          $10$ \\ \hline
\multirow{1}{*}{Coarsening}   &     $\theta$ &           $5$ \\ \hline
\multirow{2}{*}{Prolongation} &     $n_{max}$ &           $5$ \\
                              &   $\kappa_p$ &          $50$ \\
\bottomrule
\end{tabular}
\end{table}

Numerical results considering the default input parameters of \amgname{} are given in
Table \ref{tab:sstvDefRes}. In the following, we test how the performance of the method
changes when varying the input parameters relative to each of its construction phases. To
measure the preconditioner densities and the time required for all solution phases, we
introduce some symbols:
\begin{itemize}
  \item $C_{gd}$: grid complexity, i.e., the sum of the number of rows for the matrix of
    each level divided by the number of rows of the original matrix;
  \item $C_{op}$: operator complexity, i.e., the sum of the number of entries for the
    matrix of each level divided by the number of entries of the original matrix;
  \item $C_{fs}$: FSAI complexity, i.e., an average value of the aFSAI
    density over the multigrid hierarchy, as defined in \cite{MagFraJan19};
  \item $n_{it}$: number of PCG iterations needed to reach convergence;
  \item $T_{ts}$: total time (in seconds) for the test space generation;
  \item $T_{cs}$: total time (in seconds) for the coarsening phase;
  \item $T_{sm}$: total time (in seconds) for the smoother computation;
  \item $T_{pl}$: total time (in seconds) for the prolongation computation;
  \item $T_{rap}$: total time (in seconds) to perform the product RAP;
  \item $T_{p}$: time (in seconds) spent in the preconditioner setup;
  \item $T_{s}$: time (in seconds) spent in the PCG iteration phase;
  \item $T_{t}$: total time (in seconds) to solve the linear system.
\end{itemize}

\begin{table}[!htbp]
  \ra{1.4}
  \centering
  \caption{\amgname{} results considering the input parameters from Table \ref{tab:TestParam}.}
  \label{tab:sstvDefRes}
  \resizebox{\textwidth}{!}{%
  \begin{tabular}{@{}*{1}{r}|*{12}{r}@{}}
  \toprule
  \multirow{2}{*}{\thead{\bf{Test} \\ \bf{case}}} & \multicolumn{12}{c}{\bf{Results}} \\
  \cmidrule(lr){2-13}
           & $C_{gd}$ & $C_{op}$ & $C_{fs}$ & $n_{it}$ & $T_{ts}$ & $T_{cs}$ & $T_{sm}$  & $T_{pl}$ & $T_{rap}$ &  $T_p$ & $T_s$ &  $T_t$ \\
  \midrule
  {\bf T1} &    1.39 &    3.08 &    0.53 &     228 &    18.4 &     9.5 &      6.8 &      3.8 &      2.4 &   42.4 &  23.4 &  65.8 \\
  {\bf T2} &    1.38 &    2.75 &    0.55 &     239 &    49.3 &    25.1 &     33.5 &      8.4 &      4.6 &  123.0 &  82.4 & 205.4 \\
  {\bf T3} &    1.48 &    2.44 &    0.57 &      81 &     9.2 &     4.2 &      3.8 &      1.8 &      1.4 &   21.4 &   4.5 &  26.2 \\
  \bottomrule
  \end{tabular}}
\end{table}

{\bf Coarsening.}
In Table \ref{tab:sstvCoarsening}, we show the sensitivity analysis of \amgname{} with
respect to its coarsening calculation phase. We want to test the performance of the method
varying the average number of strong connections per node $\theta$ from a lower bound
equal to $2$ up to an upper bound equal to $10$. A common behavior to all numerical runs
is the decrease in the grid and operator complexities when increasing the value of
$\theta$, which is followed by an overall decrease in the total solution time with the
optimal condition reached when $\theta = 8$. Moreover, looking at the iteration counts, we
conclude that the preconditioning operator quality improves when increasing $\theta$
towards $10$ with the only exception of {\bf T3}, which reaches the minimum number of
iterations with $\theta = 4$. As expected, the parameter $\theta$ affects mostly the
coarsening and RAP times ($T_{cs}$, $T_{rap}$), with the last one being naturally related
to coarsening. However, it is worth noting that the test space and smoother computation
times ($T_{ts}$ and $T_{sm}$) are influenced indirectly by the way that the multigrid
hierarchy is built, and consequently, by the value of $\theta$. In summary, we saw that
such parameter has a significant effect on all setup phases of \amgname{} and therefore is
an important parameter in our method. Moreover, we note that the default choice of
$\theta$ lead to a satisfactory performance in all test cases being approximately $25 \%$
slower than the optimal configuration in the worst case.

\begin{table}[!htbp]
  \ra{1.4}
  \centering
  \caption{Sensitivity analysis with the coarsening parameters.}
  \label{tab:sstvCoarsening}
  \resizebox{\textwidth}{!}{%
  \begin{tabular}{@{}*{1}{r}|*{1}{r}|*{1}{r}|*{12}{r}@{}}
  \toprule
  \multirow{2}{*}{\thead{\bf{Test} \\ \bf{case}}} & \multirow{2}{*}{\thead{\bf{Run} \\ \bf{ID}}} & \multicolumn{1}{c|}{\bf{Param.}} & \multicolumn{12}{c}{\bf{Results}} \\
  \cmidrule(lr){3-15}
           &   & $\theta$ & $C_{gd}$ & $C_{op}$ & $C_{fs}$ & $n_{it}$ & $T_{ts}$ & $T_{cs}$ & $T_{sm}$ & $T_{pl}$ & $T_{rap}$ & $T_p$ & $T_s$ & $T_t$  \\
  \midrule
  {\bf T1} & 1 &  {\bf 2} &    1.66 &    4.54 &    0.64 &      316 &   22.0 &     13.2 &     9.3 &     3.7 &      2.9 &  53.1 &  44.9 &  98.7 \\
           & 2 &  {\bf 4} &    1.45 &    3.57 &    0.56 &      248 &   19.5 &     10.3 &     7.6 &     4.1 &      2.9 &  45.6 &  28.0 &  74.0 \\
           & 3 &  {\bf 6} &    1.34 &    2.79 &    0.51 &      224 &   17.4 &      8.4 &     6.2 &     3.6 &      2.3 &  39.0 &  20.7 &  60.1 \\
           & 4 &  {\bf 8} &    1.26 &    2.39 &    0.49 &      222 &   16.3 &      7.4 &     5.6 &     3.7 &      1.9 &  36.2 &  18.6 &  55.1 \\
           & 5 & {\bf 10} &    1.21 &    2.11 &    0.47 &      229 &   15.5 &      7.0 &     5.3 &     4.2 &      1.5 &  35.1 &  18.2 &  53.6 \\
  \midrule
  {\bf T2} & 1 &  {\bf 2} &    1.66 &    4.07 &    0.66 &      379 &   60.7 &     34.5 &    44.0 &     7.7 &      6.5 & 156.8 & 197.2 & 355.0 \\
           & 2 &  {\bf 4} &    1.44 &    3.10 &    0.58 &      260 &   53.0 &     27.1 &    35.8 &     8.0 &      5.4 & 131.7 & 101.2 & 233.6 \\
           & 3 &  {\bf 6} &    1.33 &    2.49 &    0.53 &      228 &   48.5 &     22.9 &    31.8 &     8.2 &      4.2 & 117.8 &  73.5 & 191.9 \\
           & 4 &  {\bf 8} &    1.26 &    2.17 &    0.50 &      221 &   44.0 &     20.1 &    29.4 &     8.6 &      3.6 & 108.0 &  64.9 & 173.6 \\
           & 5 & {\bf 10} &    1.20 &    1.92 &    0.48 &      210 &   41.2 &     18.4 &    27.8 &     9.2 &      3.2 & 101.8 &  53.6 & 155.9 \\
  \midrule
  {\bf T3} & 1 &  {\bf 2} &    1.71 &    3.01 &    0.65 &       99 &   10.8 &      5.1 &     4.5 &     1.5 &      1.4 &  24.7 &   6.2 &  31.2 \\
           & 2 &  {\bf 4} &    1.52 &    2.61 &    0.58 &       86 &    9.6 &      4.5 &     4.0 &     1.7 &      1.8 &  22.8 &   4.7 &  27.7 \\
           & 3 &  {\bf 6} &    1.45 &    2.34 &    0.55 &       91 &    9.1 &      3.9 &     4.1 &     1.8 &      1.5 &  21.8 &   4.6 &  26.5 \\
           & 4 &  {\bf 8} &    1.38 &    2.18 &    0.52 &      108 &    8.8 &      3.9 &     3.6 &     2.1 &      1.5 &  20.8 &   5.2 &  26.2 \\
           & 5 & {\bf 10} &    1.29 &    2.04 &    0.49 &      375 &    7.8 &      3.7 &     3.3 &     2.1 &      1.4 &  19.2 &  14.5 &  33.9 \\
  \bottomrule
  \end{tabular}}
\end{table}

{\bf Prolongation.}
In Table \ref{tab:sstvProlongation}, we report the sensitivity analysis of \amgname{} with
respect to its prolongation calculation phase. We vary the maximum number of interpolation
values per node, $n_{max}$ and the conditioning bound of the local least squares problems
$\kappa_p$. From the numerical runs, we note that these parameters do not affect the
average coarsening ratio, since $C_{gd}$ is practically constant. However, the operator
complexity $C_{op}$ is changed substantially, being the prolongation operator responsible
for the coarse grid operators density at all levels. Thus, it makes sense that the
most affected phases in such analysis are the coarsening and prolongation computation
itself as evidenced by the computational times $T_{pl}$, $T_{cs}$ and $T_{rap}$. From the
iteration counts, we observe that $n_{max}$ does not change too much the preconditioner
quality, so it is suggested to keep its value as low as possible, as indicated by the
default parameters, in order to maintain the operator complexity low. The usage of
$\kappa_P = 10$ conducts to sparser prolongation operators and the worst \amgname{}
operator qualities regarding convergence; however, such drawback is compensated by the
lower cost of applying the multigrid operator and constructing it. On the other hand,
$\kappa_P = 100$ provides a much more accurate fit of the test vectors which naturally
leads to higher operator complexities. However, this improved accuracy is not well paid
since the iteration counts increase instead of decreasing when compared
to the default configuration in Table \ref{tab:sstvDefRes} due to larger weights introduces
in $P$. Therefore, we conclude that the default value selected for $\kappa_P$ represents a
good balance between cost and accuracy for running the prolongation construction phase.

\begin{table}[!htbp]
  \ra{1.4}
  \centering
  \caption{Sensitivity analysis with the prolongation parameters.}
  \label{tab:sstvProlongation}
  \resizebox{\textwidth}{!}{%
  \begin{tabular}{@{}*{1}{r}|*{1}{r}|*{2}{r}|*{12}{r}@{}}
  \toprule
  \multirow{2}{*}{\thead{\bf{Test} \\ \bf{case}}} & \multirow{2}{*}{\thead{\bf{Run} \\ \bf{ID}}} & \multicolumn{2}{c|}{\bf{Param.}} & \multicolumn{12}{c}{\bf{Results}} \\
  \cmidrule(lr){3-16}
           &   & $n_{max}$ & $\kappa_p$ & $C_{gd}$ & $C_{op}$ & $C_{fs}$ & $n_{it}$ & $T_{ts}$ & $T_{cs}$ & $T_{sm}$  & $T_{pl}$ & $T_{rap}$ &  $T_p$ & $T_s$ & $T_t$  \\
  \midrule
  {\bf T1} & 6 &  {\bf 2} &        50 &     1.39 &    3.15 &    0.53 &     229 &    18.2 &      9.4 &      6.9 &     3.8 &     2.6 &   42.5 &  23.9 &  66.9 \\
           & 7 & {\bf 10} &        50 &     1.39 &    3.79 &    0.53 &     236 &    18.9 &     11.3 &      7.6 &     5.4 &     3.3 &   48.0 &  27.3 &  75.6 \\
           & 8 &      10  &  {\bf 10} &     1.39 &    2.47 &    0.53 &     298 &    17.4 &      7.9 &      6.1 &     3.3 &     1.7 &   38.0 &  27.0 &  65.4 \\
           & 9 &      10  & {\bf 100} &     1.39 &    3.90 &    0.53 &     256 &    19.4 &     11.4 &      8.0 &     5.5 &     3.2 &   49.3 &  29.9 &  79.5 \\
  \midrule
  {\bf T2} & 6 &  {\bf 2} &        50 &     1.38 &    2.73 &    0.55 &     235 &    49.1 &     24.4 &     33.8 &     8.5 &     4.7 &  122.6 &  81.0 & 204.3 \\
           & 7 & {\bf 10} &        50 &     1.38 &    3.17 &    0.55 &     244 &    50.6 &     27.6 &     35.5 &    10.9 &     5.7 &  132.7 &  94.6 & 228.1 \\
           & 8 &      10  &  {\bf 10} &     1.38 &    1.95 &    0.55 &     419 &    45.8 &     19.9 &     30.7 &     6.4 &     2.3 &  107.0 & 114.8 & 222.4 \\
           & 9 &      10  & {\bf 100} &     1.38 &    3.23 &    0.55 &     256 &    52.0 &     28.0 &     35.5 &    11.7 &     5.9 &  135.6 & 102.3 & 238.6 \\
  \midrule
  {\bf T3} & 6 &  {\bf 2} &        50 &     1.48 &    2.49 &    0.56 &      82 &     9.5 &      4.2 &      3.8 &     2.0 &     1.6 &   22.3 &   4.3 &  26.9 \\
           & 7 & {\bf 10} &        50 &     1.48 &    2.66 &    0.56 &      85 &    10.5 &      4.5 &      4.1 &     2.3 &     1.7 &   24.9 &   4.9 &  30.2 \\
           & 8 &      10  &  {\bf 10} &     1.48 &    1.67 &    0.57 &     126 &     8.6 &      2.9 &      3.1 &     1.3 &     0.9 &   17.8 &   5.2 &  23.3 \\
           & 9 &      10  & {\bf 100} &     1.48 &    2.82 &    0.56 &      91 &     9.7 &      4.8 &      4.2 &     2.4 &     1.9 &   24.6 &   5.4 &  30.4 \\
  \bottomrule
  \end{tabular}}
\end{table}

{\bf Test space generation.}
In Table \ref{tab:sstvTSpace}, we show the sensitivity analysis of \amgname{} with respect
to its test space calculation phase. In such study, we vary the number of test vectors,
$n_{tv}$, and the number of outer iterations of the SRQCG algorithm, $n_{rq}$. Due to the
low computational cost, Ritz projection is applied at every iteration, i.e., $k_{ritz} =
1$. We start by noting that the test space calculation has almost no effect on the grid
and smoother complexities. On the other hand, we see for all test cases that the operator
complexity can be reduced when a more accurate test space is computed, e.g., by increasing
the value of $n_{rq}$. This behavior can be explained by the higher parallelism of
smoothest modes with respect to other eigenvectors, which normally leads to a faster
preconditioner regarding setup time. Looking at the times distribution, the test space
computation phase, represented by $T_{ts}$, is the most affected, as the cost of each
SRQCG iteration is directly proportional to $n_{tv}$.  We note that the use of $n_{rq} =
20$ leads to a slightly better preconditioner both in terms of convergence and total time
than $n_{rq} = 10$ as used by the default configuration. However, we prefer the default
choice in order to obtain better scalability. Little or no effect is seen in the other
phases with the only exception of the RAP computation of {\bf T2} which varies
significantly due to the changes in the operator complexity. Analyzing the iteration
counts, the usage of higher values of $n_{tv}$ and $n_{rq}$ usually lead to a
preconditioner that converges slightly faster which, however, is not counterbalanced by
the high increase in the setup time $T_p$.

\begin{table}[!htbp]
  \ra{1.4}
  \centering
  \caption{Sensitivity analysis with the test space parameters.}
  \label{tab:sstvTSpace}
  \resizebox{\textwidth}{!}{%
  \begin{tabular}{@{}*{1}{r}|*{1}{r}|*{2}{r}|*{12}{r}@{}}
  \toprule
  \multirow{2}{*}{\thead{\bf{Test} \\ \bf{case}}} & \multirow{2}{*}{\thead{\bf{Run} \\ \bf{ID}}} & \multicolumn{2}{c|}{\bf{Param.}} & \multicolumn{12}{c}{\bf{Results}} \\
  \cmidrule(lr){3-16}
           &    & $n_{tv}$  &  $n_{rq}$ & $C_{gd}$ & $C_{op}$ & $C_{fs}$ & $n_{it}$ & $T_{ts}$ & $T_{cs}$ & $T_{sm}$  & $T_{pl}$ & $T_{rap}$ & $T_p$ & $T_s$ & $T_t$  \\
  \midrule
  {\bf T1} & 10 & {\bf 20} &       10 &    1.37 &    3.09 &     0.53 &     223 &    43.7 &    9.1 &      6.6 &      3.7 &      2.7 &  67.1 &  22.9 &  90.3 \\
           & 11 & {\bf 30} &       10 &    1.37 &    2.99 &     0.53 &     235 &    79.9 &    9.5 &      6.6 &      4.0 &      2.5 & 103.7 &  22.7 & 126.7 \\
           & 12 & {\bf 10} & {\bf 20} &    1.39 &    2.74 &     0.53 &     185 &    22.5 &    8.7 &      6.6 &      3.2 &      1.9 &  44.2 &  17.6 &  62.1 \\
           & 13 &       10 & {\bf 30} &    1.39 &    2.51 &     0.53 &     182 &    26.0 &    8.0 &      6.2 &      2.8 &      2.0 &  46.4 &  16.6 &  63.2 \\
  \midrule
  {\bf T2} & 10 & {\bf 20} &       10 &    1.38 &    2.36 &     0.55 &     240 &    57.7 &   22.3 &     32.2 &      6.6 &      3.5 & 124.3 &  73.7 & 198.6 \\
           & 11 & {\bf 30} &       10 &    1.36 &    2.95 &     0.54 &     191 &   233.3 &   26.9 &     33.6 &     10.5 &      5.4 & 312.1 &  71.7 & 384.6 \\
           & 12 & {\bf 10} & {\bf 20} &    1.38 &    2.38 &     0.55 &     240 &    59.0 &   22.2 &     32.9 &      7.2 &      3.7 & 127.2 &  73.9 & 201.6 \\
           & 13 &       10 & {\bf 30} &    1.38 &    2.15 &     0.55 &     282 &    68.2 &   21.6 &     32.6 &      6.0 &      3.2 & 133.8 &  82.0 & 216.4 \\
  \midrule
  {\bf T3} & 10 & {\bf 20} &       10 &    1.49 &    2.80 &     0.57 &      76 &    22.9 &    4.6 &      4.0 &      2.4 &      2.1 &  37.0 &   4.2 &  41.5 \\
           & 11 & {\bf 30} &       10 &    1.49 &    2.77 &     0.57 &      74 &    40.7 &    4.6 &      3.7 &      2.4 &      1.8 &  54.2 &   3.9 &  58.3 \\
           & 12 & {\bf 10} & {\bf 20} &    1.48 &    2.24 &     0.57 &      77 &    11.5 &    3.9 &      3.7 &      1.6 &      1.3 &  23.0 &   3.8 &  27.1 \\
           & 13 &       10 & {\bf 30} &    1.48 &    2.15 &     0.57 &      85 &    13.1 &    3.5 &      3.3 &      1.4 &      1.0 &  23.1 &   3.7 &  27.0 \\
  \bottomrule
  \end{tabular}}
\end{table}

{\bf Smoother.}
In Table \ref{tab:sstvSmoother}, we show the sensitivity analysis of \amgname{} with
respect to its smoother calculation phase. Here, we consider only the standard input
parameters of the aFSAI smoother noting that, in all configurations, we guaranteed $\omega
> 0.8$ in order to achieve good convergence. Concerning iteration count, the
preconditioner quality is nearly the same when the quantity $k_g \times \rho_g$ is fixed,
as this product mainly controls the aFSAI density. However, it is preferable to set
$\rho_g$ greater than $k_g$ in order obtain a faster setup time for the smoother $T_{sm}$
and consequently for the preconditioner as well. The choice of $\epsilon_g = 10^{-3}$ over
$0.0$ does not provoke a substantial change in the preconditioner quality when $k_g \times
\rho_g < 32$, while it limits the FSAI complexity and, thus, its accuracy when such
inequality does not hold. Finally, we note that the default smoother configuration gives a
sufficiently good convergence when comparing the number of iterations reported in Tables
\ref{tab:sstvDefRes} and \ref{tab:sstvSmoother}. In fact, slightly better convergence can
be achieved by increasing $k_g$ and $\rho_g$ which are not paid off since the values of
$T_p$ and $T_s$ increase. Analyzing the complexities, we see that changing the smoother
configuration has practically no effect on the number of degrees of freedom per level of
the multigrid hierarchy, i.e., the grid complexity, while the operator complexity is
slightly reduced as a side-effect of computing more accurate smoothers. Naturally, the
FSAI complexity, $C_{fs}$, is the one more affected in this study. In addition, we note
that the optimal configurations in terms of total solution times are always associated
with the condition $C_{fs} < C_{op}$ meaning that the application of a V(1,1)-cycle
composed by the aFSAI smoother is less costly than the regular AMG cycle based on a
Gauss-Seidel or SOR smoother.

\begin{table}[!htbp]
  \ra{1.4}
  \centering
  \caption{Sensitivity analysis with the smoother parameters.}
  \label{tab:sstvSmoother}
  \resizebox{\textwidth}{!}{%
  \begin{tabular}{@{}*{1}{r}|*{1}{r}|*{3}{r}|*{12}{r}@{}}
  \toprule
  \multirow{2}{*}{\thead{\bf{Test} \\ \bf{case}}} & \multirow{2}{*}{\thead{\bf{Run} \\ \bf{ID}}} & \multicolumn{3}{c|}{\bf{Param.}} & \multicolumn{12}{c}{\bf{Results}} \\
  \cmidrule(lr){3-17}
           &    &   $k_g$ & $\rho_g$ &   $\epsilon_g$ & $C_{gd}$ &   $C_{op}$ &   $C_{fs}$ & $n_{it}$ & $T_{ts}$ & $T_{cs}$ & $T_{sm}$  & $T_{pl}$ & $T_{rap}$ & $T_p$ & $T_s$ & $T_t$ \\
  \midrule
  {\bf T1} & 14 &       4 & {\bf 8} &       $10^{-3}$ &     1.39 &      3.04 &      1.04 &     187 &     20.7 &    9.2 &     13.8 &     3.4 &     2.3 &   51.1 &  22.4 &  74.0 \\
           & 15 &       4 &      8  &    ${\bf 0.0}$ &     1.39 &      3.02 &      1.04 &     190 &     21.0 &    9.3 &     14.4 &     3.5 &     2.2 &   51.8 &  23.0 &  75.4 \\
           & 16 & {\bf 8} & {\bf 4} & ${\bf 10^{-3}}$ &     1.39 &      3.05 &      1.04 &     189 &     20.6 &    9.4 &     27.2 &     3.5 &     2.2 &   64.5 &  23.0 &  88.0 \\
           & 17 &       8 & {\bf 8} &       $10^{-3}$ &     1.39 &      2.87 &      2.04 &     149 &     25.2 &    9.0 &     56.7 &     3.3 &     2.1 &   98.2 &  24.4 & 123.1 \\
  \midrule
  {\bf T2} & 14 &       4 & {\bf 8} &       $10^{-3}$ &     1.38 &      2.63 &      1.07 &     207 &     57.1 &   23.7 &     53.3 &     7.7 &     4.0 &  148.6 &  91.1 & 240.4 \\
           & 15 &       4 &      8  &    ${\bf 0.0}$ &     1.38 &      2.62 &      1.07 &     196 &     55.3 &   24.1 &     53.2 &     7.6 &     4.0 &  147.6 &  86.7 & 235.1 \\
           & 16 & {\bf 8} & {\bf 4} & ${\bf 10^{-3}}$ &     1.38 &      2.61 &      1.06 &     197 &     55.0 &   24.2 &    102.1 &     7.5 &     4.3 &  195.7 &  87.1 & 283.6 \\
           & 17 &       8 & {\bf 8} &       $10^{-3}$ &     1.38 &      2.48 &      2.08 &     199 &     82.6 &   40.4 &    110.6 &     6.3 &     7.0 &  250.8 & 131.1 & 382.3 \\
  \midrule
  {\bf T3} & 14 &       4 & {\bf 8} &       $10^{-3}$ &     1.49 &      2.31 &      1.10 &      79 &     10.6 &    3.7 &      6.9 &     1.7 &     1.5 &   25.6 &   4.7 &  30.7 \\
           & 15 &       4 &      8  &    ${\bf 0.0}$ &     1.49 &      2.33 &      1.11 &      79 &     10.5 &    3.9 &      7.1 &     1.6 &     1.6 &   25.8 &   4.8 &  30.9 \\
           & 16 & {\bf 8} & {\bf 4} & ${\bf 10^{-3}}$ &     1.48 &      2.28 &      1.07 &      69 &     10.8 &    4.0 &     13.2 &     1.6 &     1.3 &   31.9 &   4.1 &  36.2 \\
           & 17 &       8 & {\bf 8} &       $10^{-3}$ &     1.48 &      2.24 &      2.10 &      66 &     12.6 &    3.6 &     26.8 &     1.5 &     1.2 &   46.5 &   4.6 &  51.3 \\
  \bottomrule
  \end{tabular}}
\end{table}

\begin{figure}[!htbp]
  \centering
  \includegraphics[width=\linewidth]{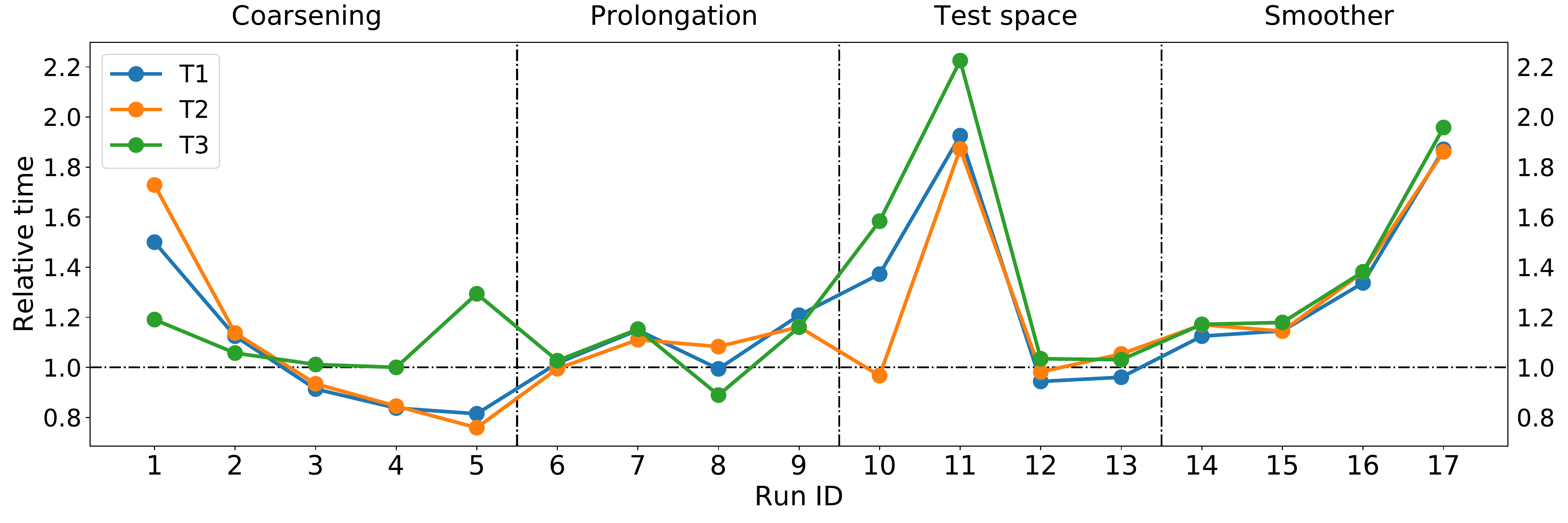}
  \caption{Distribution of the relative total time of different runs computed with respect
    to the default configuration.}
  \label{fig:sensitivity}
\end{figure}

In Figure \ref{fig:sensitivity}, we show a summary of the analysis developed here in terms
of relative total time. This measure is computed with respect to the default configuration
results from Table \ref{tab:sstvDefRes} to give an idea of the relative performance among
different runs and the default configuration for each test case. Also, we divide the
figure into four regions denoting which AMG setup phase is changed from run to run. From the
figure, we see that the test cases have a similar trend except for Run IDs $1$, $5$ and
$10$ where the selection of input parameters far from the default may impact the
performance positively or negatively, depending on the test case. Moreover, the maximum
performance degradation is equal to $2.2$ times, in case of the Run ID $11$, while most of
the other runs remain at a difference closer to $20 \%$ to the unitary baseline. This
indicates the robustness of \amgname{} setup in the range of parameters tested.

\subsection{Benchmark with state-of-the-art linear solvers}
\label{subsec:realworld}

The aim of this section is to compare the computational performance of the proposed
algorithm with other state-of-the-art AMG linear solver.  In Tables
\ref{tab:realWorldProblemsSolved1} and \ref{tab:realWorldProblemsSolved2}, we show the
grid and operator complexities ($C_{gd}$ and $C_{op}$); number of iterations ($n_{it}$);
and setup, solution and total times ($T_p$, $T_s$ and $T_t$) given by the BoomerAMG
(HYPRE) \cite{FalYan02}, GAMG (PETSc) \cite{PETSc} and \amgname{} preconditioners when
solving all the real-world test cases listed in Table \ref{tab:realWorldProblems}. The
setup parameters used for BoomerAMG are those listed as ``best practices'' in the work
\cite{FalSch14} with the exception of the strength threshold, which is set to $0.9$ here
as it provides slightly better results than $0.25$ as suggested in the aforementioned
communication. For GAMG, we use the standard input parameters from version 3.10 of the
PETSc library. Lastly, for \amgname{} we provide results for two different setup
configurations. In the first, \aspdef{}, we consider its default input parameters as
presented by Table \ref{tab:TestParam}. In the second, \aspopt{}, we use the optimal
configuration for each test case leading to the lowest total time.

\begin{table}[htbp]
\small
\ra{1.1}
\centering
\caption{Solution of real-world engineering problems (Part 1).}
\label{tab:realWorldProblemsSolved1}
\begin{tabular}{@{}*{1}{l}|*{1}{l}|*{3}{r}|*{3}{r}@{}}
\toprule
{\bf Test case}    & {\bf Preconditioner} &  $C_{gd}$ &   $C_{op}$ &   $n_{it}$ &  $T_p$ & $T_s$ &  $T_t$ \\
\midrule
{\tt heel1138k}    & \gamg{}              &     1.08 &      1.48 &        84 &   75.8 &  42.3 & 118.1 \\
                   & BoomerAMG            &     1.47 &      2.73 & $>\!1000$ &    --- &   --- &   --- \\ 
                   & \aspdef{}            &     1.39 &      3.29 &        82 &   31.2 &   7.9 &  39.1 \\ 
                   & \aspopt{}            &     1.17 &      1.98 &       103 &   18.5 &   7.5 &  26.0 \\ 
\midrule
{\tt guenda1446k}  & \gamg{}              &     1.10 &      1.70 & $>\!1000$ &    --- &   --- &   --- \\ 
                   & BoomerAMG            &     1.69 &      2.06 &       239 &   12.7 &  32.1 &  44.8 \\
                   & \aspdef{}            &     1.48 &      2.44 &        81 &   21.4 &   4.5 &  26.2 \\ 
                   & \aspopt{}            &     1.38 &      1.88 &       113 &   16.7 &   4.6 &  21.3 \\ 
\midrule
{\tt hook1498k}    & \gamg{}              &     1.11 &      1.92 &       128 &   13.6 &  16.2 &  29.8 \\
                   & BoomerAMG            &     1.61 &      2.02 &       154 &    4.4 &  16.5 &  20.9 \\
                   & \aspdef{}            &     1.37 &      2.72 &       151 &   22.1 &   8.3 &  30.4 \\ 
                   & \aspopt{}            &     1.20 &      1.81 &       279 &   12.3 &  11.3 &  23.4 \\ 
\midrule
{\tt tripod3239k}  & \gamg{}              &     1.09 &      1.56 &        77 &   17.2 &  20.7 &  38.2 \\
                   & BoomerAMG            &     1.52 &      3.04 &       181 &   13.3 &  53.1 &  66.4 \\
                   & \aspdef{}            &     1.39 &      3.08 &       228 &   42.4 &  23.4 &  65.8 \\ 
                   & \aspopt{}            &     1.26 &      1.94 &       172 &   47.7 &  15.2 &  62.9 \\ 
\midrule
{\tt thdr3559k}    & \gamg{}              &     1.09 &      1.62 &       933 &   23.2 & 237.4 & 260.6 \\
                   & BoomerAMG            &     1.71 &      2.03 &       176 &   53.1 &  58.1 & 111.2 \\
                   & \aspdef{}            &     1.46 &      2.34 &       416 &   62.7 &  62.9 & 125.6 \\ 
                   & \aspopt{}            &     1.70 &      2.18 &       136 &   71.7 &  49.3 & 121.6 \\ 
\midrule
{\tt beam6502k}    & \gamg{}              &     1.05 &      1.16 &        42 &  141.1 & 165.9 & 307.0 \\
                   & BoomerAMG            &     1.26 &      1.35 &        81 &   98.8 & 198.2 & 297.0 \\
                   & \aspdef{}            &     1.41 &      2.36 &        87 &  115.1 &  53.6 & 168.7 \\ 
                   & \aspopt{}            &     1.19 &      1.46 &        83 &   73.8 &  28.0 & 101.8 \\ 
\bottomrule
\end{tabular}
\end{table}
\begin{table}[htbp]
\small
\ra{1.1}
\centering
\caption{Solution of real-world engineering problems (Part 2).}
\label{tab:realWorldProblemsSolved2}
\begin{tabular}{@{}*{1}{l}|*{1}{l}|*{3}{r}|*{3}{r}@{}}
\toprule
{\bf Test case}    & {\bf Preconditioner} &  $C_{gd}$ &  $C_{op}$ &   $n_{it}$ &  $T_p$ &  $T_s$ &  $T_t$ \\
\midrule
{\tt jpipe6557k}   & \gamg{}              &     1.09 &     1.45 &       262 &   40.1 &  195.4 &  235.5 \\
                   & BoomerAMG            &     1.53 &     2.61 &       358 &   30.9 &  261.1 &  292.0 \\
                   & \aspdef{}            &     1.41 &     2.42 &       518 &   95.4 &  141.5 &  236.9 \\ 
                   & \aspopt{}            &     1.20 &     1.55 &       407 &   82.3 &   69.2 &  151.5 \\ 
\midrule
{\tt gear8302k}    & \gamg{}              &     1.03 &     1.17 & $>\!1000$ &    --- &    --- &    --- \\ 
                   & BoomerAMG            &     1.61 &     2.09 & $>\!1000$ &    --- &    --- &    --- \\ 
                   & \aspdef{}            &     1.38 &     2.04 &       652 &  412.6 &  586.9 &  999.5 \\ 
                   & \aspopt{}            &     1.23 &     1.65 &       628 &  374.0 &  521.1 &  895.1 \\ 
\midrule
{\tt eiffel9213k}  & \gamg{}              &     1.09 &     1.64 &       167 &   55.1 &  131.1 &  186.2 \\
                   & BoomerAMG            &     1.49 &     2.85 &       376 &   74.6 &  316.2 &  390.8 \\
                   & \aspdef{}            &     1.38 &     2.75 &       239 &  123.0 &   82.4 &  205.4 \\ 
                   & \aspopt{}            &     1.17 &     1.72 &       214 &   86.1 &   49.9 &  136.0 \\ 
\midrule
{\tt guenda11m}    & \gamg{}              &     1.08 &     1.59 & $>\!1000$ &    --- &    --- &    --- \\
                   & BoomerAMG            &     1.70 &     2.22 & $>\!1000$ &    --- &    --- &    --- \\ 
                   & \aspdef{}            &     1.49 &     2.27 &       374 &  596.6 &  471.9 & 1068.5 \\ 
                   & \aspopt{}            &     1.44 &     2.02 &       436 &  424.1 &  371.4 &  795.5 \\ 
\midrule
{\tt wrench13m}    & \gamg{}              &     1.08 &     1.52 &       139 &   41.8 &  162.2 &  204.0 \\
                   & BoomerAMG            &     1.51 &     2.98 &       183 &   40.5 &  227.9 &  268.4 \\
                   & \aspdef{}            &     1.39 &     3.12 &       326 &  183.8 &  215.2 &  399.0 \\ 
                   & \aspopt{}            &     1.17 &     1.70 &       211 &  146.2 &   82.6 &  228.8 \\ 
\midrule
{\tt agg14m}       & \gamg{}              &     1.08 &     1.59 &        30 &   82.1 &   37.8 &  119.9 \\
                   & BoomerAMG            &     1.55 &     1.84 &        89 &   61.9 &   88.7 &  150.6 \\
                   & \aspdef{}            &     1.42 &     3.73 &        29 &  310.8 &   37.7 &  348.5 \\ 
                   & \aspopt{}            &     1.18 &     1.70 &        64 &   93.8 &   24.5 &  118.3 \\ 
\midrule
{\tt M20}          & \gamg{}              &     1.03 &     1.16 &       113 & 1274.7 & 2046.2 & 3320.9 \\
                   & BoomerAMG            &     1.54 &     1.58 & $>\!1000$ &    --- &    --- &    --- \\ 
                   & \aspdef{}            &     1.38 &     2.26 &       401 &  473.6 & 1225.8 & 1699.4 \\ 
                   & \aspopt{}            &     1.14 &     1.44 &       350 &  353.2 &  779.4 & 1132.6 \\ 
\bottomrule
\end{tabular}
\end{table}

Looking at the results, we note that \amgname{} leads to an efficient solution method in
most test cases. Analyzing total time, we obtained a speedup of nearly five times in the
best scenario ({\tt heel1138k}); while, in the worst case ({\tt tripod3239k}), we were
less than two times slower than the fastest approach. It is worth noting that the
\amgname{} version configured with default input parameters proves to be more efficient
than BoomerAMG and GAMG in $8$ out of $13$ test cases, while when configured with the
optimal parameters, \aspopt{} becomes the winner in $11$ problems. Moreover, it is
particularly the best strategy in all problems above $6M$ DOFs, ensuring the method's
excellence in the solution of large-size systems. Lastly, we observe that our adaptive AMG
was able to solve $5$ problems in which at least one of the other preconditioners could
not solve within the maximum number of iterations requested (1000), being $2$ of those not
solvable by either of the other methods.

As regards complexities, we can note that the lowest values, as one can expect for
aggregation-based methods, are always provided by GAMG. BoomerAMG and \amgname{} with
default parameters provide similar complexities, while the optimized runs of \amgname{}
show an intermediate behavior. One of the main factor increasing the operator complexity
is the prolongation density. On average, BoomerAMG has a higher grid complexity with a
lower operator complexity, with respect to the default \amgname{}. This means that our
prolongation is denser and, usually, this provides a better result, allowing to solve the
linear system with fewer iterations, e.g., see {\tt guenda1446k} and {\tt
  tripod3239k}. For problems well-tackled by GAMG, the aggregation-based method is able to
solve the system with the lowest number of iterations, see {\tt beam6502k} and {\tt
  agg14m}, but when it is not able to deal with the problem, the iteration count increases
a lot, see {\tt guenda1446k} and {\tt thdr3559k}. Generally, the number of iterations
among the two runs of \amgname{}, default and optimized, are comparable. This means that
the optimized version can reduce the total time building a cheaper preconditioner
with approximately the same accuracy of the one obtained with the default parameters. When
BoomerAMG solves the problem, sometimes the number of iterations is similar to the one
obtained with the default \amgname{}, see {\tt hook1498k} and {\tt beam6502k}, while in
other cases it is larger or smaller. This means that the two approaches, on average, have
the same effectiveness. Finally, \aspopt{} solves 9 out of 13 test cases with fewer
iterations with respect to BoomarAMG, showing to be more well-suited for structural
problems, like an aggregation-based AMG.

\begin{figure}[!htbp]
  \centering
  \includegraphics[width=0.95\linewidth]{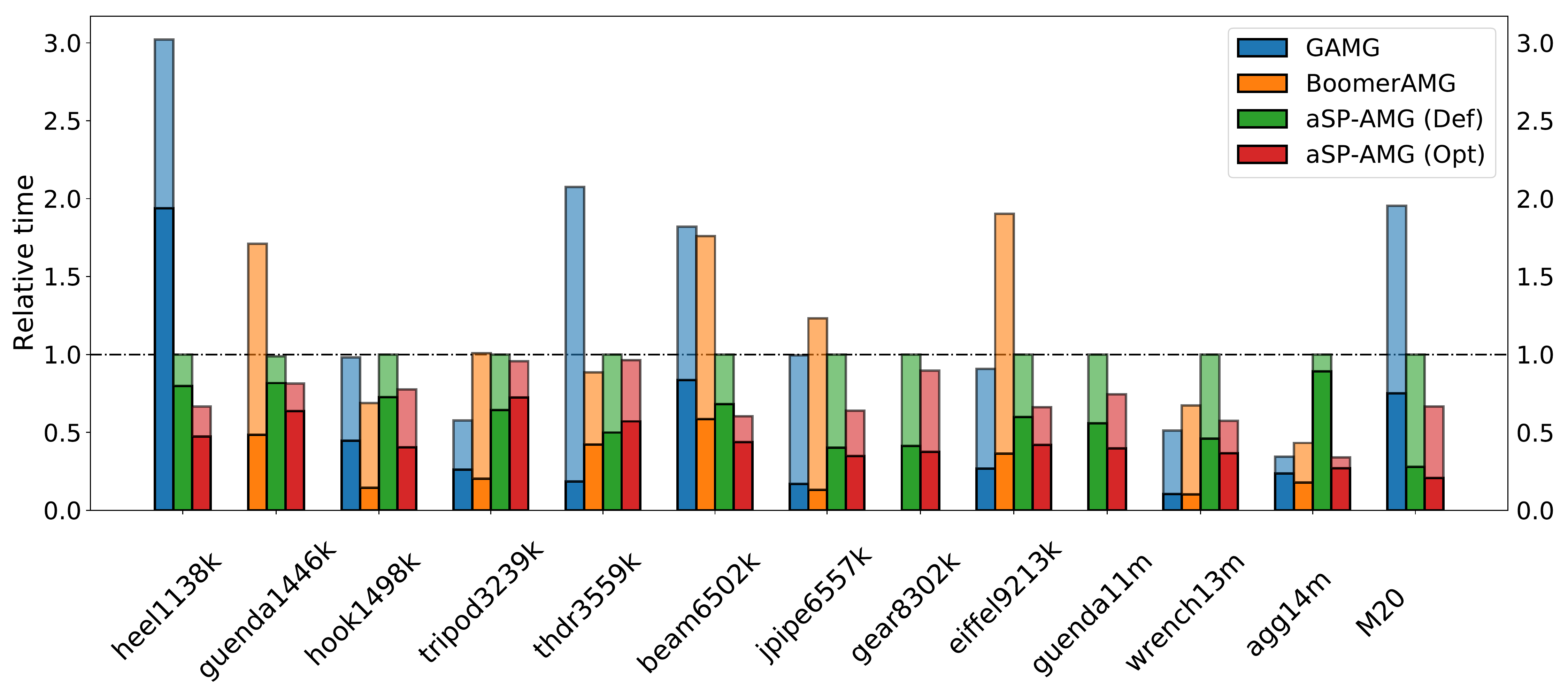}
  \caption{Relative total time for different preconditioners with respect to the \aspdef{}
    results. Darker portions of the bars represent preconditioner setup, while the lighter
    segments represent the time spent in the solver application phase.}
  \label{fig:realWorldRelTime}
\end{figure}

Figure \ref{fig:realWorldRelTime} shows a summary of the performance comparison presented
in Tables \ref{tab:realWorldProblemsSolved1} and \ref{tab:realWorldProblemsSolved2} in
terms of the relative total time for each strategy computed with respect to
\aspdef{}. Darker portions of the bars represent preconditioner setup, while the lighter
segments represent the time spent in the solver application phase. In general, one can
note that \amgname{} spends more time in its setup phase than the other preconditioners,
which is mainly caused by the smoother and test space computation that is not present in
the other strategies. However, such a difference gets balanced in the solver application
phase by many aspects such as the lower cost of the aFSAI smoother application with
respect to Gauss-Seidel and a coarse-grid correction that captures better slow-to-converge
modes of the smoother.

\subsection{Strong scalability}
\label{subsec:strongScal}

\begin{figure}[t!]
  \centering
  \subfloat[\label{fig:tripod3239k_strongScal}]{%
    \includegraphics[width=0.45\linewidth]{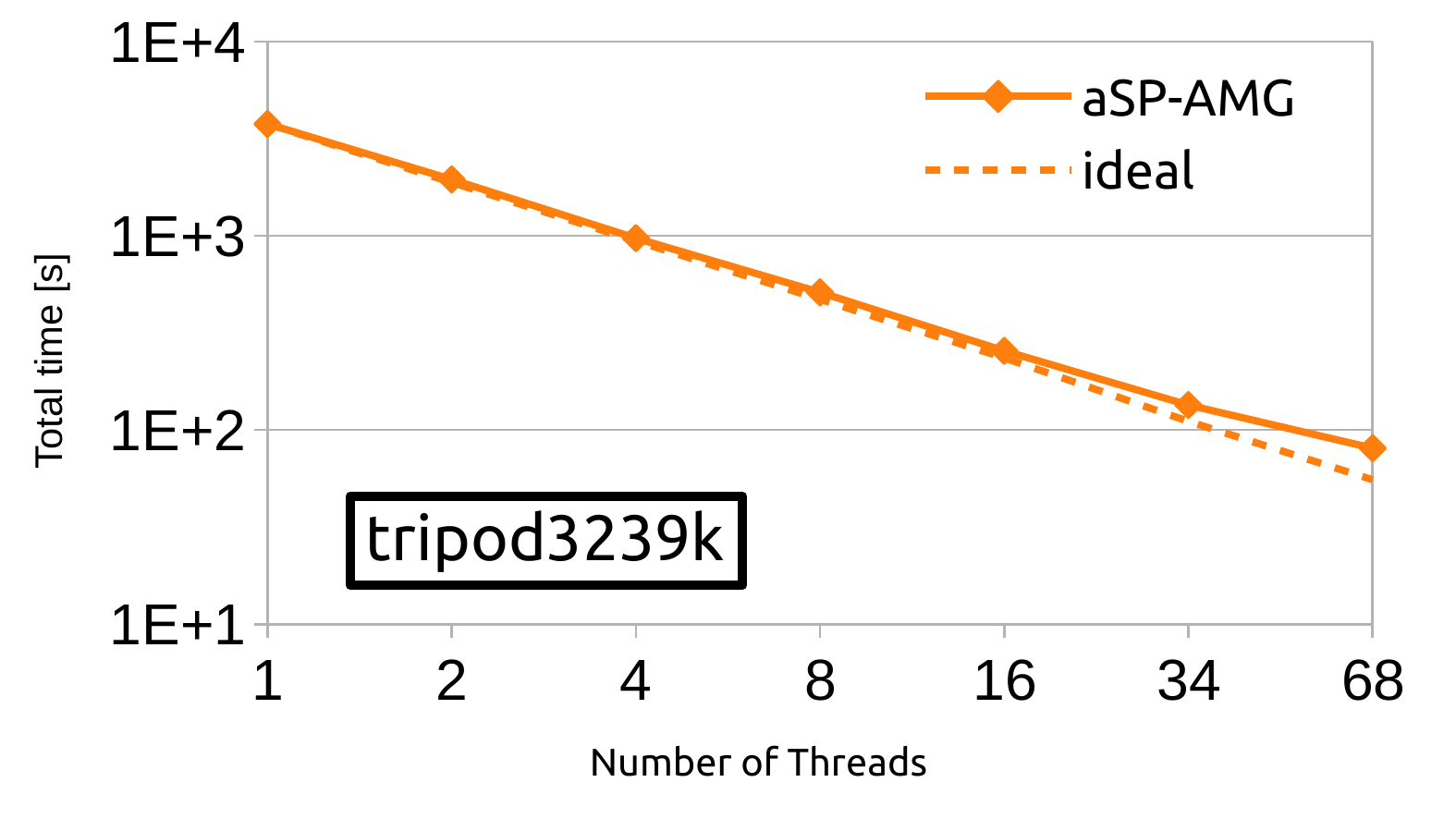}
  }
  \subfloat[\label{fig:jpipe6557k_strongScal}]{%
    \includegraphics[width=0.45\linewidth]{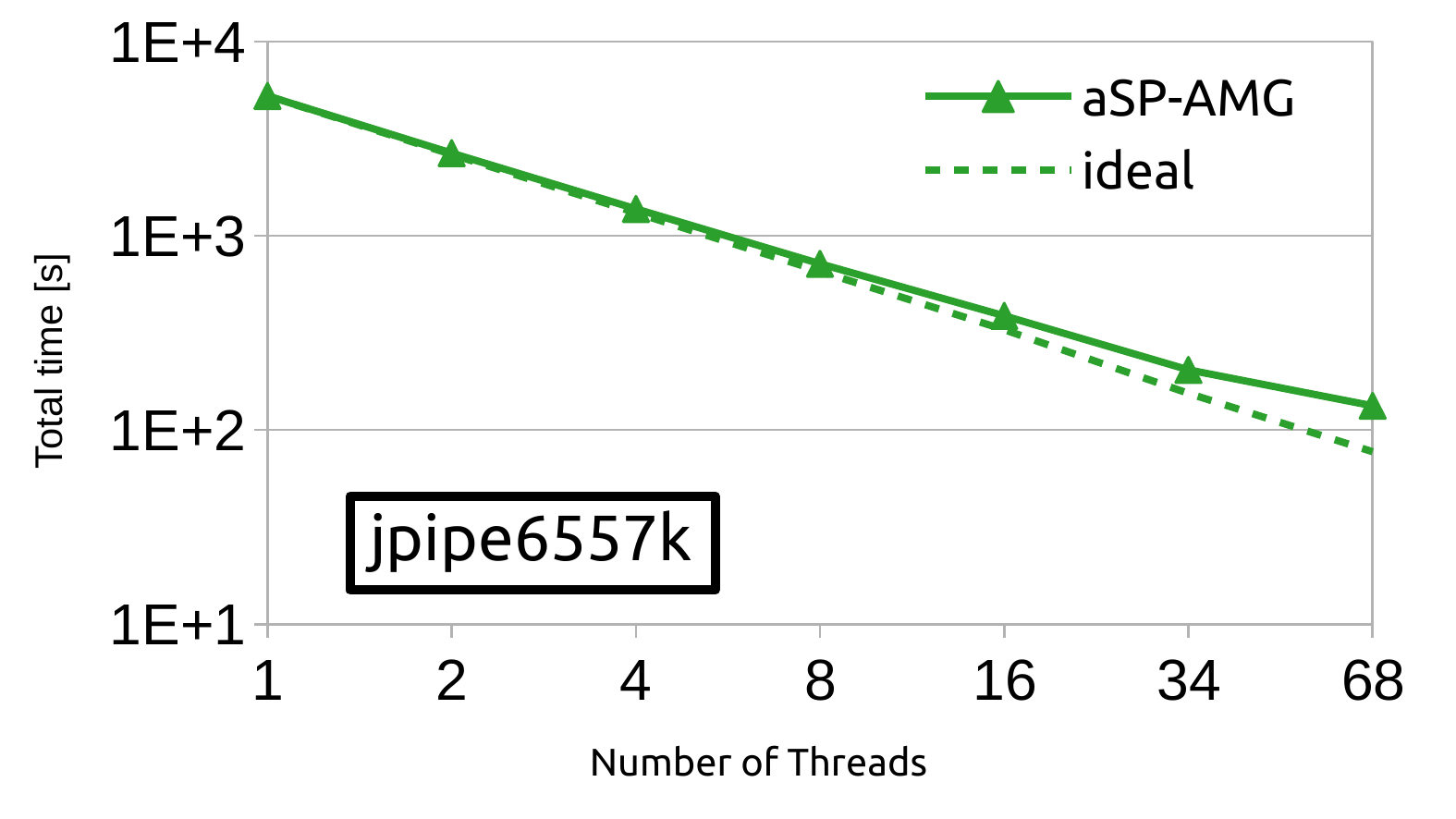}
  }
  \\
  \subfloat[\label{fig:eiffel9213k_strongScal}]{%
    \includegraphics[width=0.45\linewidth]{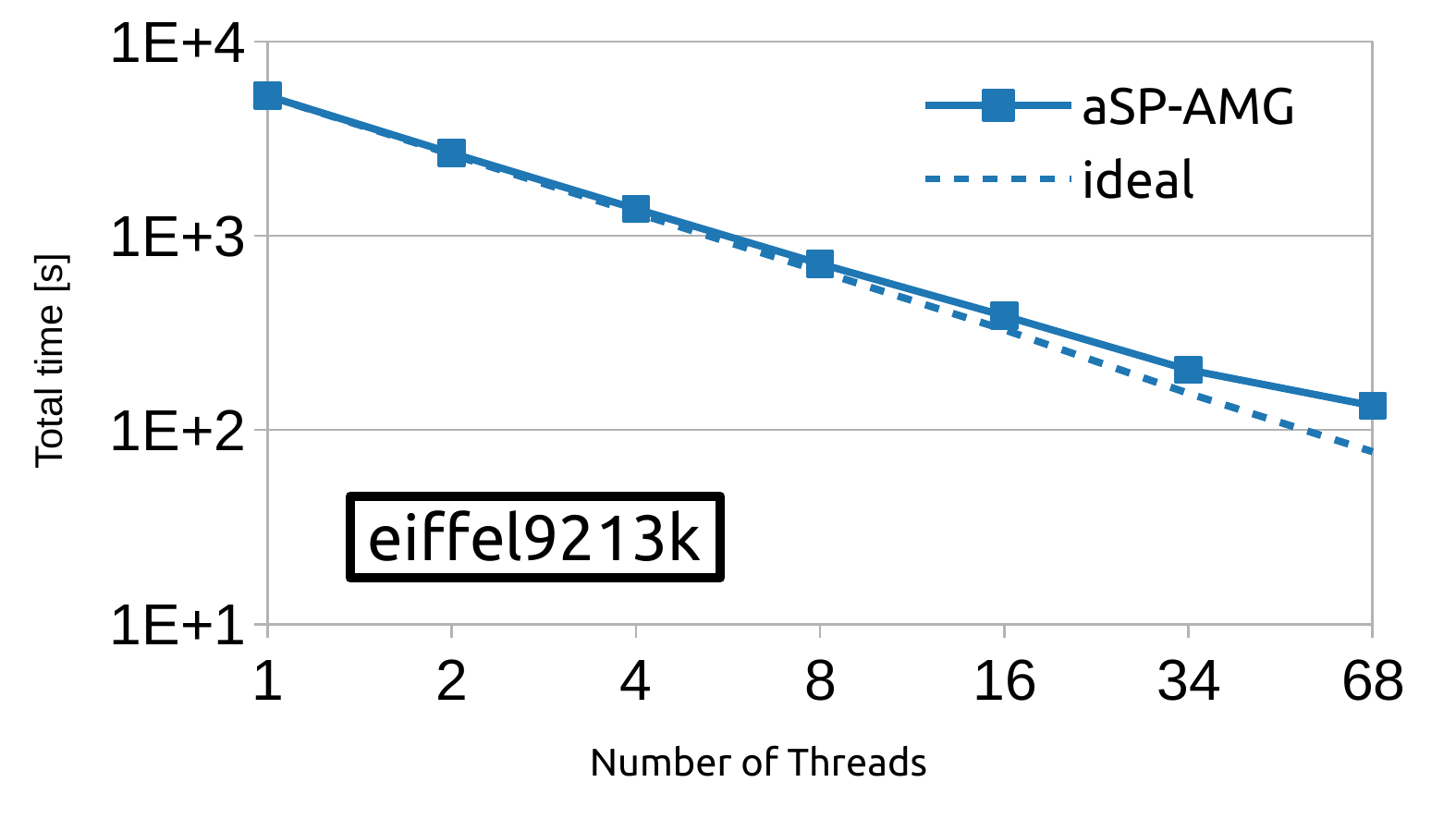}
  }
  \subfloat[\label{fig:parallelEfficiency}]{%
    \includegraphics[width=0.45\linewidth]{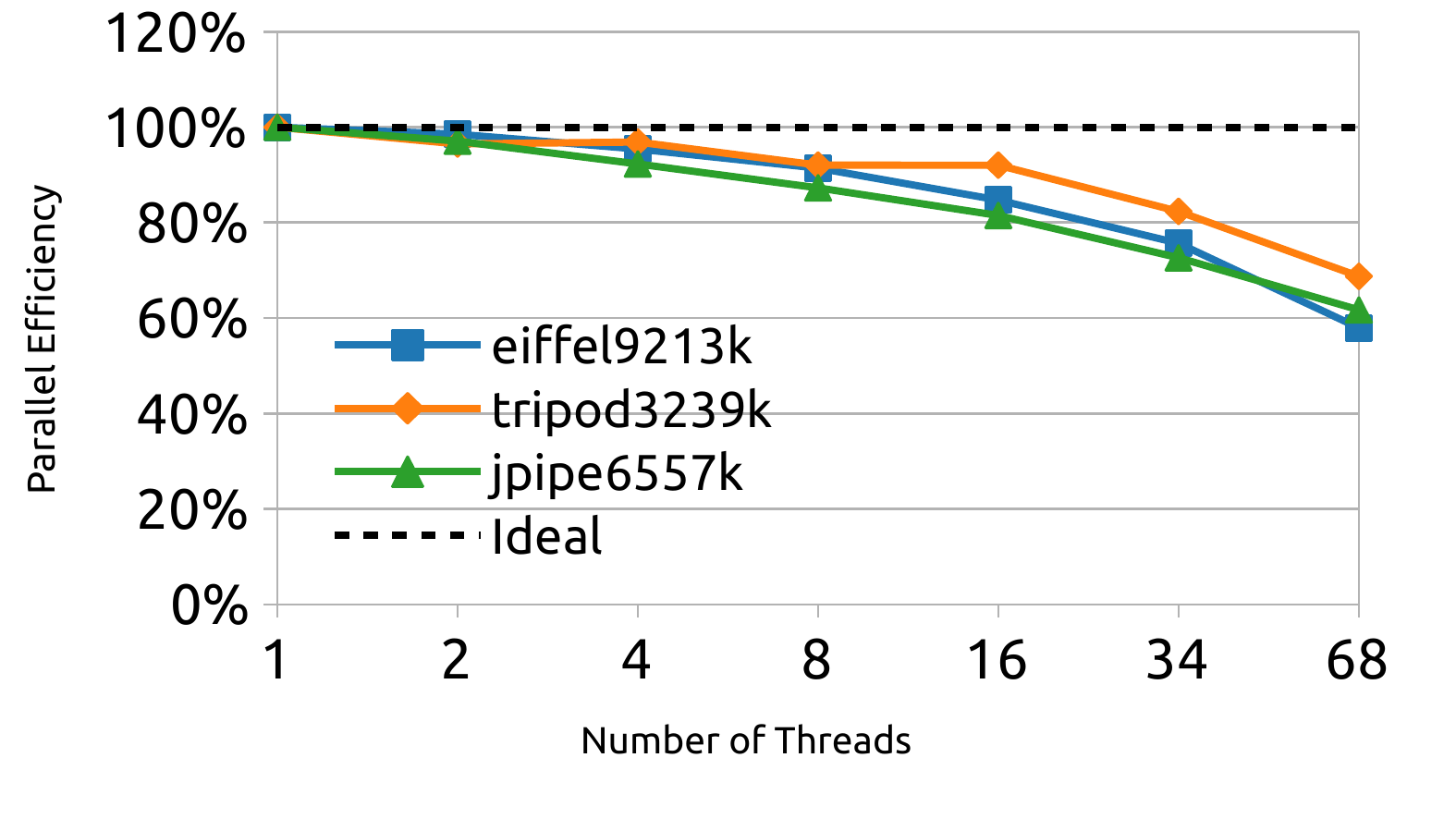}
  }
  \caption{(a-c): total solution times given by \amgname{} for test cases {\tt
      tripod3239k}, {\tt jpipe6557k}, and {\tt eiffel9213k}, respectively, when varying
    the number of computing cores from $1$ to $68$. (d): parallel efficiency curves for
    the previous test cases.}
  \label{fig:strongScal}
\end{figure}

Here, we evaluate strong scalability of our OpenMP implementation of \amgname{} in a
single node of the MARCONI-A2 cluster. For this, we consider three real-world test cases
from Table \ref{tab:realWorldProblems} that are representative of different problem sizes:
{\tt tripod3239k}, {\tt jpipe6557k} and {\tt eiffel9213k}. Figures
\ref{fig:tripod3239k_strongScal} to \ref{fig:eiffel9213k_strongScal} show total time for
solving the respective linear systems by using from $1$ up to $68$ computing cores, and,
for completeness, the ideal scaling profiles are depicted with the dashed line. Lastly,
Figure \ref{fig:parallelEfficiency} shows parallel efficiency, defined as the ratio
between ideal and wallclock time, associated with the last numerical results. As expected,
total time is reduced in all test cases with the increase of computational resources, and
this occurs with an efficiency up to $80 \%$ with respect to the ideal case when using
less than $32$ computing cores and over $60 \%$ when the total number of physical cores
are used, i.e., $68$. Such behavior fairly close to the ideal efficiency ($100 \%$) is due
to memory bandwidth limitations.

\section{Conclusions}
\label{sec:conclusions}

In this work, we introduced further modifications in the adaptive Smoothing and
Prolongation based Algebraic Multigrid method (initially proposed by \cite{MagFraJan19})
in order to improve the performance and usability of such method in the solution of
large-scale and challenging SPD linear systems arising from the discretization of linear
elasticity PDEs.

Initially, a sensitivity analysis is carried out to uncover the most important
configuration parameters for \amgname{} as well as their suitable range of usability. From
this study, we show that a large part of these parameters can be set to a default value,
provided by the analysis, obtaining performances that are comparable to the optimized
ones. After this, the \amgname{} is compared to the other multigrid methods, such as GAMG
and BoomerAMG, in the solution of real-world structural problems showing that \amgname{}
leads to the faster solution method in most of the cases in terms of both iteration time
as well as total execution time.

The next steps of the present research will be the extension to the block version, i.e.,
grouping together $x$, $y$ and $z$ unknowns of each physical node, to exploit the
supernodal nature of structural matrices. Other developments concern further improvements
in the numerical implementation as well as the extension to the hybrid MPI/OpenMP
programming paradigm.

\appendix
\section{Glimpses of meshes employed by the real-world engineering problems}
\label{sec:meshGrimpse}

This appendix provides a detailed description of the test cases presented in this work, 
as listed in Table \ref{tab:realWorldProblems}.

\subsection{Test case {\tt heel1138k}}
\label{sec:heel1138k}
The mesh derives from a 3D mechanical equilibrium problem of a human heel composed by four
physical regions. The domain is well-constrained with its top surface totally
clamped. Discretization is done via $2,247,515$ linear tetrahedral finite elements and
$379,481$ vertices. For further details the reader should refer to \cite{FonMatCarWih12,
BagFraSpiJan17}.

\subsection{Test case {\tt guenda1446k} and {\tt guenda11m}  }
\label{sec:guenda1446k}
\begin{figure}[!htbp]
  \centering
  \includegraphics[width=0.45\linewidth]{guenda_4_noRS.pdf}
  \includegraphics[width=0.45\linewidth]{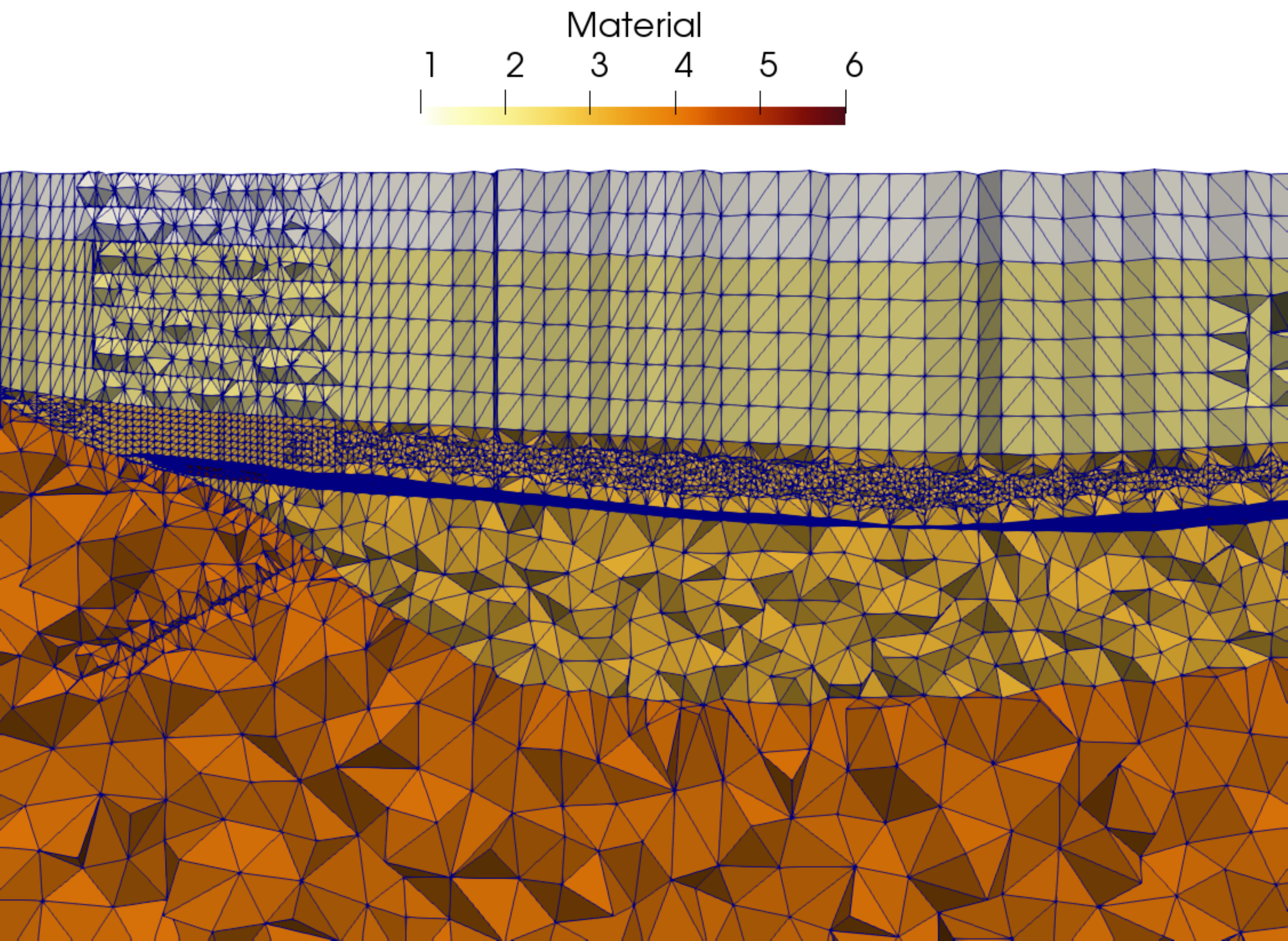}
  \caption{Left: global view of the geomechanical model solved by the {\tt guenda1446k}
  test case, the gas reservoir is in the central position where the mesh becomes more
  refined. Right: zoom in a planar cut of the mesh passing through the barycenter. Note
  that the elements become much more distorted in this region because of the high aspect
  ratio between the horizontal plane and vertical direction. The test case {\tt guenda11m}
  is obtained from running one-step of uniform mesh refinement from case {\tt
  guenda1446k}.}
  \label{fig:guenda1446k_mesh}
\end{figure}

The matrix {\tt guenda1446k} derives from a 3D geomechanical model of an underground gas storage (UGS)
site. The domain spans an area of $40 \times 40 \, km^{2} $ and extends down to $5 \,
km$ depth. To reproduce with high fidelity the real geometry of the gas reservoir, a
severely distorted mesh with $2,833,237$ linear tetrahedral elements and $482,208$
vertices is used. While fixed boundaries are prescribed on the bottom and lateral
sides, the surface is traction-free. The matrix {\tt guenda11m} is a 
refined version of the test case {\tt guenda1446k}. The highly
irregular geometry requires a distorted mesh with $22,665,896$ linear tetrahedra and
$3,817,466$ vertices. In Figure \ref{fig:guenda1446k_mesh}, we show a
representation of the problem's geometry and mesh.

\subsection{Test case {\tt hook1498k}}
\label{sec:hook1498k}

The mesh derives from a 3D mechanical equilibrium problem of a poorly constrained steel
hook discetized with $2,578,916$ linear tetrahedral finite elements giving rise to
$499,341$ vertices. For further details the reader should refer to \cite{BagFraSpiJan17}.

\subsection{Test case {\tt tripod3239k}}
\label{sec:tripod3239k}
\begin{figure}[!htbp]
  \centering
  \includegraphics[width=0.45\linewidth]{tripod_1_noRS.pdf}
  \includegraphics[width=0.45\linewidth]{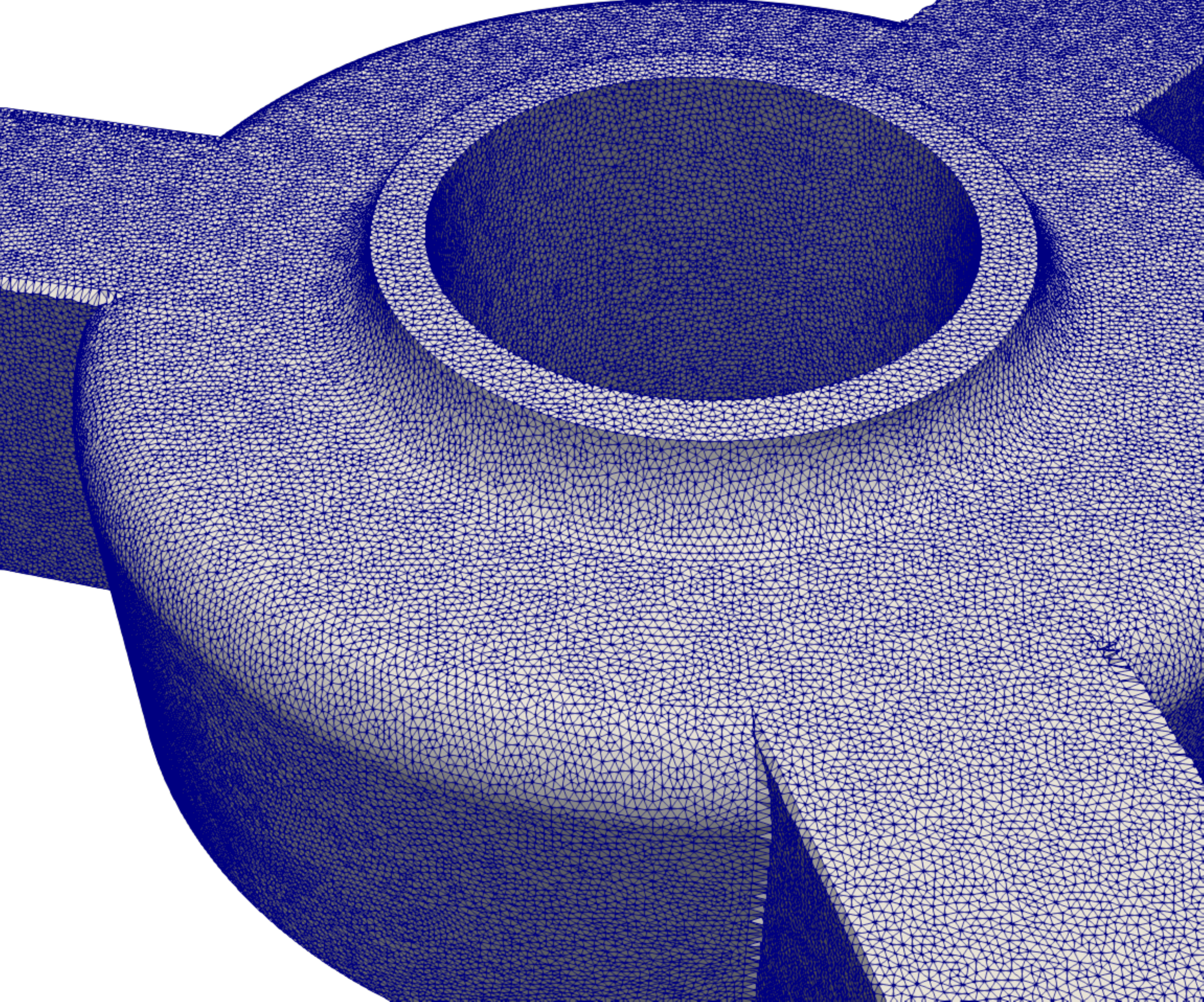}
  \caption{Left: geometry referent to the test case {\tt tripod3239k} ({\bf T1}). Right:
  zoom in the mesh of the top surface of the tripod. The elements used herein are
  regular in general, thus do not contribute heavily to the ill-conditioning of the
  system matrix.}
  \label{fig:tripod3239k_mesh}
\end{figure}

The mesh derives from a 3D mechanical equilibrium of a tripod with
clamped bases. Material is linear elastic with $(E,\nu) = (10^{6} MPa, 0.45)$. Mesh is
formed by $7,523,450$ linear tetrahedra and discretization is given by the finite
element method. A linear problem with $3,239,649$ DOFs is obtained and the system
matrix contains $142,351,857$ entries. Figure \ref{fig:tripod3239k_mesh} shows the
geometry and the mesh of the problem.

\subsection{Test case {\tt thdr3559k}}
\label{sec:thdr3559k}
\begin{figure}[!htbp]
  \centering
  \includegraphics[width=0.45\linewidth]{teod_3.pdf}
  \includegraphics[width=0.45\linewidth]{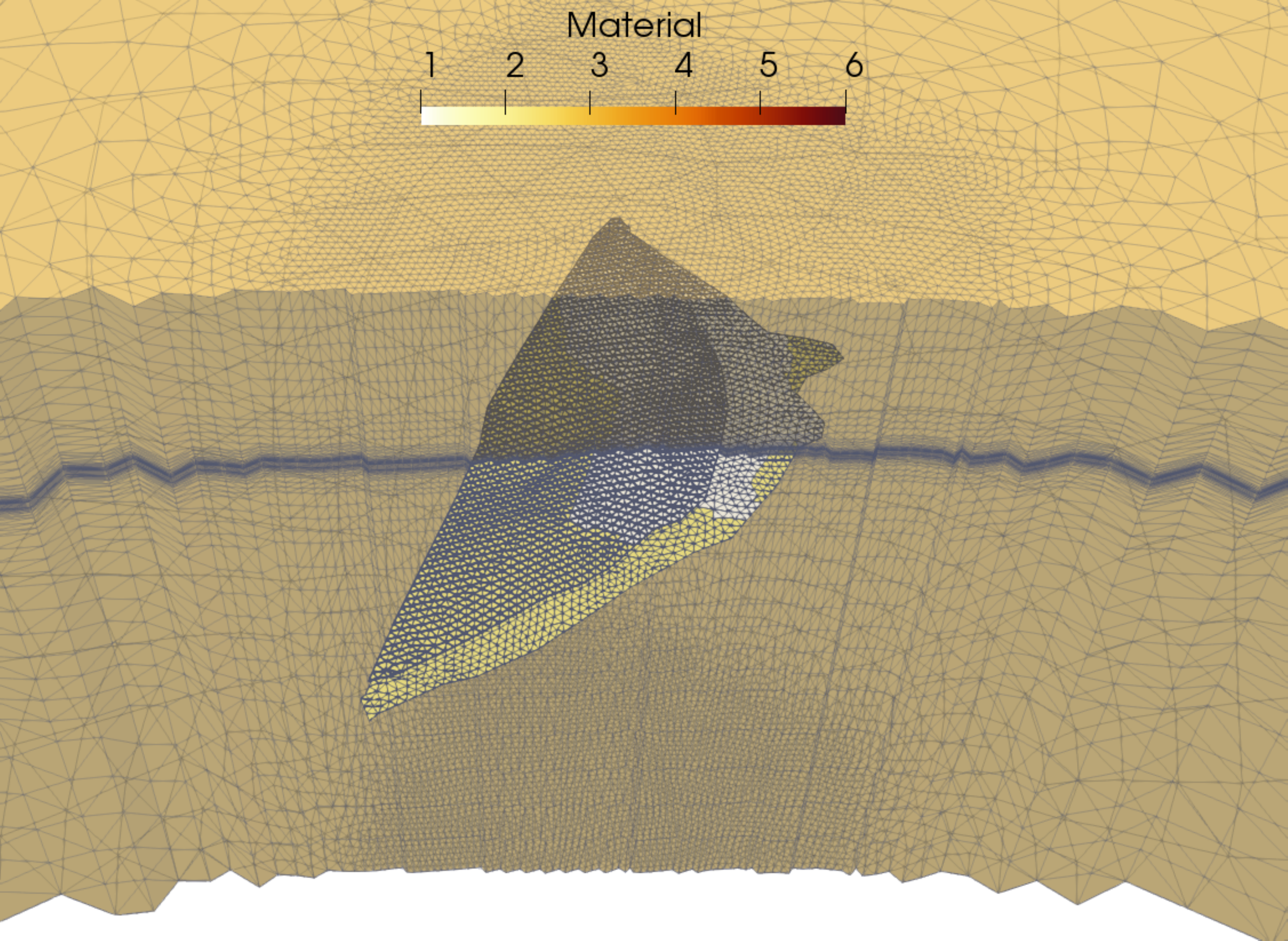}
  \caption{Left: global view of the mesh referent to test case {\tt thdr3559k}. The
  geometry is a simple box that contains a gas reservoir located in its central position,
  where the mesh becomes more refined. Right: zoom in the reservoir region showing the
  distorted elements at the middle position.}
  \label{fig:thdr3559k_mesh}
\end{figure}

The matrix derives froma a 3D geomechanical model of a vertically extruded
gas-reservoir. Geometry is a simple box with dimensions $40 \times 40 \times 5 \,
km^{2}$ composed by multiple materials with mechanical properties of the medium
varying in depth. A unstructured mesh is formed by $7,014,887$ linear tetrahedra and
$1,186,466$ vertices. The resulting stiffness matrix contains $3,559,398$ DOFs and
$81,240,330$ entries. In Figure \ref{fig:thdr3559k_mesh}, we show a representation of
the problem's geometry and mesh.

\subsection{Test case {\tt beam6502k}}
\label{sec:beam6502k}
The matrix arises from a mechanical equilibrium of a 3D multi-material cantilever beam
with its left face clamped. The body measures $1 \times 1 \times 8 \, m^3$, and it is
composed by two materials: in the first half of the beam, we have $(E_1,\nu_1) =
(10^{4} MPa, 0.2)$, while in the second, $(E_2,\nu_2) = (10^{7} MPa, 0.4)$. A structured mesh
is formed by $2,097,152$ linear hexahedra elements that are regularly shaped and
$2,167,425$ vertices. The resulting stiffness matrix contains $6,502,275$ DOFs. Such
problem is adapted from the {\tt ex2p} example driver from the MFEM library \cite{MFEM},
where the reader may find more details about this test case.

\subsection{Test case {\tt jpipe6557k}}
\label{sec:jpipe6557k}
The mesh derives from a mechanical equilibrium of a 3D structure composed
of three orthogonal pipes connected via special joints. The whole body is contained in
a box of $273 \times 339 \times 478 \, cm^3$. Mesh is unstructured and made of
$11,218,064$ linear tetrahedra elements and $2,264,181$ vertices. The characteristics of
the material are $(E,\nu) = (210,000 MPa, 0.2)$. The resulting
stiffness matrix contains $6,557,808$ DOFs and $335,451,702$ entries. Figure
\ref{fig:thdr3559k_mesh} shows the problem's geometry and mesh.
\begin{figure}[!htbp]
  \centering
  \includegraphics[width=0.45\linewidth]{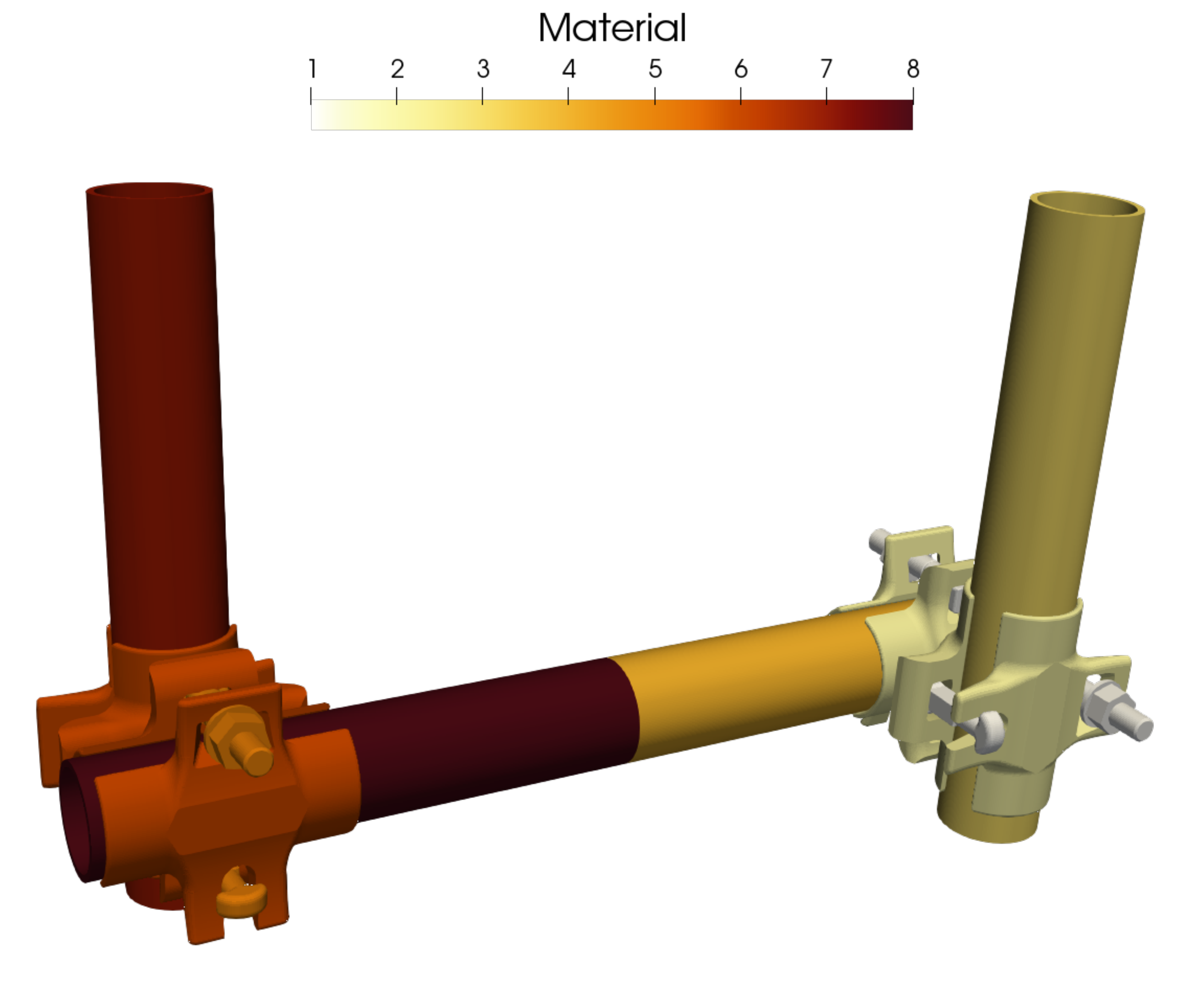}
  \includegraphics[width=0.45\linewidth]{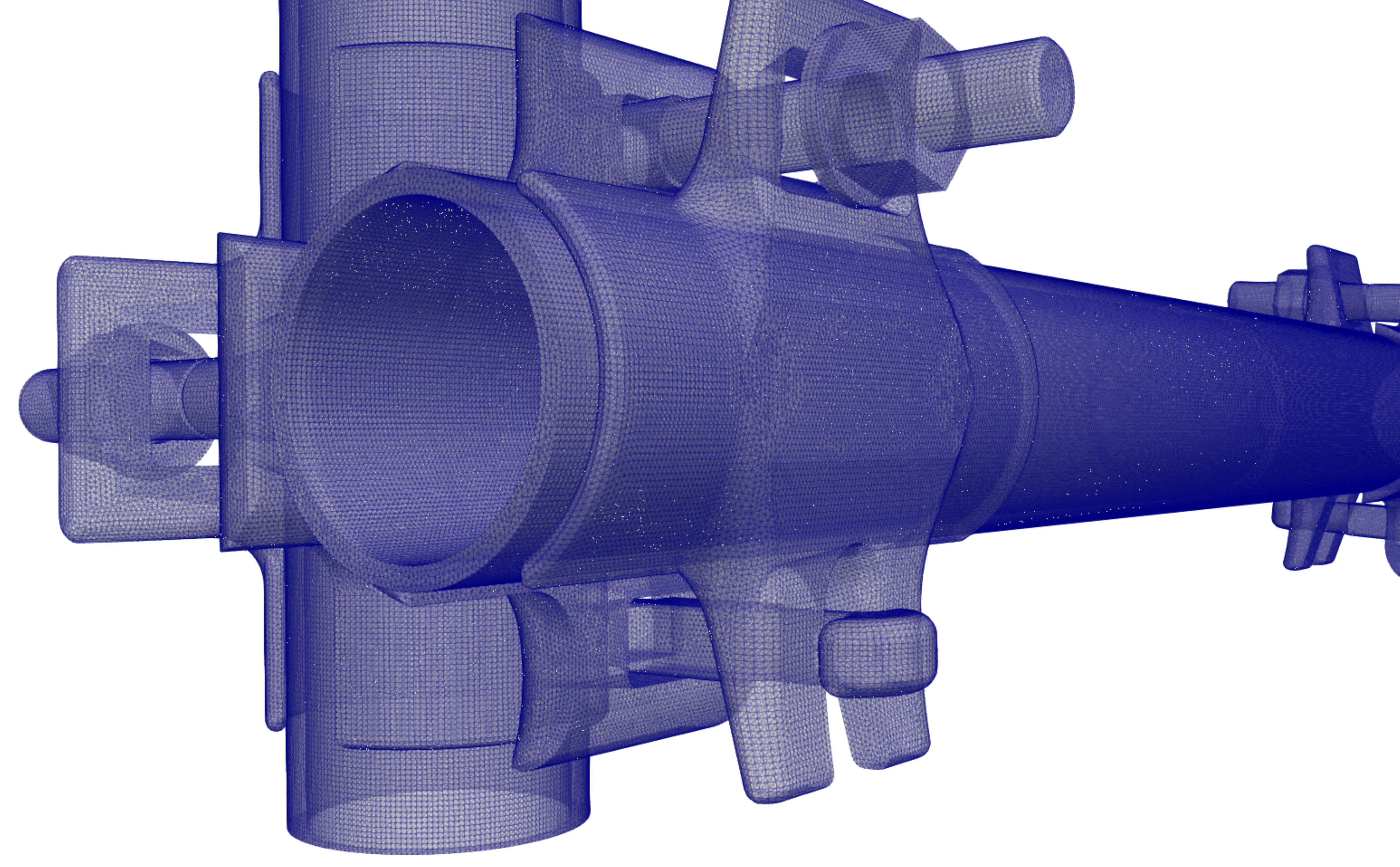}
  \caption{Left: geometry referrent to the test case {\tt jpipe6557k}. The colors
  represent different materials. Right: zoom in the left pipes connection. Although this
  geometry shows more details such as the presence of fillets, bolts and nuts, the
  elements used to represent it are well-shaped as showed in
  Figure \ref{fig:meshQuality}.}  
  \label{fig:jpipe6557k_mesh}
\end{figure}

\subsection{Test case {\tt gear8302k}}
\label{sec:gear8302k}
\begin{figure}[!htbp]
  \centering
  \includegraphics[width=0.45\linewidth]{gear_2_noRS.pdf}
  \includegraphics[width=0.45\linewidth]{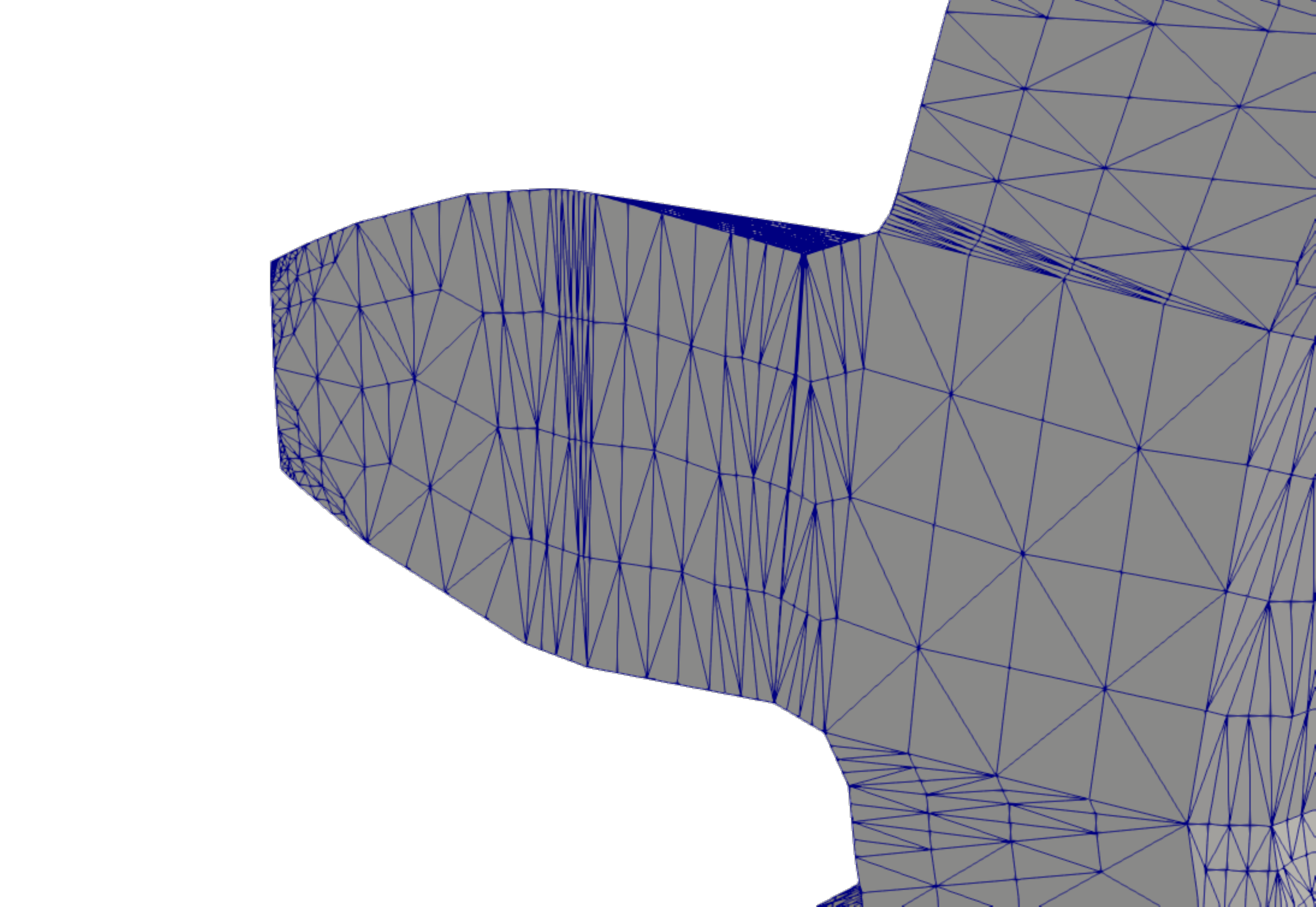}
  \caption{Left: geometry referent to the test case {\tt gear8302k}. Right: zoom in one
  of the gear teeth. A poor mesh quality is caused by elements close to the boundaries
  and those connecting the tooth to the main body of the gear.}
  \label{fig:gear8302k_mesh}
\end{figure}

The mesh derives from a 3D mechanical equilibrium of a helical gear with unitary outer
radius. Vertices belonging to the right-side boundary with $x > 0.9$ are rotated by
$10\degree$ while those belonging to the left-side boundary with $x < 0.1$ are
clamped. The mesh can be found in the examples folder of the DOLFIN library
\cite{LogWel10}. Mesh contains $380,280$ second order tetraedral elements and
$2,767,344$ vertices resulting in $8,302,032$ DOFs. Material is close to the
incompressible limit with $(E,\nu) = (10^{6} MPa,0.45)$. In Figure
\ref{fig:gear8302k_mesh}, we show a representation of the problem's geometry and mesh.

\subsection{Test case {\tt eiffel9213k}}
\label{sec:eiffel9213k}
\begin{figure}[!htbp]
  \centering
  \includegraphics[width=0.45\linewidth]{eiffel_1_noRS.pdf}
  \includegraphics[width=0.45\linewidth]{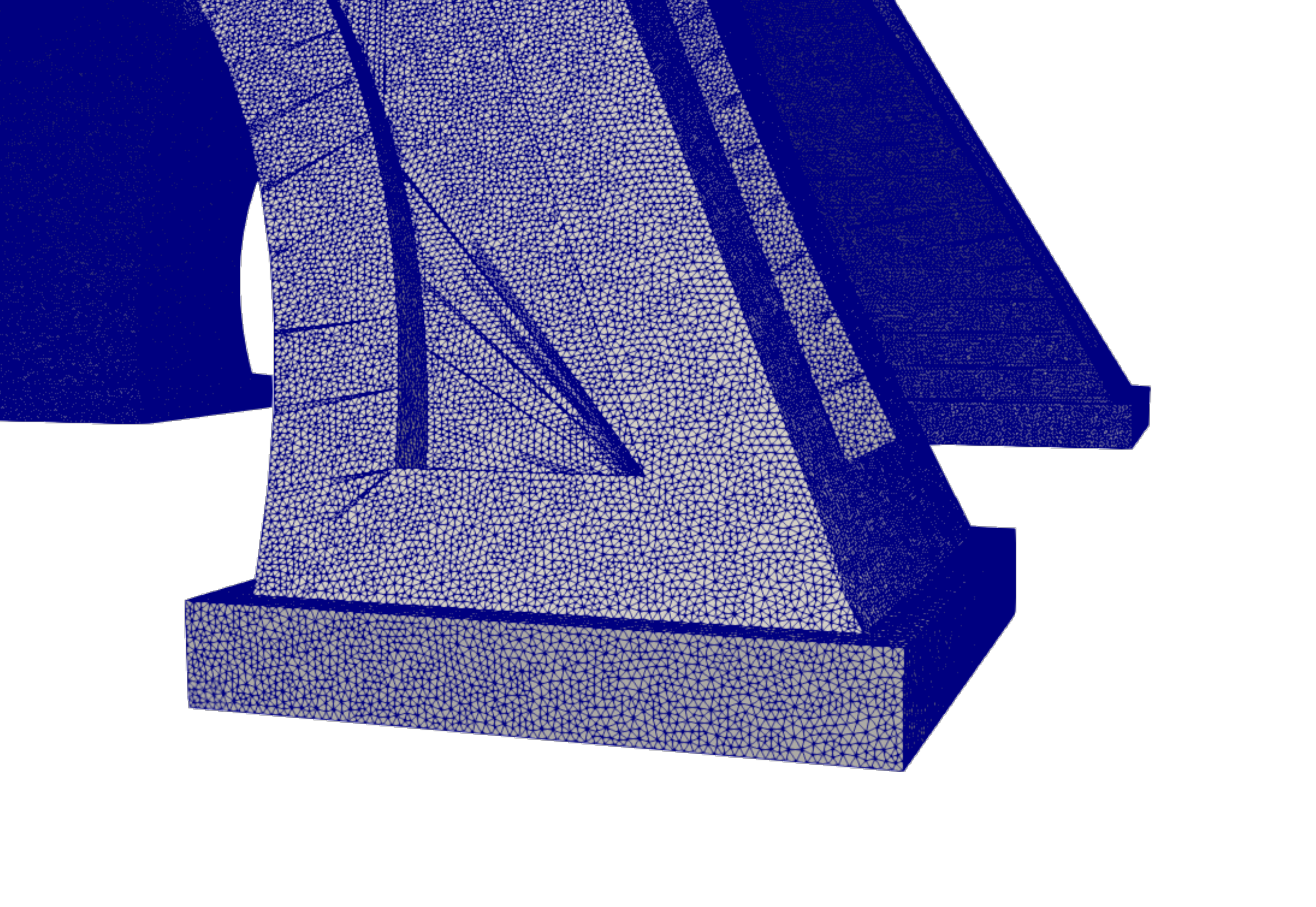}
  \caption{Left: geometry referent to the test case {\tt eiffel9213k} ({\bf T2}). Right:
  zoom in the lower right base of the tower indicating also the mesh used. Note that
  some elements close to the arc region are distorted.}
  \label{fig:eiffel9213k_mesh}
\end{figure}

The mesh derives from a mechanical equilibrium of a three-dimensional simplified model
of the Eiffel tower with clamped bases and two forces acting downwards: the first one
applied on top while the second is applied at the middle base of the tower. Material
is linear elastic with $(E, \nu) = (10^{6} MPa, 0.30)$. Mesh is formed by $16,584,015$
linear tetrahedra and discretization is performed by the finite element method. A
linear problem with $9,213,342$ DOFs is obtained and the system matrix contains
$390,108,294$ entries. In Figure \ref{fig:eiffel9213k_mesh}, we show a representation
of the problem's geometry and mesh.

\subsection{Test case {\tt wrench13m}}
\label{sec:wrench13m}
\begin{figure}[!htbp]
  \centering
  \includegraphics[width=0.45\linewidth]{wrench_1_noRS.pdf}
  \includegraphics[width=0.45\linewidth]{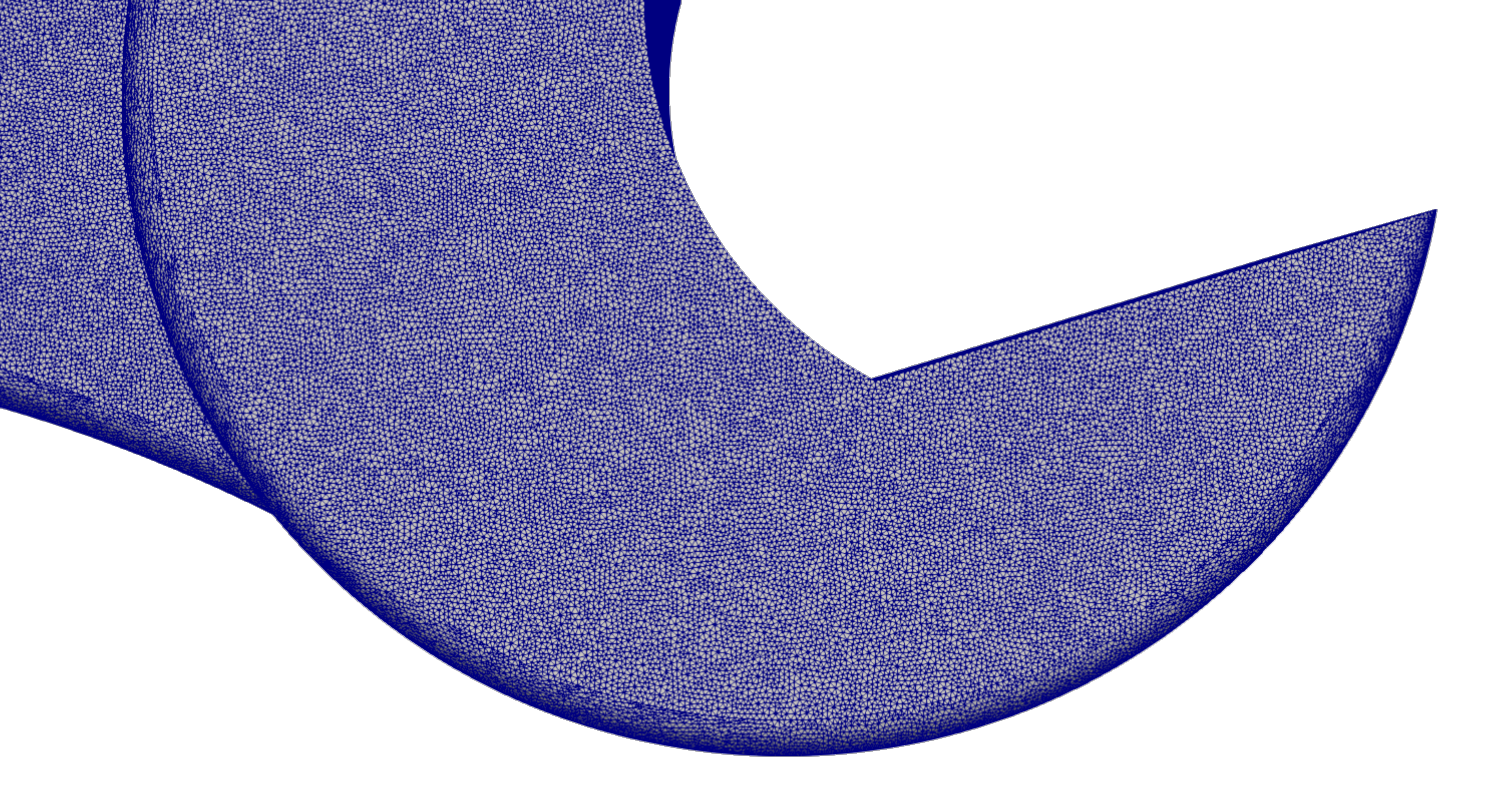}
  \caption{Left: geometry referent to the test case {\tt wrench13m}. Right: zoom in the
  right end of the wrench. Due to the simple format of this geometry, regular-shaped
  elements are able to fit it easily. The difficulty in this test case resides in the
  material properties, which are close to the incompressible limit.}
  \label{fig:wrench13m_mesh}
\end{figure}

The matrix arises from a 3D mechanical equilibrium of a wrench designed for $6/5$
bolts. One of the wrench's boundary facing the bolt is clamped while the remaining
boundaries are free. The unstructured mesh is composed by $25,692,045$ linear
tetrahedral elements and $4,331,788$ vertices. Material is elastic and close to the
incompressible limit with $(E,\nu) = (10^{6} MPa, 0.48)$. Figure \ref{fig:wrench13m_mesh}
shows the geometry and the mesh for this case.

\subsection{Test case {\tt agg14m}}
\label{sec:agg14m}
\begin{figure}[!htbp]
  \centering
  \includegraphics[width=0.45\linewidth]{agg_3_noRS.pdf}
  \includegraphics[width=0.45\linewidth]{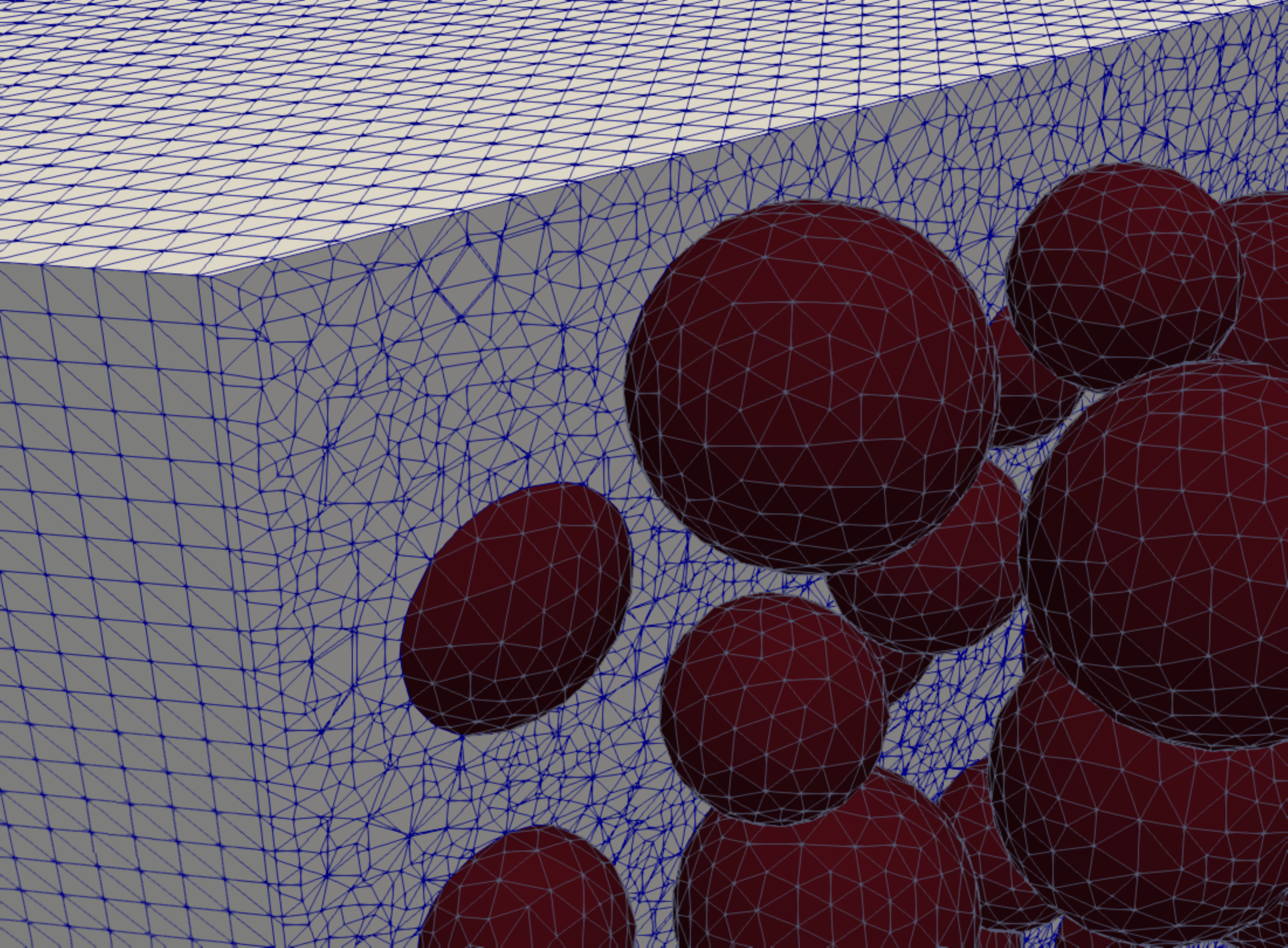}
  \caption{Left: geometry referent to the test case {\tt agg14m}, which employs a
  composite material. Right: zoom in the upper left corner of the mesh, showing the
  matrix and elipsoidal inclusions.}
  \label{fig:agg14m_mesh}
\end{figure}

The mesh derives from a 3D mesoscale simulation of a heterogeneous cube of lightened concrete. 
The domain has dimensions $50 \times 50 \times 50 \, mm^{3}$ and contains $2644$ ellipsoidal
inclusions of polystyrene. The cement matrix is characterized by $(E_1, \nu_1) = (25,000 MPa, 0.30)$,
while the polystyrene inclusions are characterized by $(E_2, \nu_2) = (5 MPa, 0.30)$. Hence, the contrast 
between the Young modules of these two materials is extremely high. 
The discretization is done via tetrahedral finite elements
\cite{MazPomSalMaj18}. In Figure \ref{fig:agg14m_mesh}, we show a representation of
the problem's geometry and mesh.

\subsection{Test case {\tt M20}}
\label{sec:M20}
The mesh derives from a 3D mechanical equilibrium of a symmetric machine cutter that is loosely
constrained. The unstructured mesh is composed by $4,577,974$ second order tetrahedra and
$6,713,144$ vertices resulting in $20,056,050$ DOFs. Material is linear elastic with
$(E,\nu) = (10^{8} MPa, 0.33)$. This problem was initially presented by \citet{KorLuGul14}
and later used in the work \cite{KorGup16}.

\section{Mesh qualities for the real-world engineering problems}
\label{sec:meshQuality}

\begin{figure}[!htbp]
  \centering
  \includegraphics[width=0.32\linewidth]{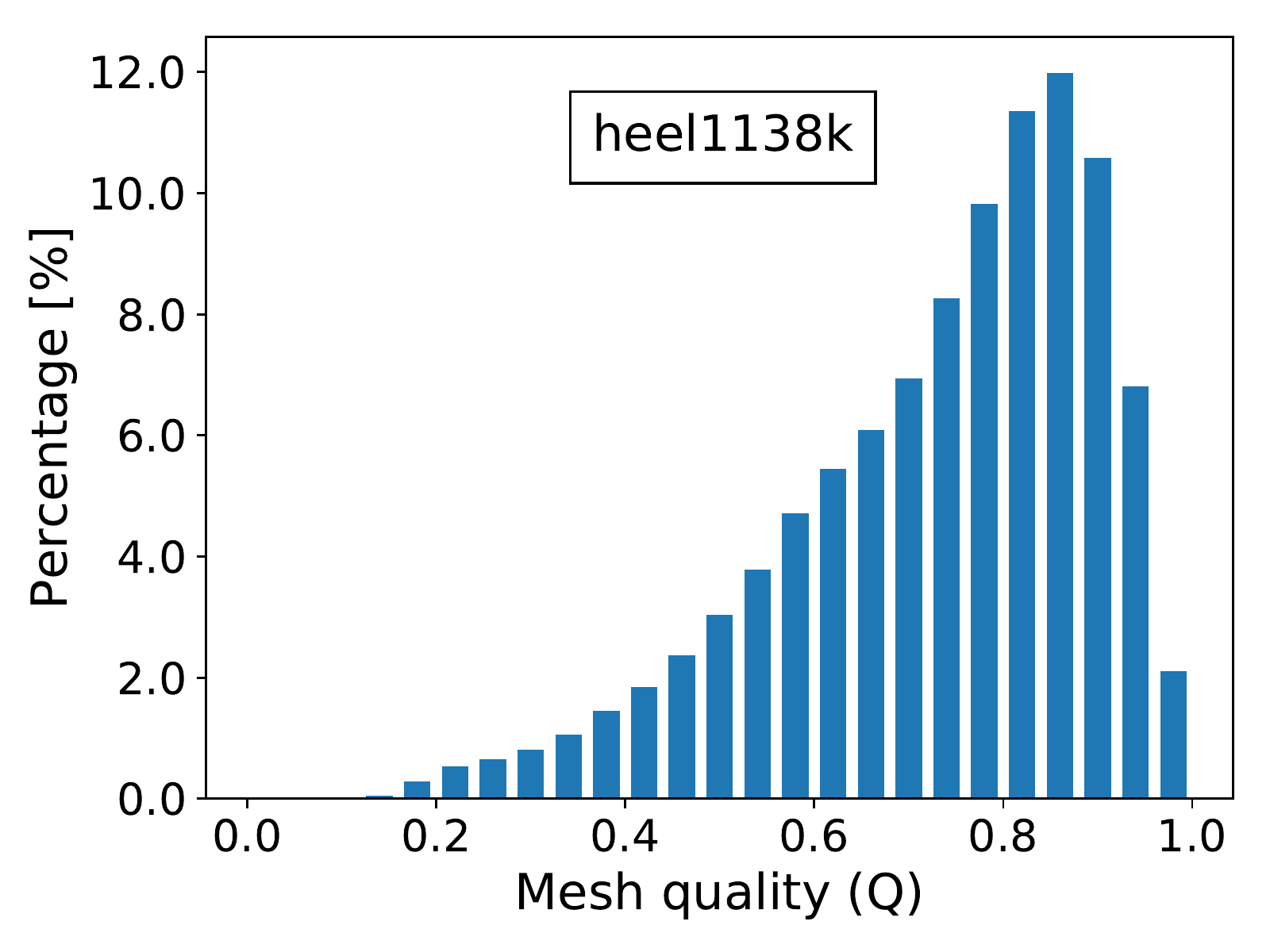}
  \includegraphics[width=0.32\linewidth]{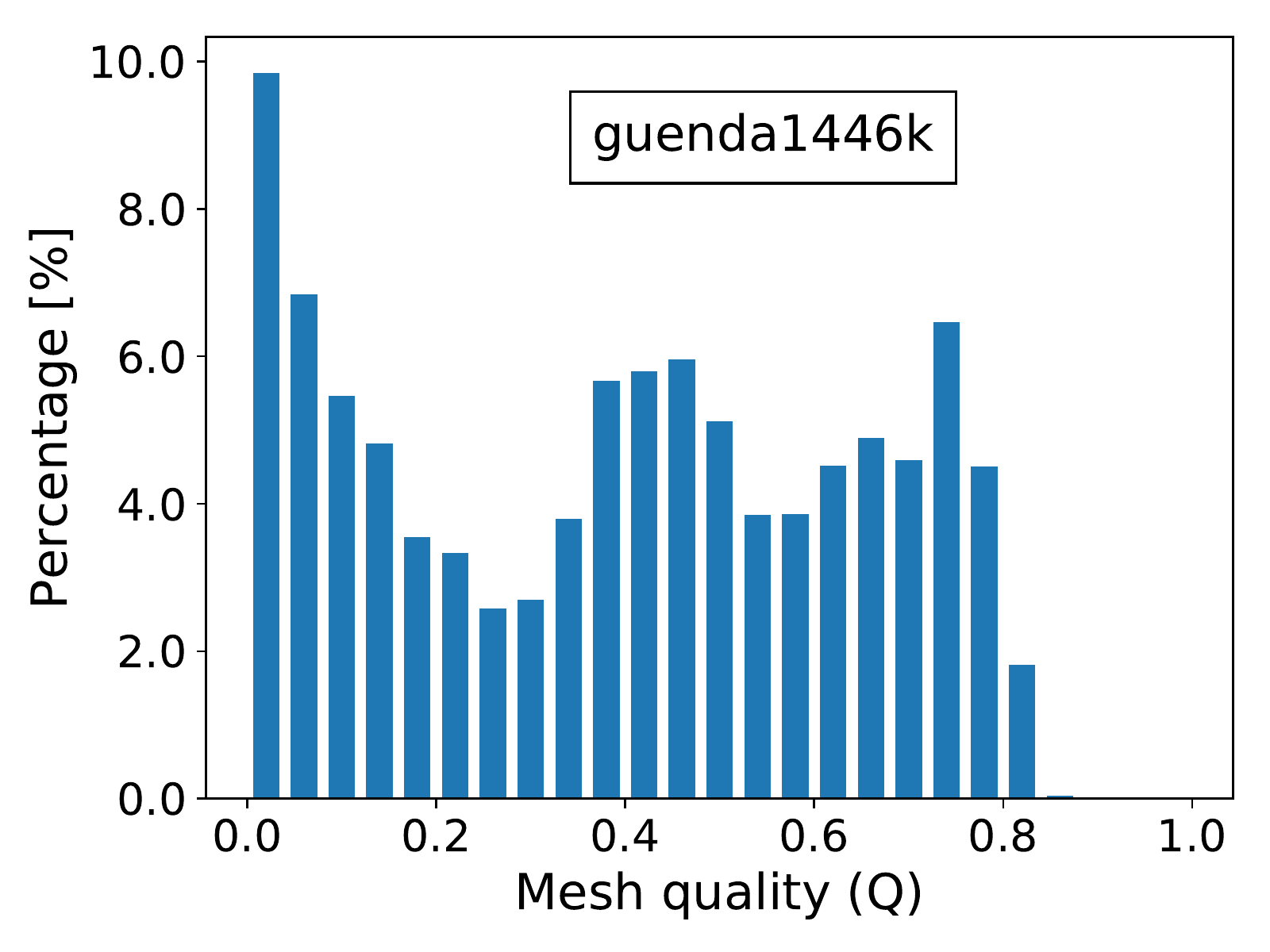}
  \includegraphics[width=0.32\linewidth]{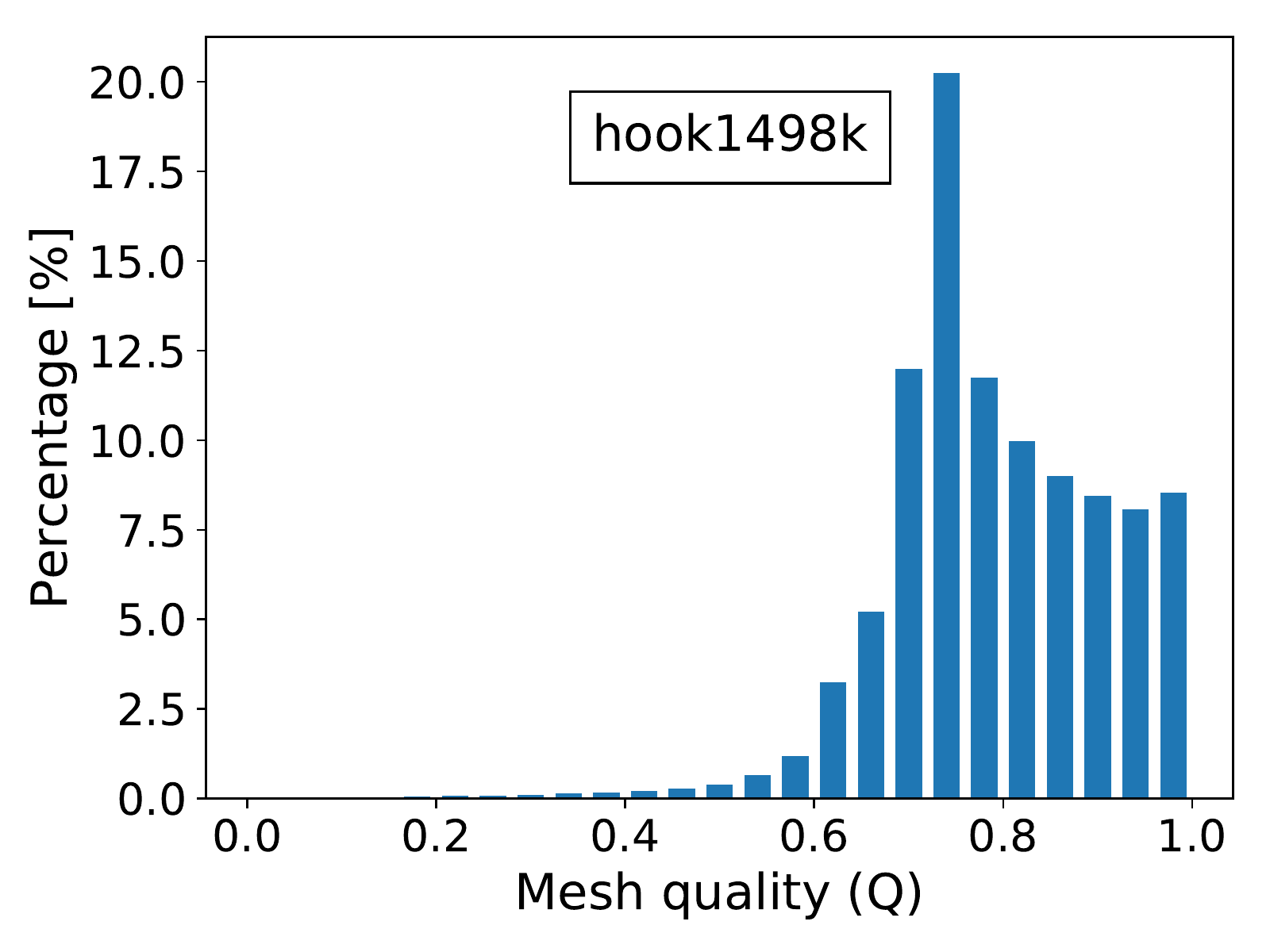}
  \includegraphics[width=0.32\linewidth]{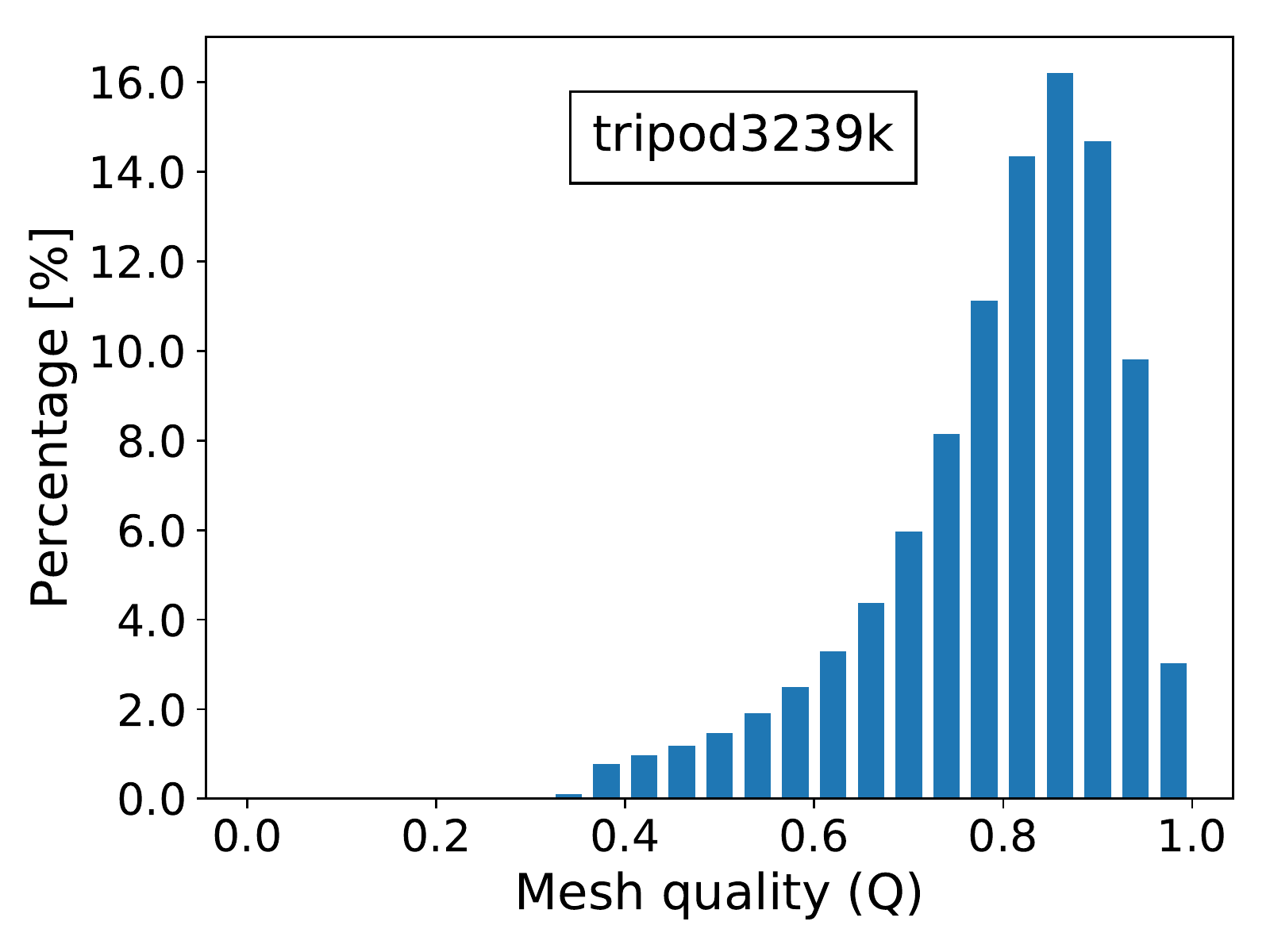}
  \includegraphics[width=0.32\linewidth]{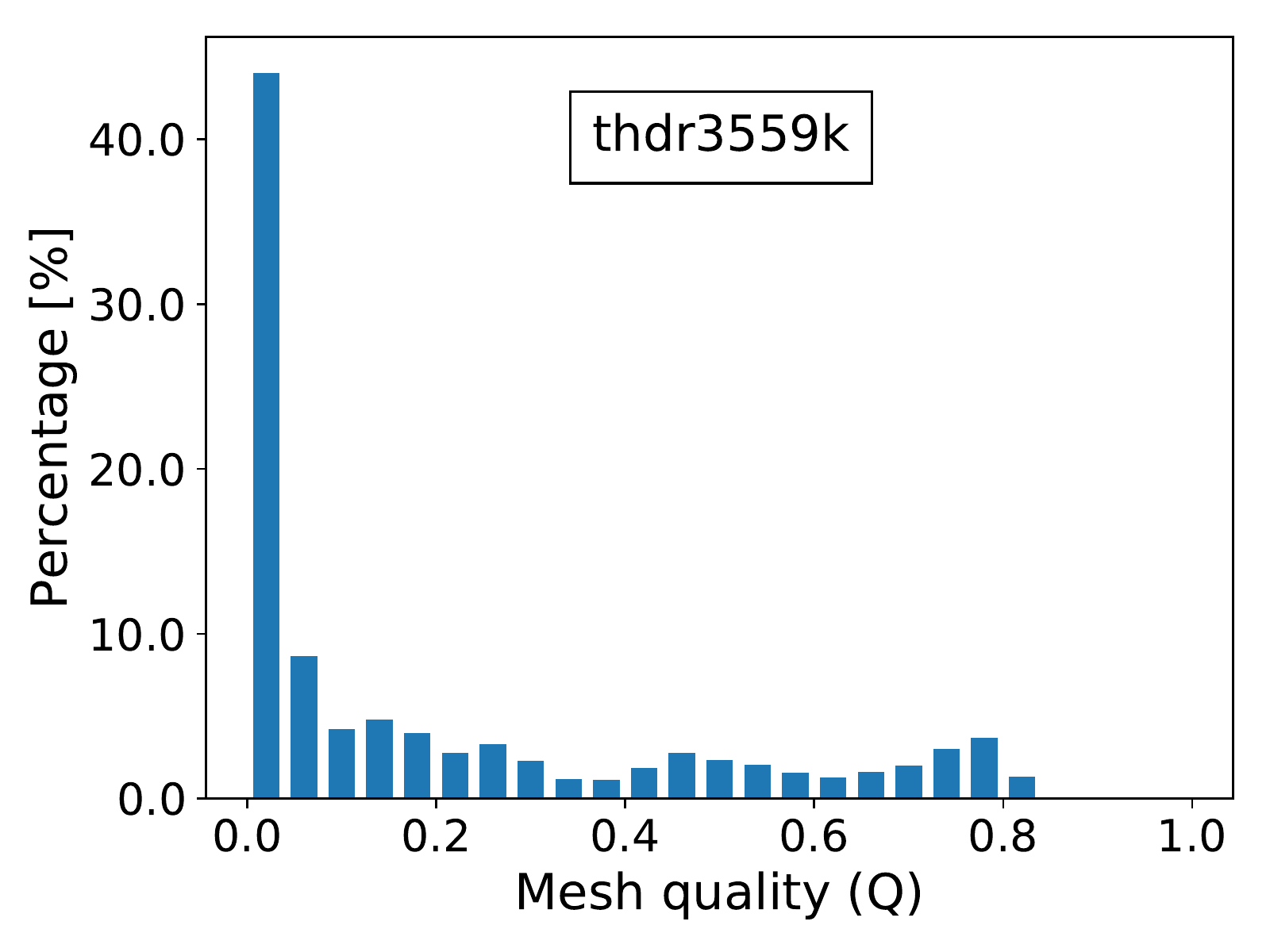}
  \includegraphics[width=0.32\linewidth]{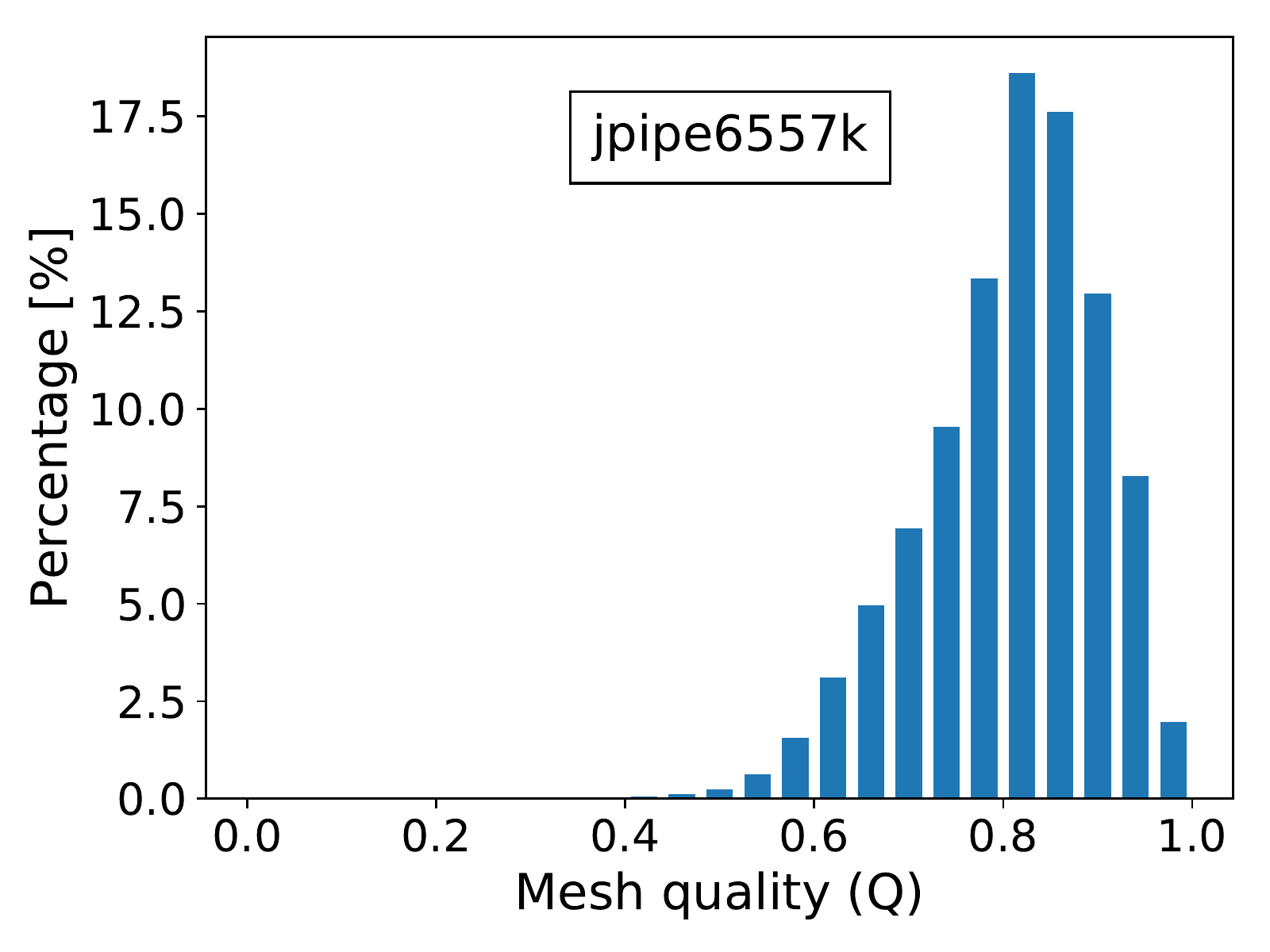}
  \includegraphics[width=0.32\linewidth]{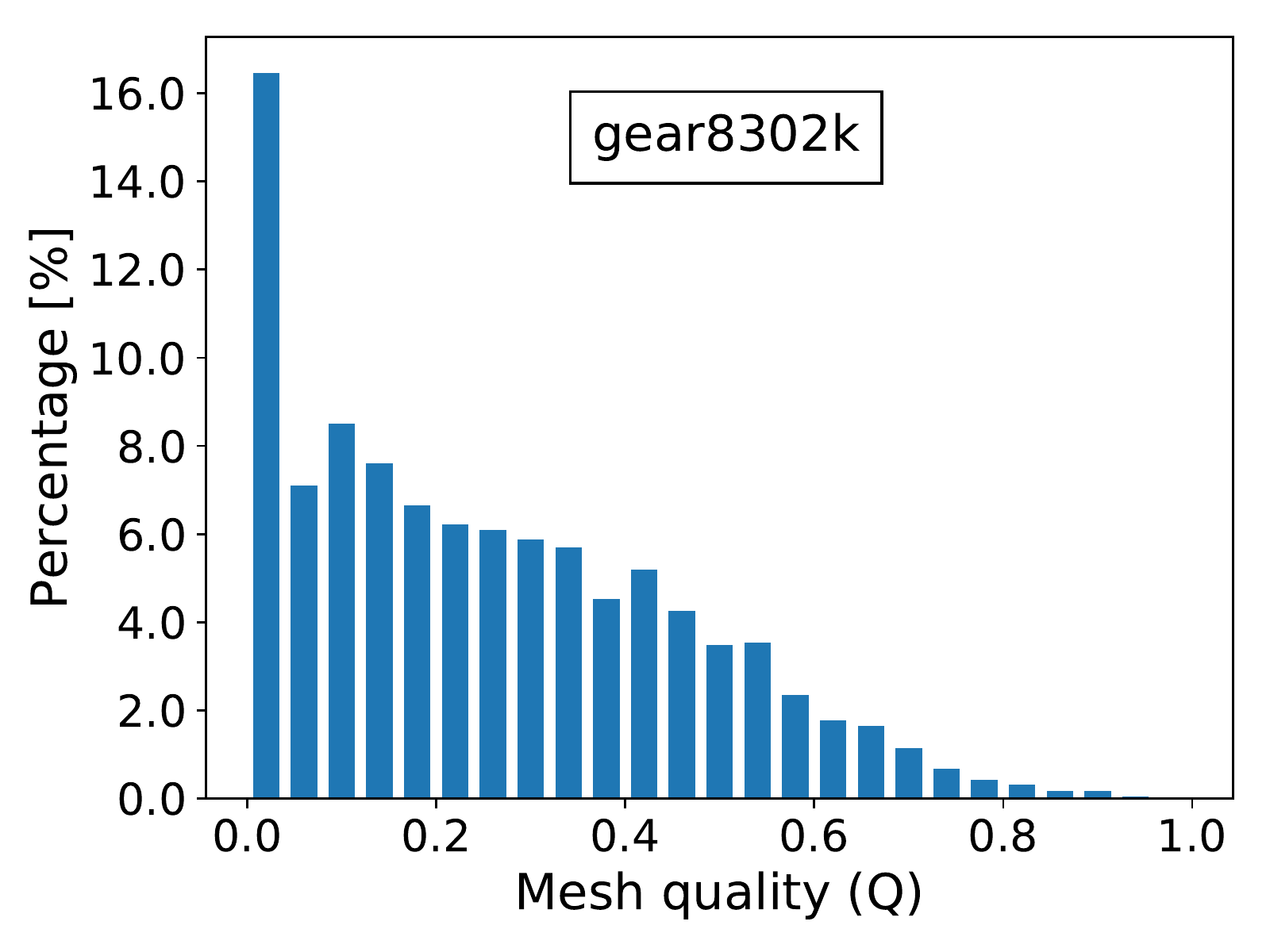}
  \includegraphics[width=0.32\linewidth]{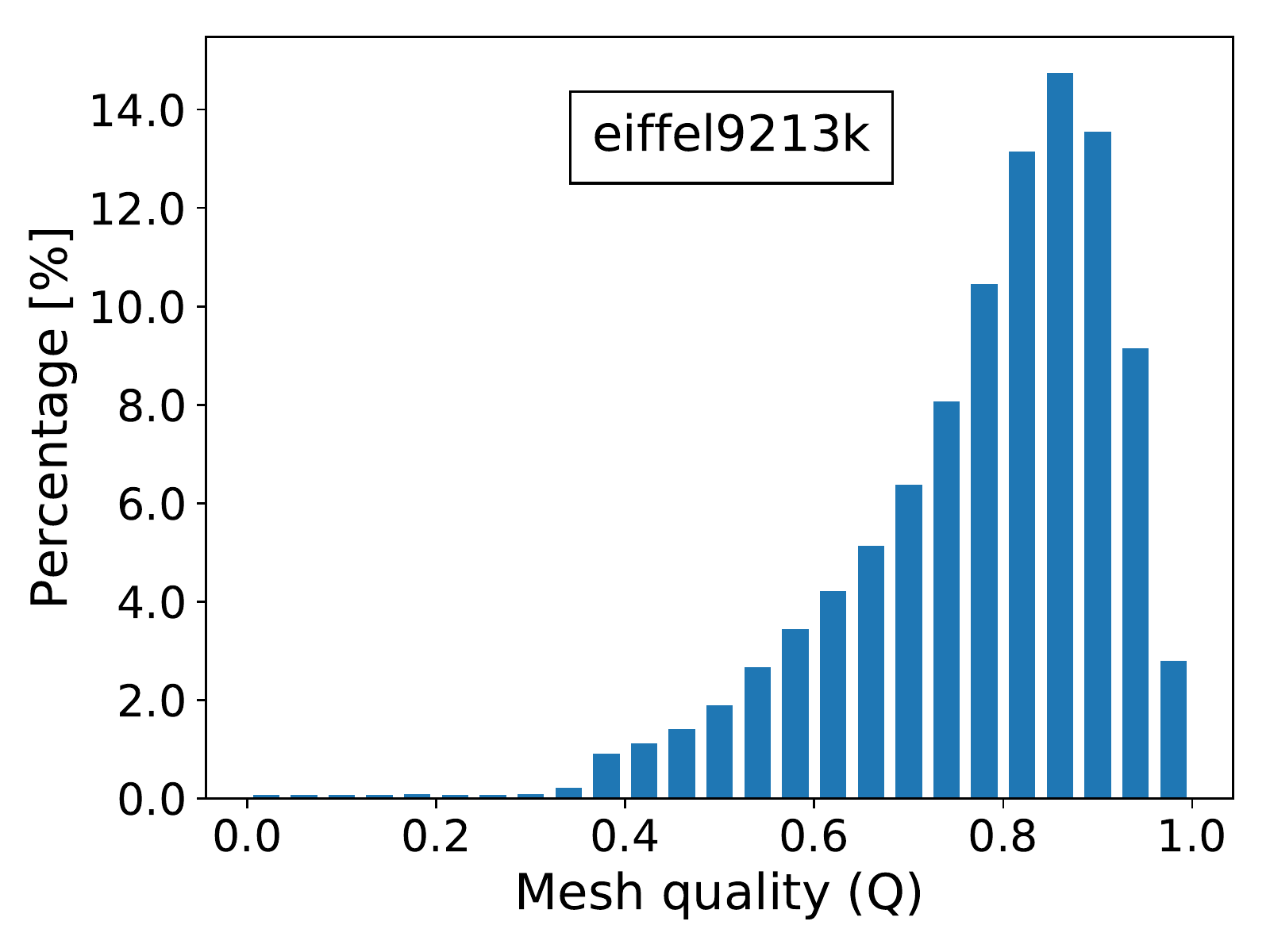}
  \includegraphics[width=0.32\linewidth]{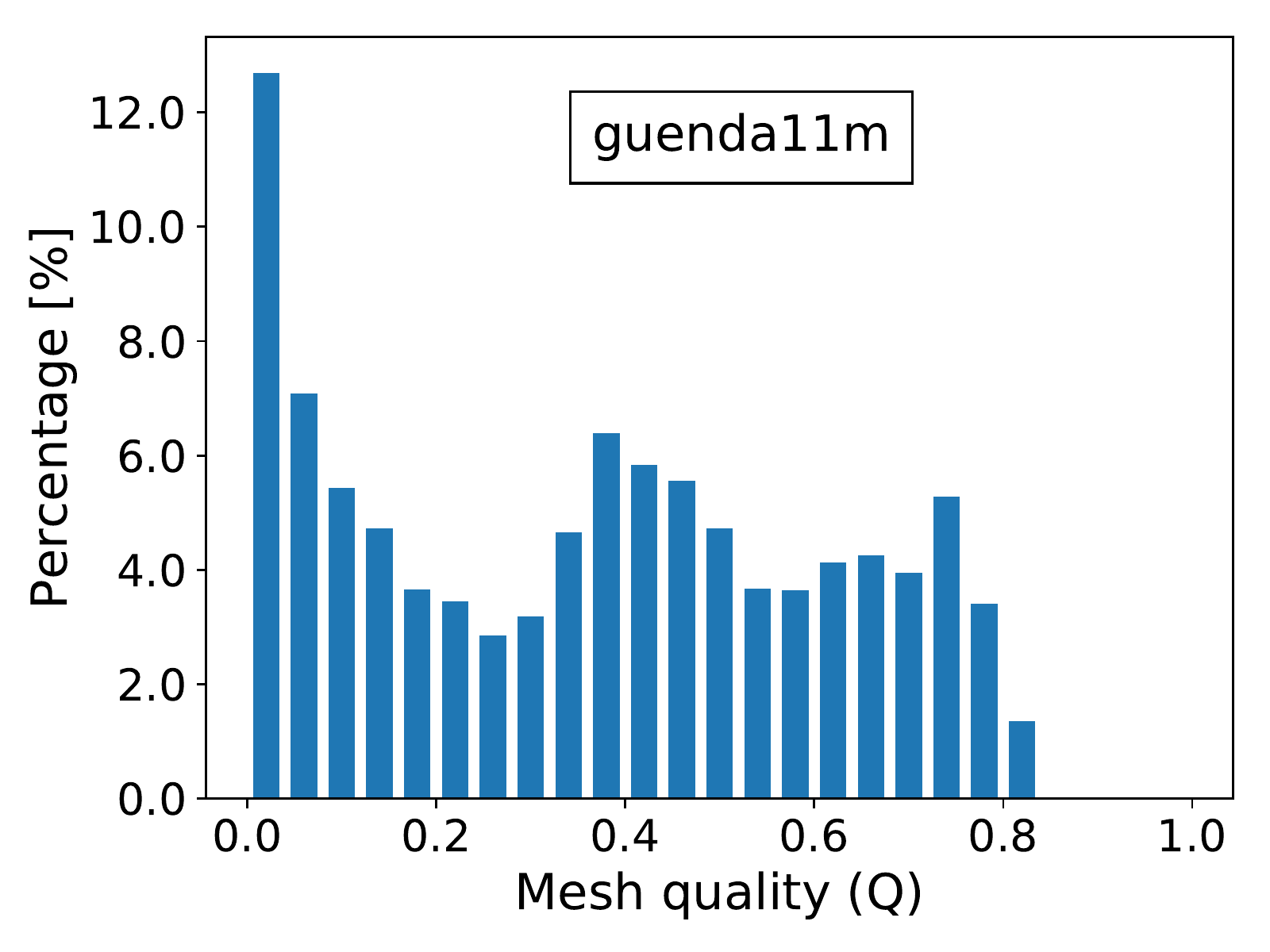}
  \includegraphics[width=0.32\linewidth]{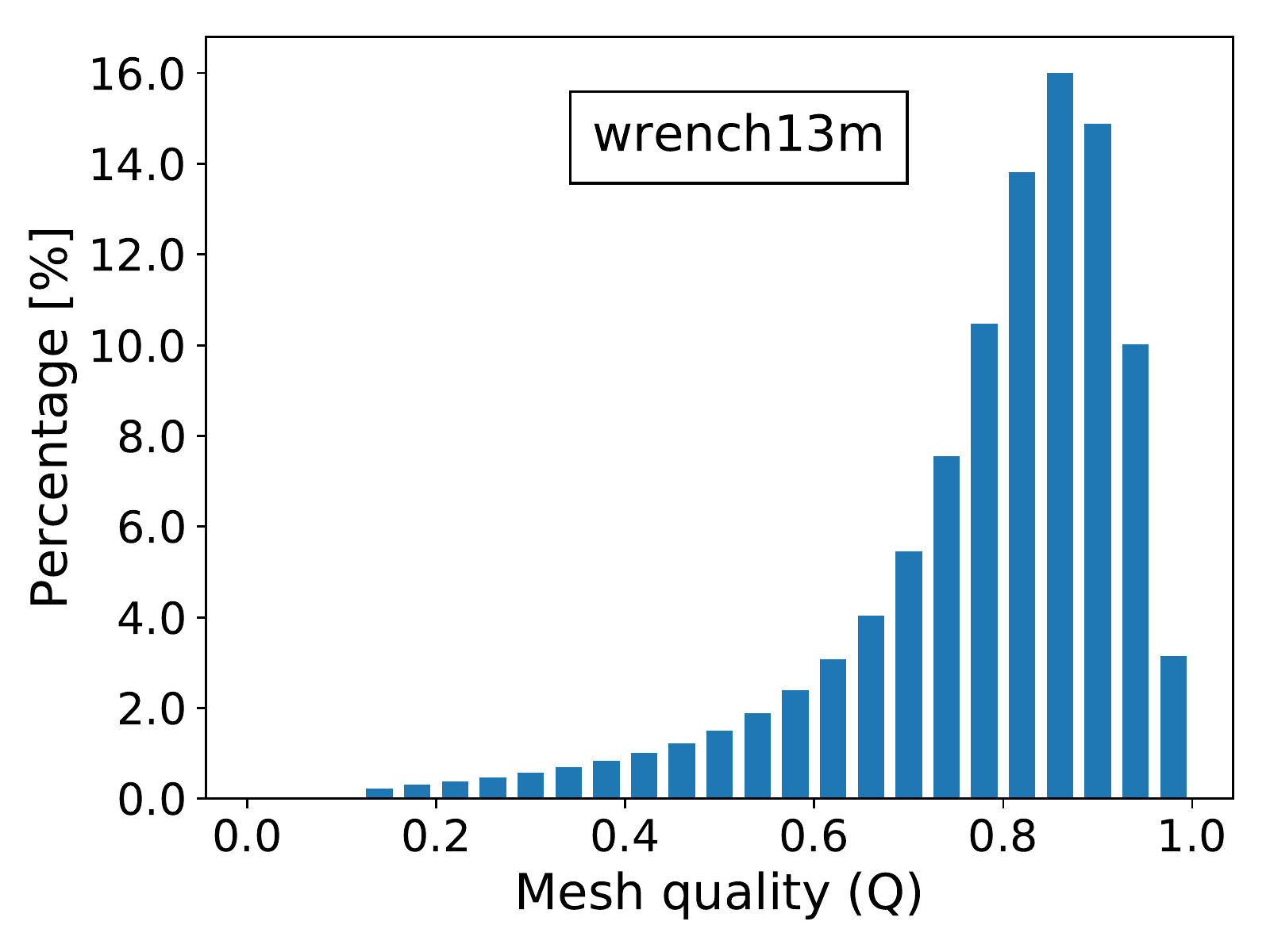}
  \includegraphics[width=0.32\linewidth]{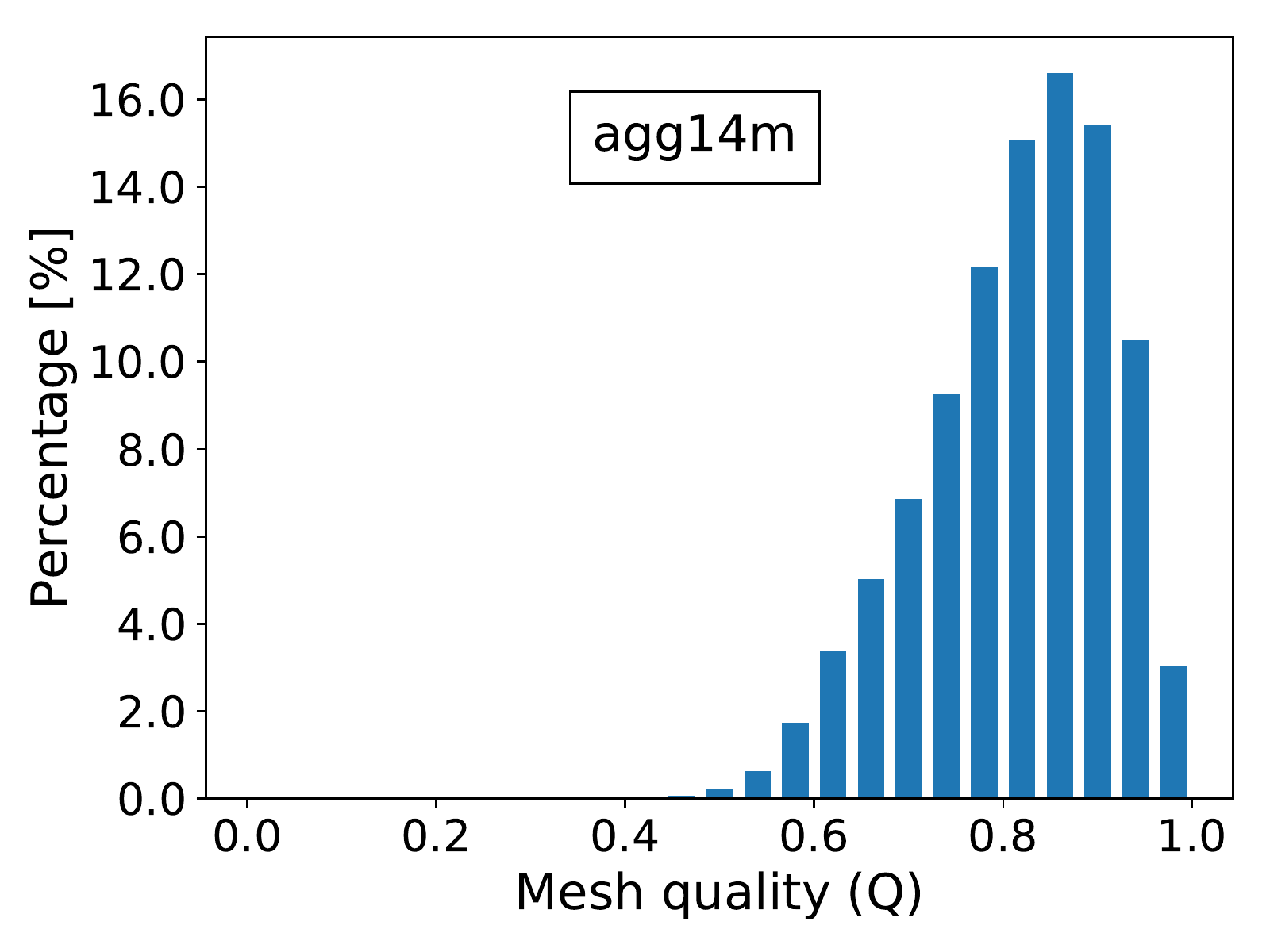}
  \includegraphics[width=0.32\linewidth]{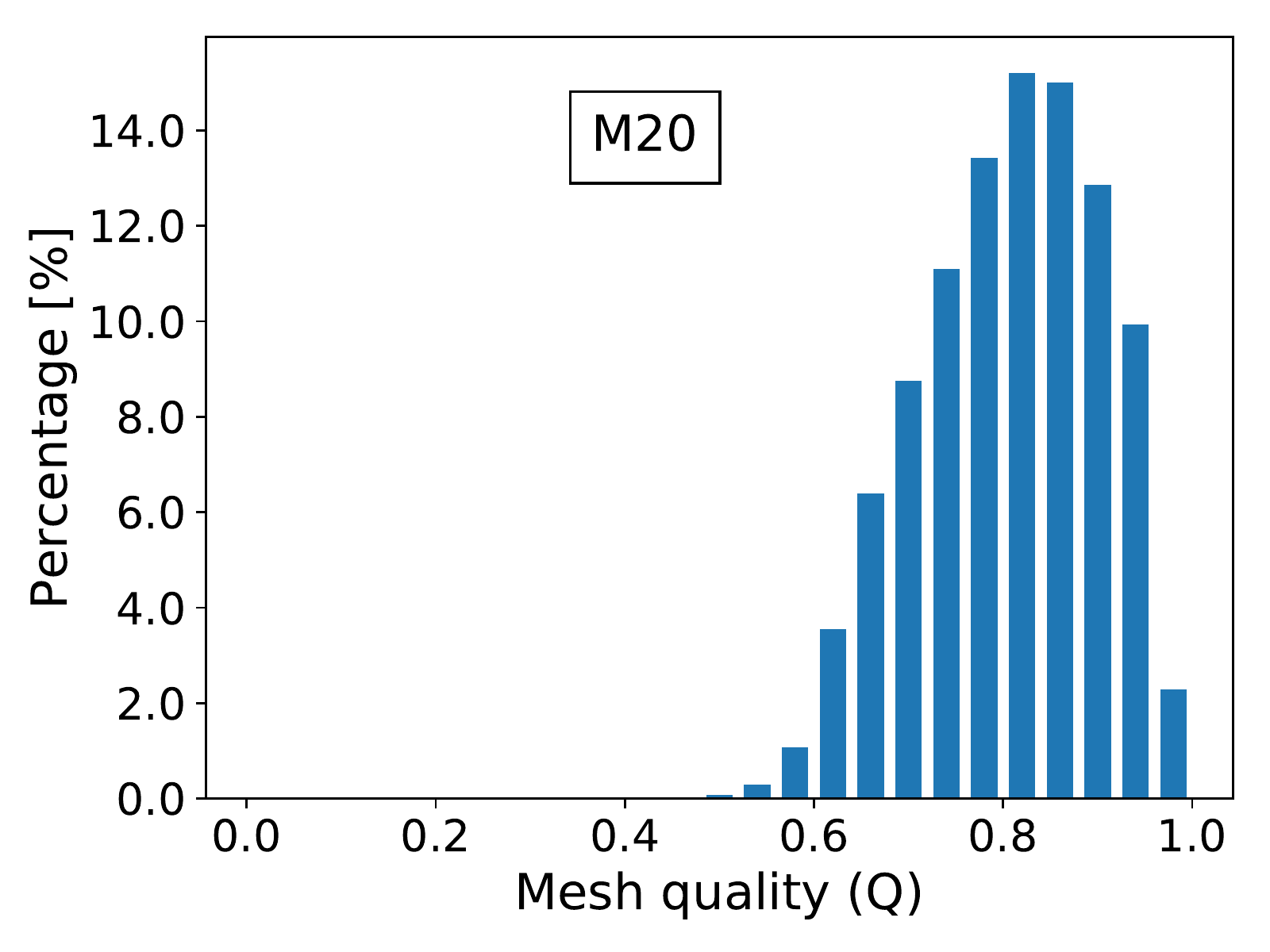}
  \caption{Percentage frequency distribution of the mesh quality indicator, $Q$, for the
  real-world engineering problems.}
  \label{fig:meshQuality}
\end{figure}

In this section, we provide an overview of the mesh qualities of the problems listed in
Table \ref{tab:realWorldProblems}. For this, we choose to quantify the mesh quality as in
the DOLFIN library \cite{LogWel10}. Given an element with index $i$, we compute the
measure
\begin{equation}
  Q(i) = d_p \cdot \dfrac{r(i)}{R(i)},
\label{eq:meshQuality}
\end{equation}
where $d_p$ is the topological dimension of the problem; $r$ is the radius of the biggest
circle/sphere that can be inscribed in the element and $R$, the radius of the smallest
circle/sphere circumscribing the element. Naturally, the elements' quality is better when
$Q(i)$ is close to unity and worst when it approximates zero, which is the limiting case
of degenerate elements. The last scenario contributes to the ill-conditioning of
discretized linear systems, causing difficulties for iterative solvers to converge. In
Figure \ref{fig:meshQuality}, we show the percentage frequency distribution of $Q(i)$ for
the real-world test cases with the only exception of case {\tt beam6502k}, since it is
composed by cubic hexahedral elements and, thus, has $Q(i) = 1$ everywhere.

\end{document}